\documentclass{article}
\usepackage[margin=1in]{geometry}
\usepackage{comment}
\usepackage{xcolor}

\newcommand\blfootnote[1]{%
  \begingroup
  \renewcommand\thefootnote{}\footnote{#1}%
  \addtocounter{footnote}{-1}%
  \endgroup
}

\usepackage{etoolbox}
\newtoggle{twocol}	
\togglefalse{twocol}
\newtoggle{proofs}	
\toggletrue{proofs}
\newtoggle{longform}	
\toggletrue{longform}
\newtoggle{thesis}	
\togglefalse{thesis}
\usepackage{mathtools,amsmath,amssymb,amsfonts}
\newcommand{\N}{\mathbb{N}}		
\newcommand{\R}{\mathbb{R}}		
\newcommand{\Rn}{\mathbb{R}^n}		
\newcommand{\Rpl}{\mathbb{R}_+}		
\newcommand{\B}{\mathcal{B}}		
\newcommand{\M}{\mathcal{M}}		
\renewcommand{\P}{\mathcal{P}}		
\newcommand{\D}{\mathcal{D}}		
\DeclareMathOperator*{\supp}{supp}	
\DeclareMathOperator*{\NE}{NE}		
\DeclareMathOperator{\lin}{Lin}		
\newcommand{\stab}{\mathfrak{S}}	

\newcommand{\TV}{\textup{TV}}		

\DeclareMathOperator*{\argmin}{arg\,min}	
\DeclareMathOperator*{\sign}{sign}		
\DeclareMathOperator*{\storage}{\Sigma}		

\usepackage{amsthm}
\theoremstyle{plain}
\newtheorem{theorem}{Theorem}
\newtheorem{proposition}{Proposition}
\newtheorem{corollary}{Corollary}
\newtheorem{lemma}{Lemma}
\theoremstyle{definition}
\newtheorem{definition}{Definition}
\newtheorem{example}{Example}
\newtheorem{assumption}{Assumption}
\theoremstyle{remark}
\newtheorem{remark}{Remark}

\usepackage[createShortEnv]{proof-at-the-end}

\usepackage[round]{natbib}

\usepackage{hyperref}			
\hypersetup{%
	colorlinks=true,%
	linkcolor=black,%
	urlcolor=black,%
	citecolor=black,%
	filecolor=black%
}					
\usepackage[all]{hypcap}		

\usepackage{graphicx,pgfplots,tikzscale}	
\pgfplotsset{compat=newest}
\pgfplotsset{plot coordinates/math parser=false}
\usetikzlibrary{external}
\usetikzlibrary{arrows,shapes,intersections}
\usepackage{tkz-euclide}	

\usepackage{booktabs}

\newcommand\hl[1]{{#1}}	

\begin{document}
	\title{Evolutionary Games on Infinite Strategy Sets: \\ Convergence to Nash Equilibria via Dissipativity}
	\author{Brendon G.\ Anderson\footnote{Brendon G.\ Anderson is with the Mechanical Engineering Department, California Polytechnic State University (email: \href{mailto:bga@calpoly.edu}{\texttt{bga@calpoly.edu}}).} \and Jingqi Li\footnote{Jingqi Li, Somayeh Sojoudi, and Murat Arcak are with the Department of Electrical Engineering and Computer Sciences, University of California, Berkeley (emails: \texttt{\string{\href{mailto:jingqili@berkeley.edu}{\texttt{jingqili}},\href{mailto:sojoudi@berkeley.edu}{\texttt{sojoudi}},\href{mailto:arcak@berkeley.edu}{\texttt{arcak}}\string}@berkeley.edu}).} \and Somayeh Sojoudi\footnotemark[2] \and Murat Arcak\footnotemark[2]}
	\date{}
	\maketitle
	\blfootnote{This material is based upon work supported in part by the U.\ S.\ Army Research Laboratory and the U.\ S.\ Army Research Office under grant number W911NF2010219. It was also supported by ONR and NSF, including NSF grant CNS-2135791.}

	\begin{abstract}
We consider evolutionary dynamics for population games in which players have a continuum of strategies at their disposal. Models in this setting amount to infinite-dimensional differential equations evolving on the manifold of probability measures. \iftoggle{thesis}{In this thesis, we}{We} generalize dissipativity theory for evolutionary games from finite to infinite strategy sets that are compact metric spaces, and derive sufficient conditions for the stability of Nash equilibria under the infinite-dimensional dynamics. The resulting analysis is applicable to a broad class of evolutionary games, and is modular in the sense that the pertinent conditions on the dynamics and the game's payoff structure can be verified independently. By specializing our theory to the class of monotone games, we recover as special cases existing stability results for the Brown-von Neumann-Nash and impartial pairwise comparison dynamics.\iftoggle{longform}{ We also extend our theory to models with dynamic payoffs, further broadening the applicability of our framework.}{} {\color{black}Throughout our analyses, we identify and elaborate on new technical conditions that are key in extending dissipativity theory from finite to infinite strategy sets, such as compactness of the set of Nash equilibria\iftoggle{longform}{ and evolution of dynamic payoffs within a compact positively invariant set}{}.} We illustrate our theory using a variety of \hl{case studies}, including a novel, continuous variant of the war of attrition game.
\end{abstract}

	\tableofcontents

		\section{Introduction}
	\label{sec: introduction}

	Population games are models in which a large number of agents \hl{interact strategically}. Examples of such models appear ubiquitously in engineering and societal-scale problems, including traffic congestion networks, decentralized control, and economic markets \citep{sandholm2010population}. Within a population game, each agent employs a strategy available to them to maximize their expected payoff. When the agents are permitted to continuously revise their strategy according to some protocol, the game gives rise to an \emph{evolutionary dynamics model} \citep{smith1982evolutionary,sandholm2010population}. Such models have a rich history within the mathematical biology literature, as reviewed in \citet{hofbauer1998evolutionary}.

	Traditional game-theoretic models are concerned with notions of (Nash) equilibrium states, in which no player is incentivized to choose a different strategy given knowledge of the payoffs. However, such notions of equilibria are static and incomplete, in the sense that they do not capture whether an evolutionary dynamics model dynamically converges to them when players revise their strategies according to some protocol. Indeed, static equilibria, such as Nash equilibria, need not be dynamically stable \citep{sato2002chaos,hart2003uncoupled,sandholm2010population}. This has led to an entire body of works concerned with assessing the\iftoggle{thesis}{ dynamic}{} stability of evolutionary games. Although many works have proven stability for specific examples of evolutionary dynamics, \iftoggle{thesis}{researchers maintain the overarching goal of proving}{it is important to prove} stability for the most general classes of models possible \citep{fox2013population}, {\color{black}as doing so would provide a principled foundation for assessing the behavior of new and emerging classes of evolutionary dynamics.}
    
	The aim of this \iftoggle{thesis}{thesis}{paper} is to prove \iftoggle{thesis}{dynamic }{}stability for a very broad class of evolutionary dynamics. The broadening of evolutionary stability theory has seen two notable directions: 1) generalizing the structural behavior of the dynamics and the game as much as possible while maintaining stability, and 2) generalizing prior stability results for specific dynamical structures to more abstract settings. We now discuss these two approaches in further depth.

	\subsection{Related Works}

	\subsubsection{Potentiality, Monotonicity, and Dissipativity}

	Potential games constitute a class of games in which the payoff is given by the gradient of a ``potential function'' \citep{monderer1996potential,sandholm2001potential}. It was shown in \citet{sandholm2001potential} that potential games satisfying the so-called ``positive correlation property'' admit the potential function as a global Lyapunov function, thereby yielding\iftoggle{thesis}{ dynamic}{} stability guarantees. \citet{hofbauer2007stable,hofbauer2009stable} introduce monotone games, also known as stable games,\iftoggle{thesis}{\footnote{We use the terminology ``monotone'' throughout this \iftoggle{thesis}{thesis}{paper}, as it more accurately represents these games' defining property \eqref{eq: monotone_game}, to come later, whereas the terminology ``stable'' may lead to confusion when discussing\iftoggle{thesis}{ dynamic}{} stability; a stable game need not give rise to stable dynamics.}}{\footnote{We use the terminology ``monotone'' throughout this \iftoggle{thesis}{thesis}{paper}, as it more accurately represents these games' defining property \eqref{eq: monotone_game}, to come later, than the alternatives ``stable,'' ``contractive,'' and ``negative semidefinite.''}} which generalize potential games with concave potential functions to allow for payoffs that act as monotone operators (like gradients of concave functions). Zero-sum games and games with an interior globally evolutionarily stable state are also known to be special cases of monotone games \citep{hofbauer2007stable}. \citet{hofbauer2007stable} show that common evolutionary dynamics, such as those of Brown-von Neumann-Nash and Smith, exhibit\iftoggle{thesis}{ dynamic}{} stability when coupled with monotone games. \citet{hofbauer2009stable} extend this result to the case of general dynamics satisfying an ``integrability'' condition on their revision protocols.

	Generalizing even further past integrability and monotonicity, \citet{fox2013population} apply notions of passivity from the systems and control literature to grant stability. The authors propose to view an evolutionary game as a nonautonomous dynamical system in feedback with inputs defined by the game's payoffs. In doing so, they prove that ``$\delta$-passive'' evolutionary dynamics coupled with monotone games yield\iftoggle{thesis}{ dynamic}{} stability. The core intuition is that, if the rate of change of internally stored energy of an evolutionary system is less than the rate of energy supplied to it by the game's payoffs, then the closed-loop system's total energy decreases. This approach was taken by \citet{mabrok2021passivity} to analyze the stability of replicator dynamics, and was further generalized by \citet{arcak2021dissipativity} to apply to more general ``$\delta$-dissipative'' dynamics. The dissipativity theory of \citet{arcak2021dissipativity} constitutes some of the broadest stability results, recovering many of the aforementioned prior results as special cases. \hl{Recent works have applied these broad theories to particular applications, such as distributed Nash equilibrium seeking \citep{martinez2023distributed} and the analysis of strategy-dependent pairwise comparison revision protocols \citep{kara2023pairwise}.} We emphasize that all of the works mentioned here are restricted to games defined over a finite number of strategies.

	\subsubsection{Games over Infinite Strategy Sets}

	Many practical games come equipped with an infinite number of strategies available to the players, e.g., pricing and generation in power systems \citep{park2001continuous}, games of timing such as the war of attrition \citep{bishop1978generalized}, plant growth models in biology \citep[Section~2.4]{bomze1989game}, \hl{and, more recently, in optimization \citep{anderson2024approximately} and multi-agent reinforcement learning with continuous action spaces \citep{mazumdar2020policy}}. Consequently, much effort has gone into abstracting\iftoggle{thesis}{ dynamic}{} stability results from the finite-strategy setting into the infinite setting. However, in doing so, the distribution of strategies being employed, termed the ``population state,'' becomes a probability measure rather than a finite-dimensional vector in a Euclidean simplex \citep{myerson1991game}. This makes the analysis much more challenging, as it requires studying differential equations evolving on the manifold of probability measures within an infinite-dimensional Banach space.

	Alongside this technical hurdle come two other key challenges. First, there are multiple standard notions of convergence for probability measures on infinite sets, and evolutionary dynamics may converge in one such notion but fail to converge in another. For example, \citet{eshel2003evolutionary} show that the\iftoggle{thesis}{ infamous}{} replicator dynamics may exhibit dynamic instability with respect to the so-called ``maximum shift topology'' even when they are stable in the weak topology.\footnote{Formal definitions of the weak and strong topologies are given in \autoref{sec: math_preliminaries}.} The second key challenge is that\iftoggle{thesis}{ dynamic}{} stability may break when moving from the finite-strategy regime to the infinite regime. For example, \citet{oechssler2002dynamic} show that even strict Nash equilibria and evolutionarily stable states may be unstable under the replicator dynamics over infinite strategy sets, even though their approximations with finitely many strategies are always\iftoggle{thesis}{ dynamically}{} stable. Similarly, we show in \autoref{sec: examples} that finite-strategy approximations of the war of attrition are guaranteed to be\iftoggle{thesis}{ dynamically}{} stable via the finite-dimensional dissipativity theory of \citet{arcak2021dissipativity}, despite the fact that the underlying infinite-dimensional game is unstable. Thus, stability guarantees for evolutionary dynamics over infinite strategy sets are not automatic from their corresponding finite-dimensional counterparts. This motivates our work in directly considering the infinite strategy set setting.

	A handful of related works have directly analyzed the stability of infinite-dimensional evolutionary dynamics. Some of the first such work was \citet{bomze1990dynamical,bomze1991cross}, which considered replicator dynamics with respect to the strong topology. A line of follow-up works on the replicator dynamics has emerged, many of which come to the consensus that convergence in the sense of the weak topology is most appropriate for evolutionary dynamics, as it better respects notions of distance between strategies \citep{oechssler2001evolutionary,oechssler2002dynamic,cressman2005stability,cressman2005measure,hingu2018evolutionary,hingu2020superiority}. These works also propose alternative notions of equilibria (beyond Nash) to ensure\iftoggle{thesis}{ dynamic}{} stability. Beyond the replicator dynamics, stability (typically of Nash equilibria) with respect to the weak topology has been assessed for the Brown-von Neumann-Nash, pairwise comparison, logit, general imitative, and perturbed best response dynamics \citep{hofbauer2009brown,cheung2014pairwise,lahkar2015logit,cheung2016imitative, lahkar2022generalized}. However, despite the applicability of these results to quite general strategy sets, all of these works are restricted to specific evolutionary dynamics and are proven in a case-by-case fashion. In comparison, the approach in this \iftoggle{thesis}{thesis}{paper} is to keep in the spirit of broadening stability guarantees, and to derive results for infinite-strategy games applicable to general classes of dissipative dynamics.

	\subsection{Contributions}
	In this \iftoggle{thesis}{thesis}{paper}, we unify the two above generalization approaches to achieve the following primary contributions:
	\begin{enumerate}
		\item We introduce novel notions of dissipativity for evolutionary dynamics over infinite strategy sets. This extension from the finite to infinite strategy sets requires new technicalities, since our models are defined on Banach spaces, and the weak topology in which we seek dynamic convergence is not equivalent to the topology induced by the total variation norm on our population states, unlike the finite-dimensional Hilbert space setting in which they are equivalent.
		\item In our main result (\autoref{thm: delta-dissipative}), \hl{we prove a new stability theorem showing that} $\delta$-dissipative evolutionary dynamics on infinite strategy sets weakly converge to Nash equilibria \hl{under decreasing energy supply rates induced by the game's payoffs} and some technical regularity conditions. Our complete identification of such technical conditions is nontrivial, as, again, our infinite-dimensional setting breaks down the topological equivalence between notions of convergence and notions of norm.
		\item \hl{We specialize our framework to prove a new stability theorem for the class of monotone games} (our \autoref{thm: monotone-passive}), and prove that this specialization recovers the main stability results of \citet[Theorem~3]{hofbauer2009brown} and \citet[Theorem~4]{cheung2014pairwise} as special cases (our \autoref{cor: recovering_monotone}).
		\iftoggle{longform}{\item We further extend the generality of our stability theory to the case in which the game's payoffs exhibit dynamic behavior (\autoref{thm: dpedm}).}{}
		\item \iftoggle{thesis}{\hl{We show that the classical war of attrition game on infinite strategy sets simultaneously fails to converge to Nash equilibria while its finite-strategy approximations succeed in convergence to Nash, and we use our theoretical framework to identify the technical stability conditions being violated.} We subsequently propose and verify the stability of a new ``continuous'' variant of the war of attrition game, and we illustrate the generality of our framework on an example of a monotone game with smoothed payoff dynamics.}{We construct a practical infinite-strategy example that fails to converge to a Nash equilibrium, even though its finite-dimensional approximations are proven to converge through prior dissipativity tools. We use our theory to identify the technical stability conditions being violated, and subsequently propose and verify the stability of a new, continuous variant of this game.}
		{\color{black}\item A variety of case studies and simulations are provided to illustrate the generality of our infinite-strategy stability theory.}
	\end{enumerate}

	\subsection{Outline}

	This \iftoggle{thesis}{thesis}{paper} is organized as follows. In \autoref{sec: math_preliminaries}, we introduce our notations and review relevant mathematical definitions and results. \iftoggle{thesis}{Population games and evolutionary dynamics are formally introduced in \autoref{sec: population_games} and \autoref{sec: evolutionary_dynamics}, along with associated definitions and results for both static and dynamic stability.}{Population games and evolutionary dynamics, along with their equilibrium states and stability, are formally introduced and discussed in \autoref{sec: population_games} and \autoref{sec: evolutionary_dynamics}.} \hl{Our primary contributions are given in \autoref{sec: dissipativity_theory} and \autoref{sec: examples}.} Namely, in \autoref{sec: dissipativity_theory}, we present our dissipativity theory for infinite strategy sets and our stability theorems. In \autoref{sec: examples}, we provide \hl{case studies} illustrating our framework and results. We \iftoggle{thesis}{conclude}{give conclusions} in \autoref{sec: conclusions}. The proofs of \iftoggle{thesis}{our}{our key} \iftoggle{thesis}{primary results (the theorems)}{theorems} are given in the main text\iftoggle{thesis}{. \iftoggle{longform}{To streamline presentation}{Due to space constraints}, all other proofs (the propositions and corollaries) are deferred to \iftoggle{longform}{\autoref{sec: proofs}}{our online technical report \citet{anderson2025dissipativity-long}}}{; \iftoggle{longform}{to streamline presentation}{due to space constraints}, proofs of the propositions and corollaries are deferred to \iftoggle{longform}{\autoref{sec: proofs}}{our online technical report \citet{anderson2025dissipativity-long}}}. \iftoggle{longform}{Supplementary definitions, results, and discussions are given in \autoref{sec: supplementary}.}{}

		\section{Mathematical Preliminaries}
	\label{sec: math_preliminaries}

	\subsection{Notations and Basic Definitions}

	The set of nonnegative real numbers is denoted by $\Rpl$. We define $\sign \colon \R \to \R$ by $\sign(x) = 1$ for $x>0$, $\sign(x)=0$ for $x=0$, and $\sign(x)=-1$ for $x<0$. The dual space of a normed vector space $X$ (i.e., the space of bounded linear functionals on $X$) is denoted by $X^*$. Let $S$ be a compact metric space. The Banach space of bounded continuous real-valued functions on $S$ endowed with the supremum norm is denoted by $(C_b(S),\|\cdot\|_\infty)$. Since $S$ is compact, $C_b(S)$ equals the set of all continuous real-valued functions on $S$, denoted by $C(S)$. The Borel $\sigma$-algebra on $S$ is denoted by $\B(S)$, and the Banach space of finite signed Borel measures on $S$ endowed with the total variation norm is denoted by $(\M(S),\|\cdot\|_\TV)$. Recall that $\|\mu\|_\TV \coloneqq |\mu|(S) = \sup_{\text{$f$ measurable}:\|f\|_\infty \le 1} \int_S f d\mu$, where $|\mu|$ is the total variation measure of $\mu$. The support of a measure $\mu\in\M(S)$ is denoted by $\supp(\mu)$.

	We denote the set of probability measures on $(S,\B(S))$ by $\P(S) = \{\mu\in\M_+(S) : \mu(X) = 1\}$, where $\M_+(S)\subseteq\M(S)$ is the set of positive Borel measures on $S$. \hl{The tangent space of $\P(S)$ is given by $T\P(S) = \{\nu\in\M(S) : \nu(S) = 0\}$, which is a linear subspace of $\M(S)$.} The Dirac measure at $s\in S$ is denoted by $\delta_s\in\P(S)$. We define the bilinear form $\left<\cdot,\cdot\right>\colon C(S) \times \M(S) \to \R$ by $\left<f,\mu\right> = \int_S f d\mu$, which is well-defined and satisfies $\left| \int_S f d\mu \right| \le \|f\|_\infty \|\mu\|_\TV$ for all $f\in C(S)$ and all $\mu\in\M(S)$. Recall that $\M(S)$ is isometrically isomorphic to the dual space of $C(S)$ \citep[Theorem~7.17]{folland1999real}, and therefore every element of $\M(S)$ can be uniquely identified with a bounded linear functional on $C(S)$. Thus, for all bounded linear functionals $I \in C(S)^*$, there exists a unique $\mu\in\M(S)$ such that $I(f) = \left<f,\mu\right>$ for all $f\in C(S)$.

    {\color{black}Our most commonly used notations---both introduced here, and in what follows---are summarized in \autoref{tab: symbols}.}

    \begin{table}[ht]
        {\color{black}
	\caption{List of commonly used symbols.}
        \begin{center}
	\begin{tabular}{l l}
	\toprule
	Symbol & Meaning \\
	\midrule %
	$S$ & Strategy set (compact metric space) \\
        $\P(S)$ & Set of probability measures on $S$ \\
        $\M(S)$ & Set of finite signed Borel measures on $S$ \\
        $C(S)$ & Set of continuous real-valued functions on $S$ \\
        $\left<f,\nu\right>$ & Integral of function $f$ against measure $\nu$, i.e., $\int_S f d\nu$ \\
        $\mu$ & Population state (probability measure) \\
        $F$ & Population game \\
        $F_\mu(s)$ & Average payoff to strategy $s$ at population state $\mu$ \\
        $E_F(\nu,\mu)$ & Average payoff to population $\nu$ relative to population $\mu$ \\
        $\NE(F)$ & Set of Nash equilibria of game $F$ \\
        $\rho$ & Payoff (continuous function) \\
        $v$ & Dynamics map (governs evolution of population $\mu$) \\
        \iftoggle{longform}{$u$ & Payoff map (governs evolution of payoff $\rho$) \\}{}
	\bottomrule
	\end{tabular}
        \end{center}
	\label{tab: symbols}
        }
    \end{table}

	\subsection{Topologies and Convergence of Measures}
	\label{sec: topologies_convergence_measures}

	Two types of convergence in $\M(S)$ will be of use. Recall that a sequence $\{\mu_n \in\M(S) : n\in\N\}$ \emph{converges weakly to $\mu\in\M(S)$} if $\lim_{n\to \infty} \int_S f d\mu_n = \int_S f d\mu$ for all $f\in C(S)$, and \emph{converges strongly to $\mu\in\M(S)$} if $\lim_{n\to\infty}\|\mu_n - \mu\|_\TV = 0$.
%
%
	Recall that strong convergence implies weak convergence. Strong and weak convergence induce topologies on $\M(S)$, termed the \emph{strong topology} and \emph{weak topology}, respectively.\footnote{The weak topology is sometimes called the ``narrow topology.'' Since $S$ is compact, the weak topology coincides with the weak-$*$ topology on $\M(S) = C(S)^*$ (i.e., the weakest topology on $C(S)^*$ making every element $f\in C(S) \subseteq C(S)^{**}$ a continuous linear functional on $C(S)^*$). In functional analysis the term ``weak topology'' on $\M(S)$ would refer to the weakest topology on $\M(S)$ making every element of the dual space $\M(S)^* = C(S)^{**}$ continuous. We stick with our definitions to remain consistent with related works.} \iftoggle{thesis}{We will also need the following product topology.

	\begin{definition}
		\label{def: weak-infinity_topology}
		Consider $\M(S)$ endowed with the weak topology and $C(S)$ endowed with its usual topology induced by $\|\cdot\|_\infty$. We call the corresponding product topology on $\M(S)\times C(S)$ the \emph{weak-$\infty$ topology}.
	\end{definition}}{We call the product topology on $\M(S)\times C(S)$ induced by the weak topology on $\M(S)$ and the norm topology on $C(S)$ the \emph{weak-$\infty$ topology}.}

	\hl{We use the following fact throughout our analyses.}

	\begin{lemma}[{\hl{\citealp[Theorem~6.4]{parthasarathy1967probability}}}]
		\label{lem: P(S)_compact}
		\iftoggle{thesis}{It holds that }{}$\P(S)$ is weakly compact.
	\end{lemma}

	\subsection{Notions of Differentiability}
	\label{sec: differentiability}

	\hl{We also need various notions of differentiability. Consider Banach spaces $(X,\|\cdot\|_X)$, $(Y,\|\cdot\|_Y)$, and $(Z,\|\cdot\|_Z)$, and open sets $U\subseteq X$ and $V\subseteq Y$.\iftoggle{thesis}{}{ The space of bounded linear operators from $X$ to $Y$ is denoted by $\lin(X,Y)$.} The Fr\'echet derivative of a map $f\colon U\to Y$ at $x\in U$, if it exists, is denoted by $Df(x)$. Recall that \iftoggle{thesis}{$Df(x) \colon X \to Y$ is a bounded linear operator}{$Df(x) \in \lin(X,Y)$}. If $f \colon U\times V \to Z$ is a map defined on $U\times V$, its first partial Fr\'echet derivative at $(x,y)\in U\times V$, if it exists, is given by $\partial_1 f(x,y) \coloneqq D(f(\cdot,y))(x)$. The second partial Fr\'echet derivative of such a map $f$ is similarly given by $\partial_2 f(x,y) \coloneqq D(f(x,\cdot))(y)$. Recall that a map $x\colon [0,\infty) \to X$ is said to be \emph{differentiable at $t\in [0,\infty)$} if there exists $\dot{x}(t)\in X$ such that
		\begin{equation*}
			\lim_{\epsilon \to 0}\left\lVert \frac{x(t+\epsilon) - x(t)}{\epsilon} - \dot{x}(t) \right\rVert_X = 0,
		\end{equation*}
		and in this case $\dot{x}(t)$ is called the \emph{derivative of $x$ at $t$}. We call a map $\mu \colon [0,\infty) \to \M(S)$ \emph{strongly differentiable at $t$} if it is differentiable at $t$, to emphasize the underlying topology on $\M(S)$ induced by $\|\cdot\|_\TV$. \iftoggle{longform}{In \autoref{sec: supplementary}, these notions of differentiability are defined more formally and are discussed further.}{}}

		\section{Population Games}
	\label{sec: population_games}

	We now describe the game-theoretic aspects of our problem. The compact set $S$ represents the (infinite) set of pure strategies of the game, and is hence called the \emph{strategy set}.\footnote{\hl{The compactness of the strategy set $S$ is standard in the literature on infinite-dimensional evolutionary games. Although this compactness is a technical condition needed for our use of Lyapunov theory, it is also an important qualitative requirement in our context of games, as it ensures that evolutionary dynamics move the population state towards distributions of strategies that are actually available to the players. For example, compactness avoids cases where there exists a ``hidden Nash equilibrium'' at a probability measure with support at strategies on the boundary of $S$ or ``at infinity'' that are inaccessible by the players. {\color{black}See \autoref{ex: at_infinity} for an example of such an equilibrium ``at infinity.''}}} A \emph{population state} is a distribution $\mu\in\P(S)$, which encodes how strategies in $S$ are being employed across the game's population. Thus, $\P(S)$ is termed the \emph{population state space}. To every population state $\mu\in\P(S)$ associates a \emph{mean payoff function} $F_\mu \in C(S)$ such that $F_\mu(s)$ quantifies the average payoff to strategy $s$ when the population is at state $\mu$. We refer to the mapping $F \colon \P(S) \to C(S)$ defined by $F(\mu) = F_\mu$ as the \emph{population game}, or simply the \emph{game}. One of the primary quantities of interest when analyzing population games over infinite strategy sets is the \emph{average mean payoff} $E_F(\nu,\mu)\in\R$ to a population state $\nu\in\P(S)$ relative to a population state $\mu\in\P(S)$, which is given by
	\begin{equation*}
		E_F(\nu,\mu) \coloneqq \left<F(\mu),\nu\right> = \int_S F_\mu d\nu.
	\end{equation*}
	The average mean payoff gives rise to a simple definition for Nash equilibria of population games.
	\begin{definition}
		\label{def: ne}
		A population state $\mu \in\P(S)$ is a \emph{Nash equilibrium of the game $F \colon \P(S) \to C(S)$} if
		\begin{equation}
			E_F(\nu,\mu) \le E_F(\mu,\mu)
			\label{eq: nash}
		\end{equation}
		for all $\nu\in\P(S)$. If, additionally, the inequality \eqref{eq: nash} holds strictly for all $\nu\in\P(S)\setminus\{\mu\}$, then $\mu$ is a \emph{strict Nash equilibrium of the game $F$}. The set of all Nash equilibria of the game $F$ is denoted by $\NE(F)$.
	\end{definition}
    
    Intuitively, a population state $\mu\in\P(S)$ is a Nash equilibrium if the average mean payoff to the population cannot be increased by moving to any other state $\nu\in\P(S)$ given the current payoffs defined by $F_\mu$. It is important to note that the notion of a Nash equilibrium is static in the sense that it does not depend on any dynamical behavior endowed to the game.\iftoggle{longform}{\hl{ Other types of relevant static equilibria are discussed in \autoref{sec: supplementary}.}}{} {\color{black} To illustrate the notions of population states, mean payoffs, and Nash equilibria in population games, we now discuss the classical war of attrition game.
    
    \begin{example}
    \label{ex: war_of_attrition}
    \hl{The ``war of attrition'' game is motivated by animal conflict and studied at length in \citet[Chapter~3]{smith1982evolution}.} We adopt the formalism of the game from \citet{bishop1978generalized} and \citet[Example~6]{hofbauer2009brown}. Consider a contest being carried out on a time interval $S \coloneqq [0,T] \subseteq \R$, with a common value of $V>0$ awarded to the winner. The winner is the one who decides to compete in the contest for the longest amount of time. In this case, a population state of the game is a probability distribution on $[0,T]$, encoding the proportion of players deciding to remain in the game for any given length of time.
    
    Formally, the game's mean payoffs are given by $F_\mu(s) = \int_S f(s,s') d\mu(s')$, where
\begin{equation*}
	f(s,s') = \begin{aligned}
		\begin{cases}
			V - s' & \text{if $s' < s$}, \\
			\frac{V}{2} - s & \text{if $s'=s$}, \\
			-s & \text{if $s'>s$},
		\end{cases}
	\end{aligned}
\end{equation*}
defines the payoff to a player employing strategy $s$ when their opponent employs strategy $s'$. It is assumed that $T > V / 2$, so that there may be incentive to resigning from the contest before time $T$. Thus, when a player's opponent quits the game earlier, at time $s'<s$, the player wins the value $V$, minus the ``damage'' $s'$ incurred by the length of the contest. Similarly, when a player quits the game earlier than their opponent, at time $s<s'$, then they incur the ``damage'' $-s$, without gaining any portion of the value $V$. Finally, if a player and their opponent quit at the same time $s=s'$, then they each gain half of the value, $\frac{V}{2}$, and incur ``damage'' $s=s'$.

This game $F$ has a unique Nash equilibrium $\mu^\star \in \P(S)$ with cumulative distribution given by
\begin{equation}
	\mu^\star([0,s]) = \begin{aligned}
		\begin{cases}
			1 - e^{-s / V} & \text{if $s \in [0,s^\star)$}, \\
			1 - e^{- s^\star / V} & \text{if $s\in [s^\star, T)$}, \\
			1 & \text{if $s = T$},
		\end{cases}
	\end{aligned}
	\label{eq: war_nash}
\end{equation}
where $s^\star = T - V / 2$ (cf., \citealt{bishop1978generalized,hofbauer2009brown}). Since this measure has a point mass at $s=T$, it does not have a density function with respect to the Lebesgue measure. Despite this, it is useful to think of the region $[s^\star,T)$, on which the above cumulative distribution is constant, as having zero density. With this in mind, the Nash equilibrium $\mu^\star$ intuitively corresponds to the case where players randomly select times $s\in[0,T]$ to play the game for, which tend to either be in the more conservative regime with $s\in [0,s^\star)$, or be ``all in'' with $s = T$. On the other hand, the population is not at equilibrium when players put probability mass on times $s\in (s^\star,T)$, since those players may either decrease the time that they play in order to reduce their incurred damage on average, or increase their play time to $s=T$ in order to increase the value they gain on average.
\end{example}
}

    The following result gives equivalent characterizations of Nash equilibria, which are used throughout our proofs. \hl{Such characterizations are sometimes taken as alternative definitions in the literature, e.g., in \citet{hofbauer2009brown,cheung2014pairwise}, albeit without proof of equivalence.}

	\iftoggle{longform}{
	\begin{propositionE}[][end,restate,text link=]
	}{
	\begin{propositionE}[][end,text link=]
	}
		\label{prop: equivalent_nash}
		Consider a game $F\colon \P(S) \to C(S)$, and let $\mu\in \P(S)$. The following are equivalent:
		\begin{enumerate}
			\item $\mu$ is a Nash equilibrium of the game $F$.
			\item $E_F(\delta_s,\mu) \le E_F(\mu,\mu)$ for all $s\in S$.
			\item $F_\mu(s) \le F_\mu(s')$ for all $s\in S$ and all $s' \in \supp(\mu)$.
		\end{enumerate}
	\end{propositionE}

	\iftoggle{longform}{
	\begin{proofE}
		Suppose that the third condition holds, so that $F_\mu(s) \le F_\mu(s')$ for all $s\in S$ and all $s' \in \supp(\mu)$. Then, for all $s\in S$, it holds that $E_F(\delta_s,\mu) = F_\mu(s) \le F_\mu(s')$ for all $s'\in \supp(\mu)$ and consequently that $E_F(\delta_s,\mu) = \int_S F_\mu(s) d\mu(s') \le \int_S F_\mu(s') d\mu(s') = E_F(\mu,\mu)$. Thus, the second condition holds. Furthermore, if $\nu\in\P(S)$, then $E_F(\nu,\mu) = \int_S F_\mu(s) d\nu(s) = \int_S E_F(\delta_s,\mu) d\nu(s) \le \int_S E_F(\mu,\mu) d\nu(s) = E_F(\mu,\mu)$, so the first condition holds as well.

		To complete the proof, we show that the first condition implies the third. Suppose that the first condition holds, so that $E_F(\nu,\mu) \le E_F(\mu,\mu)$ for all $\nu\in\P(S)$. Notice that $\sup_{s\in\supp(\mu)} E_F(\delta_s,\mu) \le E_F(\mu,\mu)$, and also that
		\begin{align*}
			\sup_{s\in\supp(\mu)} E_F(\delta_s,\mu) &= \int_S \left(\sup_{s\in\supp(\mu)}E_F(\delta_s,\mu)\right) d\mu(s') \\
			&= \int_{\supp(\mu)} \left(\sup_{s\in\supp(\mu)}E_F(\delta_s,\mu)\right) d\mu(s') \\
			&\ge \int_{\supp(\mu)} E_F(\delta_{s'},\mu) d\mu(s') \\
			&= \int_S F_\mu(s') d\mu(s') \\
			&= E_F(\mu,\mu).
		\end{align*}
		Hence, $\sup_{s\in\supp(\mu)}E_F(\delta_s,\mu) = E_F(\mu,\mu)$. Suppose for the sake of contradiction that there exists $s'\in\supp(\mu)$ such that $E_F(\delta_{s'},\mu) < \sup_{s\in\supp(\mu)} E_F(\delta_s,\mu) = E_F(\mu,\mu)$. Since $F_\mu$ is a continuous real-valued function on $S$, the preimage $U \coloneqq F_\mu^{-1}((-\infty,E_F(\mu,\mu))) = \{s\in S : F_\mu(s) < E_F(\mu,\mu)\}$ is open and contains $s'$, and hence it must be the case that $\mu(U) > 0$ by definition of $\supp(\mu)$. Thus, since the Lebesgue integral of a positive function over a set of positive measure is positive, we find that
		\begin{align*}
			0 &= E_F(\mu,\mu) - E_F(\mu,\mu) \\
			  &= \int_S (E_F(\mu,\mu) - F_\mu(s)) d\mu(s) \\
			  &= \int_U (E_F(\mu,\mu) - F_\mu(s))d\mu(s) + \int_{S\setminus U}(E_F(\mu,\mu) - E_F(\delta_s,\mu))d\mu(s) \\
			  &\ge \int_U (E_F(\mu,\mu) - F_\mu(s))d\mu(s) \\
			  &> 0,
		\end{align*}
		which is a contradiction. Hence, it must be the case that $E_F(\delta_{s'},\mu) = \sup_{s \in \supp(\mu)} E_F(\delta_s,\mu) = E_F(\mu,\mu)$ for all $s'\in\supp(\mu)$. Therefore, $F_\mu(s') = E_F(\delta_{s'},\mu) = E_F(\mu,\mu) \ge E_F(\nu,\mu)$ for all $\nu\in\P(S)$ and all $s'\in\supp(\mu)$, and in particular, we find that $F_\mu(s') \ge E_F(\delta_s,\mu) = F_\mu(s)$ for all $s\in S$ and all $s'\in\supp(\mu)$, so the third condition holds.
	\end{proofE}
	}{}

	\autoref{prop: equivalent_nash} shows that at a Nash equilibrium state $\mu\in\P(S)$, every strategy $s' \in S$ that is in use (meaning that $s'\in\supp(\mu)$) must have maximal average payoff $F_\mu(s')$ compared to all other possible strategies $s\in S$. From the contrapositive viewpoint, this shows that a strategy $s' \in S$ whose average payoff $F_\mu(s')$ is strictly less than that of some other strategy will not be employed at a Nash equilibrium state $\mu$.

	\hl{In general, there may be more than one Nash equilibrium of a game $F$. Even in this case, the following result unveils advantageous topological characteristics of $\NE(F)$.}

	\iftoggle{longform}{
	\begin{propositionE}[][end,restate,text link=]
	}{
	\begin{propositionE}[][end,text link=]
	}
		\label{prop: ne_closed}
		Consider a game $F\colon \P(S) \to C(S)$. If $\theta_\nu \colon \P(S) \to \R$ defined by $\theta_\nu(\mu) = E_F(\nu,\mu) - E_F(\mu,\mu)$ is weakly continuous for all $\nu\in\P(S)$, then $\NE(F)$ is weakly compact.
	\end{propositionE}

	\iftoggle{longform}{
	\begin{proofE}
		It holds that \iftoggle{thesis}{
		\begin{align*}
			\NE(F) &= \{\mu\in\P(S) : \text{$E_F(\nu,\mu) - E_F(\mu,\mu) \le 0$ for all $\nu\in\P(S)$}\} \\
			       &= \bigcap_{\nu\in\P(S)}\{\mu\in\P(S) : E_F(\nu,\mu) - E_F(\mu,\mu) \le 0\}.
		\end{align*}
		}{$\NE(F) = \{\mu\in\P(S) : \text{$E_F(\nu,\mu) - E_F(\mu,\mu) \le 0$ for all $\nu\in\P(S)$}\} = \bigcap_{\nu\in\P(S)}\{\mu\in\P(S) : E_F(\nu,\mu) - E_F(\mu,\mu) \le 0\}$.} For all $\nu\in\P(S)$, the set $\{\mu\in \P(S) : E_F(\nu,\mu) - E_F(\mu,\mu) \le 0\}$ is the preimage of the closed set $(-\infty,0]$ under the map $\theta_\nu$. Hence, if this map is weakly continuous, then $\NE(F)$ is weakly closed. Since $\P(S)$ is weakly compact by \autoref{lem: P(S)_compact}, the weakly closed subset $\NE(F) \subseteq \P(S)$ must also be weakly compact.
	\end{proofE}
	}{}

	Together with \autoref{prop: ne_closed}, the following result shows that $\NE(F)$ is weakly compact whenever the game $F$ is weakly continuous. This compactness result is of particular technical importance in our stability proofs of \autoref{sec: dissipativity_theory}, {\color{black}as it gives a mild sufficient condition for ensuring that evolutionary games do not drift off towards unattainable equilibria ``at infinity'' (see \autoref{ex: at_infinity} for such an equilibrium ``at infinity'').}

	\iftoggle{longform}{
	\begin{propositionE}[][end,restate,text link=]
	}{
	\begin{propositionE}[][end,text link=]
	}
		\label{prop: weakly_continuous_payoff}
		Consider a game $F\colon \P(S) \to C(S)$. If $F$ is weakly continuous, then $\theta_\nu \colon \P(S) \to \R$ defined by $\theta_\nu(\mu) = E_F(\nu,\mu) - E_F(\mu,\mu)$ is weakly continuous for all $\nu\in\P(S)$.
	\end{propositionE}

	\iftoggle{longform}{
	\begin{proofE}
		Suppose that $F$ is weakly continuous and let $\nu\in\P(S)$. Since $S$ is a metric space, the weak topology on $\P(S)$ is metrizable \citep[Theorem~11.3.3]{dudley2002real}. Therefore, the weak topology on $\P(S)$ is first-countable and hence functions with domain $\P(S)$ are weakly continuous if they are weakly sequentially continuous. Thus, to prove the claim, it suffices to show that $\theta_\nu$ is weakly sequentially continuous. To this end, let $\{\mu_n \in \P(S) : n\in\N\}$ be a sequence that converges weakly to $\mu\in\P(S)$. Then we have that
		\iftoggle{thesis}{
		\begin{align*}
			|E_F(\nu,\mu_n) - E_F(\nu,\mu)| &= |\left<F(\mu_n),\nu\right> - \left<F(\mu),\nu\right>| \\
			&= |\left<F(\mu_n)-F(\mu),\nu\right>| \\
			&\le \|F(\mu_n) - F(\mu)\|_\infty \|\nu\|_\TV \\
			&\to 0
		\end{align*}
		}{
		\begin{equation*}
			|E_F(\nu,\mu_n) - E_F(\nu,\mu)| = |\left<F(\mu_n),\nu\right> - \left<F(\mu),\nu\right>| = |\left<F(\mu_n)-F(\mu),\nu\right>| \le \|F(\mu_n) - F(\mu)\|_\infty \|\nu\|_\TV \to 0
		\end{equation*}
		}
		since $\|\nu\|_\TV = 1$ and $F(\mu_n) \to F(\mu)$ in $C(S)$ with the topology induced by $\|\cdot\|_\infty$ due to weak continuity of $F$. Furthermore, we have that
		\begin{align*}
			|E_F(\mu_n,\mu_n) - E_F(\mu,\mu)| &= |\left<F(\mu_n),\mu_n\right> - \left<F(\mu),\mu\right>| \\
			&\le |\left<F(\mu_n),\mu_n\right> - \left<F(\mu),\mu_n\right>| + |\left<F(\mu),\mu_n\right> - \left<F(\mu),\mu\right>| \\
			&= |\left<F(\mu_n) - F(\mu),\mu_n\right>| + |\left<F(\mu),\mu_n - \mu\right>| \\
			&\le \|F(\mu_n) - F(\mu)\|_\infty \|\mu_n\|_\TV + |\left<F(\mu),\mu_n - \mu\right>| \\
			&= \|F(\mu_n) - F(\mu)\|_\infty + |\left<F(\mu),\mu_n-\mu\right>| \\
			&\to 0
		\end{align*}
		since again $F(\mu_n) \to F(\mu)$ by weak continuity of $F$, and since $\left<F(\mu),\mu_n - \mu\right> \to 0$ by definition of weak convergence of $\mu_n$ to $\mu$. Therefore, we conclude that
		\begin{equation*}
			\theta_\nu(\mu_n) = E_F(\nu,\mu_n) - E_F(\mu_n,\mu_n) \to E_F(\nu,\mu) - E_F(\mu,\mu) = \theta_\nu(\mu),
		\end{equation*}
		which proves the claim.
	\end{proofE}
	}{}

	\iftoggle{longform}{\hl{See \autoref{prop: unique_gess}, \autoref{cor: unique_nash}, and \autoref{prop: monotone_convex_ne} in \autoref{sec: supplementary} for conditions under which the (Nash) equilibria of a game are unique or constitute a convex set.}}{}

		\section{Evolutionary Dynamics}
	\label{sec: evolutionary_dynamics}
    
	In this section, the population game $F$ is endowed with dynamical behavior. Such dynamics are used to model the evolutionary aspects of a population playing out a game, wherein players revise their strategies over time according to the game's current payoff profile. Our infinite-strategy analogue of the evolutionary dynamics models considered in \citet{fox2013population} and \citet{arcak2021dissipativity} is formalized as follows.

	\begin{definition}
		\label{def: edm}
		Consider a game $F\colon \P(S) \to C(S)$. Let $\mu_0\in\P(S)$ and let $v\colon \P(S) \times C(S) \to \M(S)$. The differential equation
		\begin{equation}
		\begin{aligned}
			\dot{\mu}(t) &= v(\mu(t),\rho(t)), \\
			\rho(t) &= F(\mu(t)), \\
			\mu(0) &= \mu_0,
		\end{aligned}
			\label{eq: edm}
		\end{equation}
		is called an \emph{evolutionary dynamics model (EDM)}. The measure $\mu_0$ is called the \emph{initial state} and the mapping $v$ is called the \emph{dynamics map}. A strongly differentiable mapping $\mu\colon [0,\infty) \to \P(S)$ satisfying \eqref{eq: edm} is called a \emph{solution to the EDM}.
	\end{definition}

    {\color{black}The EDM \eqref{eq: edm} can be intuitively understood as a feedback interconnected control system, with $\mu$ being the state of the system to be controlled, $v$ being the nonautonomous ``open-loop'' dynamics (i.e., the plant), and $\rho(t) = F(\mu(t))$ being the static feedback controller that steers the state $\mu$ via the game's incentives. \autoref{fig: diagram} graphically illustrates this control-theoretic perspective that we employ throughout the paper. In what follows, we will leverage this framework to give \emph{separate} conditions on the open-loop dynamics $v$ and on the controller $F$ that ensure stability of the closed-loop evolutionary dynamics. Before doing so, we give examples of some of the most commonly studied evolutionary dynamics models, and also formalize the notions of\iftoggle{thesis}{ dynamic}{} stability to be considered.}

    \begin{figure}[ht]
    \centering
    \begin{tikzpicture}
    [plant/.style={very thick,black,fill=orange!10,draw,rounded corners},
    control/.style={very thick,black,fill=orange!10,draw,rounded corners}]

    \node [plant] (plant) at (0,0) {
		\begin{tabular}{c}
		\textbf{Open-loop dynamics}\\
		$\dot{\mu}(t) = v(\mu(t),\rho(t))$
		\end{tabular}
		};

    \node [control,below=1.5em of plant] (control) {
		\begin{tabular}{c}
		\textbf{Feedback payoffs}\\
		$\rho(t) = F(\mu(t))$
		\end{tabular}
		};

    \node [below=0.75em of plant,shape=coordinate] (phantommid) {};
    \node [left=7em of phantommid,shape=coordinate] (phantomleft) {};
    \node [right=7em of phantommid,shape=coordinate] (phantomright) {};
    \node [right=3.5em of plant,shape=coordinate] (phantomfar) {};

    \path[-latex,draw,thick] (plant.east) -| (phantomright) |- (control.east);
    \path[-latex,draw,thick] (plant.east) -- (phantomfar) node [right] {$\mu(t)$};
    \path[-latex,draw,thick] (control.west) -| (phantomleft) node [left] {$\rho(t)$} |- (plant.west);

\end{tikzpicture} 
    \caption{{\color{black} We study the evolutionary dynamics model \eqref{eq: edm} from a control-theoretic lens by viewing the evolution as ``open-loop'' dynamics $v$ controlled by the game's feedback payoffs $\rho(t) = F(\mu(t))$.}}
    \label{fig: diagram}
    \end{figure}

    {\color{black}
	\begin{example}
		\label{ex: bnn}
		Let $\lambda\in\P(S)$ be a fixed reference probability measure with full support. This reference measure is commonly taken as that of a uniform distribution, but more general \color{black}probability measures may be used to model the case where strategies are chosen at nonuniform revision rates. The \emph{Brown-von Neumann-Nash (BNN) dynamics} are given by the EDM \eqref{eq: edm} with closed-loop dynamics defined by
		\iftoggle{twocol}{
		\begin{align*}
			v(\mu,F(\mu))(B) &= \int_B \sigma_+(s,\mu) d\lambda(s) \\
			&\qquad - \mu(B)\int_S \sigma_+(s,\mu) d\lambda(s)
		\end{align*}
		}{
		\begin{equation*}
			v(\mu,F(\mu))(B) = \int_B \sigma_+(s,\mu) d\lambda(s) - \mu(B)\int_S \sigma_+(s,\mu) d\lambda(s)
		\end{equation*}
		}
		for all $B\in\B(S)$, where
		\begin{equation*}
			\sigma_+(s,\mu) = \max\{0,E_F(\delta_s,\mu) - E_F(\mu,\mu)\}
		\end{equation*}
            for all $s\in S$ and all $\mu\in\P(S)$. Here, the mean payoff function associated to $\mu\in\P(S)$ takes the form
		\begin{equation*}
			F_\mu(s) = \int_S f(s,s') d\mu(s'),
		\end{equation*}
		with $f\colon S\times S\to\R$ being a bounded measurable function that gives the payoff $f(s,s')$ to a player choosing strategy $s$ when an opponent chooses strategy $s'$. Intuitively, $\sigma_+(s,\mu)$ can be thought of as quantifying how well strategy $s$ performs relative to the average performance of the entire population $\mu$. Therefore, one may interpret the BNN dynamics as steering the population towards strategies that perform better than average, at a rate proportional to how much better those strategies are performing. 
        
        It is easy to see that the BNN dynamics may be written in the interconnected form \eqref{eq: edm} with the dynamics map given by
		\iftoggle{twocol}{
		\begin{align*}
			v(\mu,\rho)(B) &= \int_B \max\{0,\left<\rho,\delta_s\right> - \left<\rho,\mu\right>\} d\lambda(s) \\
			&\qquad - \mu(B) \int_S \max\{0,\left<\rho,\delta_s\right> - \left<\rho,\mu\right>\} d\lambda(s).
		\end{align*}
		}{
		\begin{equation*}
			v(\mu,\rho)(B) = \int_B \max\{0,\left<\rho,\delta_s\right> - \left<\rho,\mu\right>\} d\lambda(s) - \mu(B) \int_S \max\{0,\left<\rho,\delta_s\right> - \left<\rho,\mu\right>\} d\lambda(s).
		\end{equation*}
		}
		See \citet{hofbauer2009brown} for a thorough study on the BNN dynamics over infinite strategy sets.
	\end{example}
    }

        {\color{black}
	\begin{example}
		\label{ex: pcd}
		Let $\lambda$ be a fixed reference measure as in \autoref{ex: bnn}. Furthermore, let $\gamma \colon S\times S\times C(S) \to \Rpl$ be a continuous and bounded map, termed the \emph{conditional switch rate}, such that $\gamma(s,s',\rho)$ encodes the rate at which players switch from strategy $s\in S$ to strategy $s'\in S$ whenever the strategies' payoffs are described by the function $\rho\in C(S)$. Assume that the conditional switch rate satisfies \emph{sign-preservation}, given by
		\begin{equation*}
			\sign(\gamma(s,s',\rho)) = \sign(\max\{0,\rho(s')-\rho(s)\})
		\end{equation*}
		for all $s,s'\in S$ and all $\rho\in C(S)$. Sign-preservation ensures that the conditional switch rate from strategy $s$ to strategy $s'$ is positive if and only if $s'$ has higher payoff than $s$ according to the function $\rho$. The \emph{pairwise comparison dynamics} are given by the EDM \eqref{eq: edm} with
		\iftoggle{twocol}{
		\begin{align*}
			v(\mu,\rho)(B) &= \int_S \int_B \gamma(s,s',\rho) d\lambda(s') d\mu(s) \\
			&\qquad - \int_S \int_B \gamma(s',s,\rho) d\mu(s') d\lambda(s)
		\end{align*}
		}{
		\begin{equation*}
			v(\mu,\rho)(B) = \int_S \int_B \gamma(s,s',\rho) d\lambda(s') d\mu(s) - \int_S \int_B \gamma(s',s,\rho) d\mu(s') d\lambda(s)
		\end{equation*}
		}
		for all $B\in\B(S)$. Notice that the nonautonomous portion of the dynamics, defined by this dynamics map $v$, is entirely determined by the conditional switch rate $\gamma$. Intuitively, the first term in the dynamics $v$ models the average ``inflow'' of players switching to strategies in the set $B$, whereas the second term models the average ``outflow'' of players switching from strategies in $B$ to potentially different strategies in all of $S$.

		When $\gamma$ takes the form $\gamma(s,s',\rho) = \max\{0,\rho(s')-\rho(s)\}$, the pairwise comparison dynamics reduce to the famous \emph{Smith dynamics}, introduced in \citet{smith1984stability}. If, for all $s'\in S$, there exists some continuous function $\phi_{s'}\colon\R\to\Rpl$ such that the conditional switch rate satisfies
		\begin{equation*}
			\gamma(s,s',\rho) = \phi_{s'}(\rho(s') - \rho(s))
		\end{equation*}
		for all $s\in S$ and all $\rho\in C(S)$, then the pairwise comparison dynamics are said to be \emph{impartial}. The Smith dynamics are seen to be impartial by taking $\phi_{s'}(\cdot) = \max\{0,\cdot\}$ for all $s'$. See \citet{cheung2014pairwise} for a thorough study on the pairwise comparison dynamics over infinite strategy sets.
	\end{example}
        }

    {\color{black}The idea of studying evolutionary dynamics models, like those in \autoref{ex: bnn} and \autoref{ex: pcd}, as a feedback control system} was proposed in \citet{fox2013population} and further studied in \citet{arcak2021dissipativity}. {\color{black}This was done} as a means to derive finite-strategy stability results based on the idea that interconnections of energy-dissipating systems result in a \iftoggle{thesis}{closed-loop system that is dynamically stable}{stable closed-loop system}. This allows one to prove stability of the overall evolutionary dynamics model by studying the dissipativity properties of the (nonautonomous) system and input in a modular fashion. \iftoggle{thesis}{To the best of our knowledge, our}{Our} work is the first to generalize this modular dissipativity approach to evolutionary games with infinite strategy sets---prior works on infinite strategy sets primarily prove stability in a closed-loop black-box fashion on a case-by-case basis, e.g., \citet{oechssler2001evolutionary,hingu2020superiority} for replicator dynamics, \citet{hofbauer2009brown} for Brown-von Neumann-Nash dynamics, and \citet{cheung2014pairwise} for pairwise comparison dynamics, as well as the references therein and subsequent works.

	\iftoggle{thesis}{
	\hl{We make use of the following existence and uniqueness assumption throughout the remainder of the \iftoggle{thesis}{thesis}{paper}.}

	\begin{assumption}
		\label{ass: exist_edm_solution}
		Consider a dynamics map $v\colon\P(S) \times C(S) \to \M(S)$. For every initial state $\mu_0\in\P(S)$, there exists a unique solution to the EDM \eqref{eq: edm} with initial state $\mu_0$ and dynamics map $v$.
	\end{assumption}

	\autoref{ass: exist_edm_solution}}{Throughout {\color{black}the remainder of} this \iftoggle{thesis}{thesis}{paper}, we will always assume that, for every initial state $\mu_0\in \P(S)$, the EDM \eqref{eq: edm} admits a unique solution. This assumption} holds for the BNN dynamics of \autoref{ex: bnn} \citep[Theorem~1]{hofbauer2009brown}, and also holds for the pairwise comparison dynamics of \autoref{ex: pcd} \citep[Theorem~1]{cheung2014pairwise} under some mild regularity conditions on $F$. It is also easy to see that, for these dynamics, $v(\mu,\rho)(S) = 0$ for all $\mu\in\P(S)$ and all $\rho\in C(S)$, and hence the codomain of these dynamics maps can be taken to be $T\P(S)$.\iftoggle{longform}{\hl{ For more discussion on the characteristics and existence of solutions to EDMs, see \autoref{sec: supplementary}.}}{}

	\subsection{Dynamic Notions of Stability}

	We now formally introduce the notions of dynamic equilibria and stability with which we are concerned. \iftoggle{thesis}{
	\begin{definition}
		\label{def: rest_point}
		A population state $\mu \in \P(S)$ is a \emph{rest point of the EDM \eqref{eq: edm}} if $v(\mu, F(\mu)) = 0$.
	\end{definition}}{%
	A population state $\mu\in\P(S)$ is said to be a \emph{rest point of the EDM \eqref{eq: edm}} if $v(\mu,F(\mu)) = 0$.} The following condition, which is solely a property of the nonautonomous dynamics defined by the dynamics map $v$, ensures that the rest points and Nash equilibria coincide for the EDM under the feedback interconnection \eqref{eq: edm}.

	\begin{definition}
		\label{def: ns}
		A map $v\colon \P(S) \times C(S) \to \M(S)$ is \emph{Nash stationary} if, for all $\mu\in\P(S)$ and all $\rho\in C(S)$, it holds that
		\begin{equation*}
			v(\mu,\rho) = 0
		\end{equation*}
		if and only if
		\begin{equation*}
			\left<\rho,\nu\right> \le \left<\rho,\mu\right> ~ \text{for all $\nu\in\P(S)$}.
		\end{equation*}
	\end{definition}
	\iftoggle{longform}{
	\begin{propositionE}[][end,restate,text link=]
	}{
	\begin{propositionE}[][end,text link=]
	}
		\label{prop: ns_implies_rest_equals_nash}
		Consider a game $F\colon\P(S) \to C(S)$ and let $v \colon \P(S) \times C(S) \to \M(S)$. If $v$ is Nash stationary, then the set of rest points of the EDM \eqref{eq: edm} equals $\NE(F)$.
	\end{propositionE}

	\iftoggle{longform}{
	\begin{proofE}
		Suppose that $v$ is Nash stationary. Let $\mu\in\P(S)$ be a rest point of the EDM \eqref{eq: edm} with dynamics map $v$. Then $v(\mu,F(\mu)) = 0$. Since $v$ is Nash stationary, this implies that $E_F(\nu,\mu) - E_F(\mu,\mu) = \left<F(\mu),\nu\right> - \left<F(\mu),\mu\right> \le 0$ for all $\nu\in\P(S)$. Thus, $\mu\in\NE(F)$. On the other hand, if $\mu\in\NE(F)$, then $\left<F(\mu),\nu\right> - \left<F(\mu),\mu\right> = E_F(\nu,\mu) - E_F(\mu,\mu) \le 0$, so $v(\mu,F(\mu))=0$ since $v$ is Nash stationary. Thus, $\mu$ is a rest point of the EDM \eqref{eq: edm} with dynamics map $v$.
	\end{proofE}
	}{}

	\autoref{prop: ns_implies_rest_equals_nash} shows that if an evolutionary game's population state converges to a rest point under the EDM with a Nash stationary dynamics map, then the population state converges to a Nash equilibrium. In other words, Nash stationarity ensures that all\iftoggle{thesis}{ dynamically}{} stable rest points have game-theoretic importance. We now recall that the popular BNN dynamics and pairwise comparison dynamics both satisfy Nash stationarity.\iftoggle{longform}{\footnote{Technically, our definition of Nash stationarity, which is a property of the dynamics map $v$ viewed as a nonautonomous system, is slightly different than the definitions used in \citet{hofbauer2009brown} and \citet{cheung2014pairwise}, which are properties of the closed-loop interconnection \eqref{eq: edm}. For self-containedness, we prove \autoref{prop: bnn_pcd_satisfy_ns} in \iftoggle{proofs}{\autoref{sec: proofs}}{the online technical report \citet{anderson2025dissipativity-long}} using our definition.}}{}

	\iftoggle{longform}{
	\begin{propositionE}[\citealp{hofbauer2009brown,cheung2014pairwise}][end,restate,text link=]
	}{
	\begin{propositionE}[\citealp{hofbauer2009brown,cheung2014pairwise}][end,text link=]
	}
		\label{prop: bnn_pcd_satisfy_ns}
		If $v\colon \P(S) \times C(S) \to \M(S)$ is the dynamics map for either the BNN dynamics of \autoref{ex: bnn} or the pairwise comparison dynamics of \autoref{ex: pcd}, then $v$ is Nash stationary.
	\end{propositionE}

	\iftoggle{longform}{
	\begin{proofE}
		Let $\mu\in\P(S)$ and let $\rho\in C(S)$. First consider the BNN dynamics of \autoref{ex: bnn}. We have that
		\begin{equation*}
			v(\mu,\rho)(B) = \int_B \max\{0,\left<\rho,\delta_s\right> - \left<\rho,\mu\right>\} d\lambda(s) - \mu(B) \int_S \max\{0,\left<\rho,\delta_s\right> - \left<\rho,\mu\right>\} d\lambda(s)
		\end{equation*}
		for all $B\in\B(S)$. If $\left<\rho,\nu\right> \le \left<\rho,\mu\right>$ for all $\nu\in\P(S)$, then it follows immediately that $v(\mu,\rho)(B) = 0$ for all $B\in\B(S)$, and hence $v(\mu,\rho) = 0$.

		On the other hand, suppose that $v(\mu,\rho) = 0$. \iftoggle{thesis}{Suppose that $$\int_S \max\{0,\left<\rho,\delta_s\right> - \left<\rho,\mu\right>\} d\lambda(s) = 0.$$}{Suppose that $\int_S \max\{0,\left<\rho,\delta_s\right> - \left<\rho,\mu\right>\} d\lambda(s) = 0$.} Then we find that
		\begin{equation*}
			\int_B \max\{0,\left<\rho,\delta_s\right> - \left<\rho,\mu\right>\} d\lambda(s) = 0
		\end{equation*}
		for all $B\in\B(S)$. Hence, $\max\{0,\left<\rho,\delta_s\right> - \left<\rho,\mu\right>\} = 0$ for $\lambda$-almost every $s\in S$. Since $\lambda$ has full support by assumption and $s\mapsto \max\{0,\left<\rho,\delta_s\right> - \left<\rho,\mu\right>\}$ is continuous, this shows that $\max\{0,\left<\rho,\delta_s\right> - \left<\rho,\mu\right>\} = 0$ for all $s\in S$. Hence, $\left<\rho,\delta_s\right> \le \left<\rho,\mu\right>$ for all $s\in S$. Since $S$ is compact and $\rho$ is continuous, the optimization $\sup_{s\in S}\rho(s)$ is attained by some $s'\in S$. Therefore, for all $\nu\in\P(S)$ it holds that $\left<\rho,\nu\right> = \int_S \rho(s) d\nu(s) \le \int_S \rho(s') d\nu(s) = \rho(s') = \left<\rho,\delta_{s'}\right> \le \left<\rho,\mu\right>$. Now suppose that the other case holds, namely, that $\int_S \max\{0,\left<\rho,\delta_s\right> - \left<\rho,\mu\right>\}d\lambda(s) > 0$. Then it holds that
		\begin{equation*}
			\mu(B) = \int_B \frac{\max\{0,\left<\rho,\delta_{s'}\right> - \left<\rho,\mu\right>\}}{\int_S \max\{0,\left<\rho,\delta_{s}\right> - \left<\rho,\mu\right>\} d\lambda(s)} d\lambda(s')
		\end{equation*}
		for all $B\in\B(S)$. Suppose for the sake of contradiction that there exists $\tilde{s}\in S$ such that $\left<\rho,\delta_{\tilde{s}}\right> - \left<\rho,\mu\right> > 0$. Then, by continuity of $s'\mapsto \left<\rho,\delta_{s'}\right> - \left<\rho,\mu\right>$, the preimage $\{s'\in S : \left<\rho,\delta_{s'}\right> - \left<\rho,\mu\right> > 0\}$ is open and contains $\tilde{s}$, and hence it must be the case that $\lambda(\{s'\in S : \left<\rho,\delta_{s'}\right> - \left<\rho,\mu\right> > 0\}) > 0$ by definition of $\supp(\lambda)$ and the fact that $\lambda$ has full support. Therefore, since the Lebesgue integral of a positive function over a set of positive measure is positive, we find that
		\begin{align*}
			\left<\rho,\mu\right> &= \int_S \rho d\mu \\
			&= \int_S \rho(s') \frac{\max\{0,\left<\rho,\delta_{s'}\right> - \left<\rho,\mu\right>\}}{\int_S \max\{0,\left<\rho,\delta_{s}\right> - \left<\rho,\mu\right>\} d\lambda(s)} d\lambda(s') \\
			&= \int_{\{s'\in S : \left<\rho,\delta_{s'}\right> - \left<\rho,\mu\right> > 0\}} \left<\rho,\delta_{s'}\right> \frac{\max\{0,\left<\rho,\delta_{s'}\right>- \left<\rho,\mu\right>\}}{\int_S \max\{0,\left<\rho,\delta_{s}\right> - \left<\rho,\mu\right>\} d\lambda(s)} d\lambda(s') \\
			&> \int_{\{s'\in S : \left<\rho,\delta_{s'}\right> - \left<\rho,\mu\right> > 0\}} \left<\rho,\mu\right> \frac{\max\{0,\left<\rho,\delta_{s'}\right> - \left<\rho,\mu\right>\}}{\int_S \max\{0,\left<\rho,\delta_{s}\right> - \left<\rho,\mu\right>\} d\lambda(s)} d\lambda(s') \\
			&= \left<\rho,\mu\right>,
		\end{align*}
		which is a contradiction. Hence, it must be that $\left<\rho,\delta_{\tilde{s}}\right> \le \left<\rho,\mu\right>$ for all $\tilde{s}\in S$. Arguing as in the prior case, this yields that $\left<\rho,\nu\right> \le \left<\rho,\mu\right>$ for all $\nu\in\P(S)$. Since this exhausts all cases to be considered, we conclude that indeed $v$ is Nash stationary.

		Now consider the pairwise comparison dynamics of \autoref{ex: pcd}. We have that
		\begin{equation*}
			v(\mu,\rho)(B) = \int_S \int_B \gamma(s,s',\rho)d\lambda(s') d\mu(s) - \int_S\int_B \gamma(s',s,\rho) d\mu(s') d\lambda(s)
		\end{equation*}
		for all $B\in\B(S)$. By \autoref{lem: equivalent_nash-general}, which we prove after completing the current proof, it holds that $\left<\rho,\nu\right> \le \left<\rho,\mu\right>$ for all $\nu\in\P(S)$ if and only if $\rho(s) \le \rho(s')$ for all $s\in S$ and all $s'\in\supp(\mu)$. Thus, if $\left<\rho,\nu\right> \le \left<\rho,\mu\right>$ for all $\nu\in\P(S)$, then $\max\{0,\rho(s) - \rho(s')\} = 0$ for all $s\in S$ and all $s'\in\supp(\mu)$, implying that $\sign(\max\{0,\rho(s) - \rho(s')\}) = 0$ for all $s\in S$ and all $s'\in\supp(\mu)$. Hence, since the conditional switch rate $\gamma$ satisfies sign-preservation by assumption, we find that
		\begin{equation*}
			\sign(\gamma(s',s,\rho)) = 0
		\end{equation*}
		for all $s\in S$ and all $s'\in\supp(\mu)$. This implies that $v(\mu,\rho)(B) = 0$ for all $B\in\B(S)$, and therefore that $v(\mu,\rho) = 0$.

		On the other hand, suppose that $v(\mu,\rho) = 0$. Define the measures $v_1(\mu,\rho),v_2(\mu,\rho)\in\M(S)$ by
		\begin{align*}
			v_1(\mu,\rho)(B) &\coloneqq \int_S \int_B \gamma(s,s',\rho) d\lambda(s')d\mu(s) = \int_B \int_S \gamma(s,s',\rho) d\mu(s) d\lambda(s'), \\
			v_2(\mu,\rho)(B) &\coloneqq \int_S \int_B \gamma(s',s,\rho) d\mu(s')d\lambda(s) = \int_B \int_S \gamma(s',s,\rho) d\lambda(s) d\mu(s'),
		\end{align*}
		so that $v(\mu,\rho) = v_1(\mu,\rho) - v_2(\mu,\rho)$. Since $v(\mu,\rho) = 0$, it holds that $v_1(\mu,\rho) = v_2(\mu,\rho)$, and hence $\left<\rho,v_1(\mu,\rho)\right> = \left<\rho,v_2(\mu,\rho)\right>$. Therefore,
		\begin{equation*}
			\int_S \rho(s') \int_S \gamma(s,s',\rho) d\mu(s) d\lambda(s') = \int_S \rho(s') \int_S \gamma(s',s,\rho)d\lambda(s) d\mu(s').
		\end{equation*}
		Hence,
		\begin{equation*}
			\int_S\int_S \rho(s') \gamma(s,s',\rho) d\mu(s) d\lambda(s') = \int_S\int_S \rho(s') \gamma(s',s,\rho) d\lambda(s) d\mu(s'),
		\end{equation*}
		implying that
		\begin{equation*}
			\int_S \int_S (\rho(s) - \rho(s')) \gamma(s',s,\rho) d\mu(s') d\lambda(s) = 0.
		\end{equation*}
		By sign-preservation of the conditional switch rate $\gamma$, it holds that \iftoggle{thesis}{$$\sign(\gamma(s',s,\rho)) = \sign(\max\{0,\rho(s) - \rho(s')\})$$}{$\sign(\gamma(s',s,\rho)) = \sign(\max\{0,\rho(s) - \rho(s')\})$} for all $s,s'\in S$, and therefore, if $\rho(s) \ge \rho(s')$, we find that $\gamma(s',s,\rho) \ge 0$ so that $(\rho(s)-\rho(s'))\gamma(s',s,\rho) \ge 0$, and similarly if $\rho(s) \le \rho(s')$, we find that $\gamma(s',s,\rho) = 0$ so that $(\rho(s) - \rho(s'))\gamma(s',s,\rho)=0$. Hence, $(\rho(s) - \rho(s'))\gamma(s',s,\rho) \ge 0$ for all $s,s'\in S$, and also $\int_S (\rho(s) - \rho(s'))\gamma(s',s,\rho)d\mu(s') \ge 0$ for all $s\in S$. Since $s\mapsto \int_S (\rho(s) - \rho(s'))\gamma(s',s,\rho) d\mu(s')$ is continuous (which follows from compactness of $S$ and continuity of $s'\mapsto (\rho(s)-\rho(s'))\gamma(s',s,\rho)$, together with the dominated convergence theorem), the preimage $\{s\in S : \int_S(\rho(s)-\rho(s'))\gamma(s',s,\rho)d\mu(s') > 0\}$ is open and therefore must be empty, for otherwise $\int_S \int_S (\rho(s) - \rho(s'))\gamma(s',s,\rho) d\mu(s') d\lambda(s) > 0$ as $\lambda$ has full support. Hence,
		\begin{equation*}
			\int_S (\rho(s) - \rho(s'))\gamma(s',s,\rho) d\mu(s') = 0 ~ \text{for all $s\in S$}.
		\end{equation*}
		Similarly, since $s' \mapsto (\rho(s) - \rho(s'))\gamma(s',s,\rho)$ is continuous for all $s\in S$, the preimage $\{s'\in S : (\rho(s) - \rho(s'))\gamma(s',s,\rho) > 0\}$ is open for all $s\in S$, and hence for all $s'\in\supp(\mu)$ it must be the case that
		\begin{equation*}
			(\rho(s) - \rho(s'))\gamma(s',s,\rho) = 0
		\end{equation*}
		for all $s\in S$. Thus, for all $s\in S$ and all $s'\in\supp(\mu)$, either $\rho(s) = \rho(s')$, or $\gamma(s',s,\rho) = 0$. In the latter case, we see by sign-preservation of the conditional switch rate that \iftoggle{thesis}{$$\sign(\max\{0,\rho(s)-\rho(s')\}) = \sign(\gamma(s',s,\rho)) = 0,$$}{$\sign(\max\{0,\rho(s)-\rho(s')\}) = \sign(\gamma(s',s,\rho)) = 0$,} and hence $\rho(s) \le \rho(s')$. Therefore, we conclude that
		\begin{equation*}
			\rho(s) \le \rho(s') ~ \text{for all $s\in S$ and all $s'\in\supp(\mu)$}.
		\end{equation*}
		By \autoref{lem: equivalent_nash-general}, this proves that $\left<\rho,\nu\right> \le \left<\rho,\mu\right>$ for all $\nu\in\P(S)$, and consequently that $v$ is Nash stationary.
	\end{proofE}
	}{}

	\iftoggle{longform}{
	\begin{lemmaE}[][all end,text link=]
		\label{lem: equivalent_nash-general}
		It holds that $\left<\rho,\nu\right> \le \left<\rho,\mu\right>$ for all $\nu\in\P(S)$ if and only if $\rho(s) \le \rho(s')$ for all $s\in S$ and all $s'\in\supp(\mu)$.
	\end{lemmaE}

	\begin{proofE}
		Suppose first that $\rho(s) \le \rho(s')$ for all $s\in S$ and all $s'\in\supp(\mu)$. Then, it holds that
		\begin{equation*}
			\rho(s) = \int_S \rho(s) d\mu(s') \le \int_S \rho(s') d\mu(s') = \left<\rho,\mu\right>
		\end{equation*}
		for all $s\in S$. Therefore, for all $\nu\in\P(S)$, we conclude that
		\begin{equation*}
			\left<\rho,\nu\right> = \int_S \rho(s) d\nu(s) \le \int_S \left<\rho,\mu\right> d\nu(s) = \left<\rho,\mu\right>,
		\end{equation*}
		which proves one direction of the lemma.

		On the other hand, suppose that $\left<\rho,\nu\right> \le \left<\rho,\mu\right>$ for all $\nu\in\P(S)$. Then we have that $\rho(s) = \left<\rho,\delta_s\right> \le \left<\rho,\mu\right>$ for all $s\in S$. Furthermore,
		\begin{align*}
			\sup_{s\in\supp(\mu)} \rho(s) &= \int_S \left(\sup_{s\in\supp(\mu)}\rho(s)\right) d\mu(s') \\
			&= \int_{\supp(\mu)} \left(\sup_{s\in\supp(\mu)}\rho(s)\right) d\mu(s') \\
			&\ge \int_{\supp(\mu)} \rho(s') d\mu(s') \\
			&= \int_S \rho(s') d\mu(s') \\
			&= \left<\rho,\mu\right>.
		\end{align*}
		Hence, $\sup_{s\in\supp(\mu)}\rho(s) = \left<\rho,\mu\right>$. Suppose for the sake of contradiction that there exists $s'\in\supp(\mu)$ such that $\rho(s') < \sup_{s\in\supp(\mu)}\rho(s) = \left<\rho,\mu\right>$. Since $\rho$ is a continuous real-valued function on $S$, the preimage $U\coloneqq \rho^{-1}((-\infty,\left<\rho,\mu\right>)) = \{s\in S : \rho(s) < \left<\rho,\mu\right>\}$ is open and contains $s'$, and hence it must be the case that $\mu(U) > 0$ by definition of $\supp(\mu)$. Thus, since the Lebesgue integral of a positive function over a set of positive measure is positive, we find that
		\begin{align*}
			0 &= \left<\rho,\mu\right> - \left<\rho,\mu\right> \\
			&= \int_S (\left<\rho,\mu\right> - \rho(s)) d\mu(s) \\
			&= \int_U(\left<\rho,\mu\right> - \rho(s))d\mu(s) + \int_{S\setminus U}(\left<\rho,\mu\right> - \rho(s))d\mu(s) \\
			&\ge \int_U (\left<\rho,\mu\right> - \rho(s))d\mu(s) \\
			&> 0,
		\end{align*}
		which is a contradiction. Hence, it must be the case that $\rho(s') = \sup_{s\in\supp(\mu)}\rho(s) = \left<\rho,\mu\right>$ for all $s'\in\supp(\mu)$. Therefore, $\rho(s') = \left<\rho,\mu\right> \ge \left<\rho,\nu\right>$ for all $\nu\in\P(S)$ and all $s'\in\supp(\mu)$, and in particular, we find that $\rho(s') \ge \left<\rho,\delta_s\right> = \rho(s)$ for all $s\in S$ and all $s'\in\supp(\mu)$. This concludes the proof.
	\end{proofE}
	}{}

	Although a population state being at a rest point ensures that the population's distribution of strategies remains constant for all time, the definition of rest point does not itself come equipped with an adequate notion of stability. For this, we turn to the classical definitions of Lyapunov stability and attraction. Since ultimately we are interested in convergence of a game's dynamics to \emph{some} Nash equilibrium, we deal with such definitions in the sense of sets. \hl{We present the definitions for general Banach spaces with topologies on them that may not be induced by the space's norm. This level of abstraction is particularly useful for extending our dissipativity theory to the case of dynamic payoff models, {\color{black}as we do in \iftoggle{longform}{\autoref{sec: dynamic_payoff}}{our online technical report \citet{anderson2025dissipativity-long}}.}}

	\begin{definition}
		\label{def: tau-lyapunov_stable}
		Consider a Banach space $X$ and a topology $\tau$ on $X$. Let $Y\subseteq X$ and let $v \colon Y \to X$. A $\tau$-compact set $P\subseteq Y$ is \emph{$\tau$-Lyapunov stable under $v$} if, for all relatively $\tau$-open sets $Q\subseteq Y$ containing $P$, there exists a relatively $\tau$-open set $R\subseteq Y$ containing $P$ such that every solution $x\colon [0,\infty) \to Y$ to the differential equation $\dot{x}(t) = v(x(t))$ with $x(0) = x_0 \in Y$ satisfies $x(t) \in Q$ for all $t\in[0,\infty)$ whenever $x(0) \in R$.
	\end{definition}
    
	\begin{definition}
		\label{def: tau-attracting}
		Consider a Banach space $X$ and a topology $\tau$ on $X$. Let $Y\subseteq X$ and let $v \colon Y \to X$. A $\tau$-compact set $P\subseteq Y$ is \emph{$\tau$-attracting under $v$ from $x_0\in Y$} if, for all relatively $\tau$-open sets $Q\subseteq Y$ containing $P$ and for all solutions $x\colon [0,\infty) \to Y$ to the differential equation $\dot{x}(t) = v(x(t))$ with $x(0) = x_0$, there exists $T\in [0,\infty)$ such that
		\begin{equation*}
			x(t) \in Q ~ \text{for all $t\in[T,\infty)$}.
		\end{equation*}
		The $\tau$-compact set $P$ is \emph{globally $\tau$-attracting under $v$} if it is $\tau$-attracting under $v$ from $x_0$ for all $x_0\in Y$.
	\end{definition}

    {\color{black}\autoref{def: tau-lyapunov_stable} and \autoref{def: tau-attracting} generalize the standard notions of Lyapunov stability and attraction to the infinite-dimensional setting. We will most often be concerned with applying these definitions to the set $P = \NE(F)$, as this is the set of population states we would like to ensure stability for. In doing so, \autoref{def: tau-lyapunov_stable} ensures that the game's population state remains near some Nash equilibrium whenever it starts near some Nash equilibrium, whereas \autoref{def: tau-attracting} ensures that the population state eventually converges towards $\NE(F)$.}
        
	In our main results, \autoref{thm: delta-dissipative} and \autoref{thm: monotone-passive} {\color{black}that follow}, we are concerned with stability with respect to the weak topology on $\P(S)$. Therein, $\tau$ coincides with the weak topology, which in our setting is strictly weaker than the topology induced by $\|\cdot\|_\TV$. \iftoggle{thesis}{For this case, we specialize the above stability definitions to the setting of evolutionary dynamics of the form \eqref{eq: edm}.

	\begin{definition}
		\label{def: weakly_lyapunov_stable}
		Consider a game $F\colon \P(S) \to C(S)$ and let $v\colon \P(S) \times C(S) \to \M(S)$. A weakly compact set $P\subseteq \P(S)$ is \emph{weakly Lyapunov stable under the EDM \eqref{eq: edm}} if $P$ is $\tau$-Lyapunov stable under $\mu \mapsto v(\mu,F(\mu))$ according to \autoref{def: tau-lyapunov_stable} with $X = \M(S)$, $Y = \P(S)$, and $\tau$ being the weak topology.
	\end{definition}

	\begin{definition}
		\label{def: attracting}
		Consider a game $F\colon \P(S) \to C(S)$ and let $v \colon \P(S) \times C(S) \to \M(S)$. A weakly compact set $P\subseteq \P(S)$ is \emph{weakly attracting under the EDM \eqref{eq: edm} from initial state $\mu_0\in\P(S)$} if $P$ is $\tau$-attracting under $\mu\mapsto v(\mu,F(\mu))$ from $\mu_0$ according to \autoref{def: tau-attracting} with $X = \M(S)$, $Y = \P(S)$, and $\tau$ being the weak topology. The weakly compact set $P$ is \emph{globally weakly attracting under the EDM \eqref{eq: edm}} if it is weakly attracting under the EDM \eqref{eq: edm} with initial state $\mu_0$ for all $\mu_0\in\P(S)$.
	\end{definition}

 	\autoref{def: weakly_lyapunov_stable} and \autoref{def: attracting} are equivalent to the\iftoggle{thesis}{ dynamic}{} stability notions used in \citet{cheung2014pairwise}, which are defined in terms of the Prokhorov metric on $\P(S)$ that metrizes the weak topology.
}{\hl{In this case, we use the terminology \emph{weakly Lyapunov stable} to mean $\tau$-Lyapunov stable, and \emph{weakly attracting} to mean $\tau$-attracting. Such definitions of weak Lyapunov stability and weak attraction are equivalent to the\iftoggle{thesis}{ dynamic}{} stability notions used in \citet{cheung2014pairwise}, which are defined in terms of the Prokhorov metric on $\P(S)$ that metrizes the weak topology.}}

		\section{Dissipativity Theory}
	\label{sec: dissipativity_theory}

	\hl{In this section, we present our main theoretical contributions. We begin by defining notions of dissipativity in our setting of infinite strategy sets, which we then use to prove our main results characterizing the\iftoggle{thesis}{ dynamic}{} stability of Nash equilibria. \iftoggle{thesis}{To the best of our knowledge, the}{The} following definition is new to the literature on evolutionary games over infinite strategy sets.}

	\begin{definition}
		\label{def: delta-dissipative}
		A map $v\colon \P(S) \times C(S) \to \M(S)$ is \emph{$\delta$-dissipative with supply rate $w\colon \M(S) \times C(S) \to \R$} if there exist $\sigma\colon \P(S) \times C(S) \to \Rpl$ and $\storage \colon \P(S) \times C(S) \to \Rpl$ that extends to a map $\overline{\storage}\colon U\times C(S) \to \R$ with strongly open $U\subseteq \M(S)$ containing $\P(S)$, such that the following conditions hold:
		\begin{enumerate}
			\item $\overline{\storage}$ is weak-$\infty$-continuous.
			\item $\overline{\storage}$ is Fr\'echet differentiable.
			\item For all $\mu\in\P(S)$, all $\rho\in C(S)$, and all $\eta\in C(S)$, it holds that
			\iftoggle{twocol}{
			\begin{equation}
				\begin{aligned}
				&\partial_1 \overline{\storage}(\mu,\rho) v(\mu,\rho) + \partial_2 \overline{\storage}(\mu,\rho)\eta \\
				&\qquad \le -\sigma(\mu,\rho) + w(v(\mu,\rho),\eta).
				\end{aligned}
				\label{eq: delta-dissipative-1}
			\end{equation}
			}{
			\begin{equation}
				\partial_1 \overline{\storage}(\mu,\rho) v(\mu,\rho) + \partial_2 \overline{\storage}(\mu,\rho)\eta \le -\sigma(\mu,\rho) + w(v(\mu,\rho),\eta).
				\label{eq: delta-dissipative-1}
			\end{equation}
			}
			\item For all $\mu\in\P(S)$ and all $\rho\in C(S)$, it holds that
			\begin{equation}
				\text{$\storage(\mu,\rho)=0$ if and only if $v(\mu,\rho) = 0$}.
				\label{eq: delta-dissipative-2}
			\end{equation}
		\end{enumerate}
		If, additionally, $(\mu,\rho) \mapsto \partial_1 \overline{\storage}(\mu,\rho)$ and $(\mu,\rho) \mapsto \partial_2 \overline{\storage}(\mu,\rho)$ are weak-$\infty$-continuous, every partial Fr\'echet derivative $\partial_1 \overline{\storage}(\mu,\rho)$ is weakly continuous, and
		\begin{equation}
			\text{$\sigma(\mu,\rho) = 0$ if and only if $v(\mu,\rho) = 0$}
			\label{eq: delta-dissipative-3}
		\end{equation}
		for all $\mu\in\P(S)$ and all $\rho\in C(S)$, then $v$ is \emph{strictly $\delta$-dissipative with supply rate $w$}.
	\end{definition}

	\iftoggle{thesis}{The ``$\delta$'' in the definition of $\delta$-dissipative is short for ``differentially;'' it does not refer to a quantity $\delta$. Also note}{Note} that $\delta$-dissipativity is solely a property of the nonautonomous {\color{black}(``open-loop'')} dynamics defined by the dynamics map $v$. {\color{black}As mentioned in \autoref{ex: pcd}}, the pairwise comparison dynamics are entirely determined by some conditional switch rate function $\gamma \colon S\times S\times C(S) \to \Rpl$, and therefore $\delta$-dissipativity may be viewed as a property of the conditional switch rate function in such a setting. {\color{black}Intuitively, the open-loop evolutionary dynamics are $\delta$-dissipative if the rate of change of internally stored energy, characterized by the ``storage function'' $\storage$, is less than the energy supply rate $w$.}

	\iftoggle{thesis}{We also remark that \eqref{eq: delta-dissipative-1}, \eqref{eq: delta-dissipative-2}, and \eqref{eq: delta-dissipative-3} are the \hl{key criteria} for $\delta$-dissipativity; the continuity and differentiability conditions in \autoref{def: delta-dissipative} are regularity conditions \hl{that are indispensable technical criteria for our stability proofs that follow}. In addition, we need the following regularity conditions on the game $F$ and the dynamics map $v$, which, as we will see in the \hl{case study} of \autoref{sec: failure_finite}, are of utmost importance in ensuring that\iftoggle{thesis}{ dynamic}{} stability actually holds.}{We need the following regularity conditions on the game $F$ and the dynamics map $v$, which, as we will see in the \hl{case study} of \autoref{sec: failure_finite}, are of utmost importance in ensuring that\iftoggle{thesis}{ dynamic}{} stability actually holds.}

	\begin{assumption}
		\label{ass: game}
		Consider a game $F\colon \P(S) \to C(S)$. The following hold:
		\begin{enumerate}
			\item $F$ is weakly continuous.
			\item $F$ extends to a weakly continuous Fr\'echet differentiable map $\overline{F}\colon U' \to C(S)$ defined on a strongly open set $U' \subseteq \M(S)$ containing $\P(S)$.
		\end{enumerate}
	\end{assumption}
    
	\begin{assumption}
		\label{ass: game_derivative}
		Consider a game $F\colon \P(S) \to C(S)$ that satisfies \autoref{ass: game}. It holds that \iftoggle{thesis}{$D\overline{F}$}{$D \overline{F} \colon U' \to \lin(\M(S),C(S))$} and every Fr\'echet derivative \iftoggle{thesis}{$D \overline{F}(\mu)$}{$D\overline{F}(\mu) \colon \M(S) \to C(S)$} are weakly continuous.
	\end{assumption}
    
	\begin{assumption}
		\label{ass: dynamics_map}
		The dynamics map $v\colon \P(S) \times C(S) \to T\P(S)$ is $\|\cdot\|_\TV$-bounded on weak-$\infty$ compact subsets of $\P(S) \times C(S)$, and is continuous with respect to the weak-$\infty$ topology on its domain and the weak topology on its codomain.
	\end{assumption}

        {\color{black}
        The first part of \autoref{ass: game} ensures that $\NE(F)$ is compact (per \autoref{prop: ne_closed} and \autoref{prop: weakly_continuous_payoff}), which is key in preventing populations from drifting off towards an unattainable equilibrium ``at infinity.'' We illustrate this point in the following example.

	\begin{example}
		\label{ex: at_infinity}
		Consider the \emph{noncompact} strategy set $S=\R$, and suppose that $F \colon \P(S) \to C(S)$ is a game that awards the highest payoffs when the strategies being played are distinct from one another. In other words, $F$ incentivizes ingenuity and individuality. In this case, players will tend towards choosing strategies distinct from the the rest of the population, resulting in the population's probability measure to diffuse downward and outward towards a ``uniform distribution'' on $\R$ with zero density everywhere. Of course, this limiting distribution is not a valid probability distribution, and hence cannot be attained by the population. Weak compactness of $\NE(F)$ prevents this type of unattainability of the game's equilibria, and this example shows that the attainability of Nash equilibria is a distinct notion from the evolution's stable behavior.
	\end{example}

	The second part of \autoref{ass: game}, as well as \autoref{ass: game_derivative} and \autoref{ass: dynamics_map}, are regularity conditions used to employ Lyapunov theory for stability analysis. Specifically, the second part of \autoref{ass: game} ensures that our Lyapunov function, should one exist, is smooth enough to have well-defined derivatives in order to quantify a non-increase in ``energy'' along trajectories of the game's evolution. On the other hand, \autoref{ass: game_derivative} and \autoref{ass: dynamics_map} are what allow us to quantify a \emph{strict} ``energy'' decrease in a well-defined manner, should such asymptotic convergence of the population occur.
        }

	We now present our main result, which shows that, when the nonautonomous portion of the dynamics is Nash stationary and $\delta$-dissipative, and when the feedback portion of the dynamics \hl{induces decreasing energy supply rates}, the interconnected closed-loop evolutionary dynamics model has\iftoggle{thesis}{ dynamically}{} stable Nash equilibria.

	\begin{theorem}[Main Result]
		\label{thm: delta-dissipative}
		Consider a game $F\colon \P(S) \to C(S)$, let $v\colon \P(S) \times C(S) \to T\P(S)$, and assume that \iftoggle{thesis}{\autoref{ass: exist_edm_solution} and \autoref{ass: game} both hold}{\autoref{ass: game} holds}. If $v$ is Nash stationary and $\delta$-dissipative with supply rate $w \colon \M(S) \times C(S) \to \R$ and it holds that
		\begin{equation}
			w\left(\nu , D\overline{F}(\mu)\nu\right) \le 0 ~ \text{for all $\mu\in\P(S)$ and all $\nu\in T\P(S)$},
			\label{eq: supply_inequality}
		\end{equation}
		then $\NE(F)$ is weakly Lyapunov stable under the EDM \eqref{eq: edm}. If, additionally, \autoref{ass: game_derivative} and \autoref{ass: dynamics_map} both hold and $v$ is strictly $\delta$-dissipative, then $\NE(F)$ is globally weakly attracting under the EDM \eqref{eq: edm}.
	\end{theorem}

	\begin{proof}
        We break the proof up into two parts: first, the proof of weak Lyapunov stability, and second, the proof of global weak attraction.
        
        \paragraph{Proof of weak Lyapunov stability\iftoggle{longform}{.}{}}
        Since $v$ is $\delta$-dissipative with supply rate $w\colon \M(S) \times C(S) \to \R$, there exist $\sigma \colon \P(S) \times C(S) \to \Rpl$ and $\storage\colon\P(S) \times C(S) \to \Rpl$ with $\storage$ having an appropriate extension $\overline{\storage} \colon U\times C(S) \to \R$ as in \autoref{def: delta-dissipative}. Define $V \colon \P(S) \to \Rpl$ by $V(\mu) = \storage(\mu,F(\mu))$. By \autoref{prop: ne_closed} and \autoref{prop: weakly_continuous_payoff}, $\NE(F)$ is weakly compact.\iftoggle{longform}{ Thus, by \autoref{lem: lyapunov}, it suffices to show that $V$ is a global Lyapunov function for $\NE(F)$ under $\mu\mapsto v(\mu,F(\mu))$ (according to \autoref{def: lyapunov_function}).}{ Thus, it suffices to show that $V$ is a global Lyapunov function for $\NE(F)$ under $\mu\mapsto v(\mu,F(\mu))$ (see Appendix B in our online technical report \citet{anderson2025dissipativity-long}).} Let $\overline{V}\colon U\cap U' \to \R$ be defined by $\overline{V}(\mu) = \overline{\storage}(\mu,\overline{F}(\mu))$. Note that $U\cap U'$ is strongly open and contains $\P(S)$, and that $\overline{V}$ is weakly continuous and Fr\'echet differentiable since $\overline{\storage}$ is weak-$\infty$-continuous and $\overline{F}$ is weakly continuous, and both $\overline{\storage}$ and $\overline{F}$ are Fr\'echet differentiable. Also note that $\overline{V}(\mu) = V(\mu)$ for all $\mu\in \P(S)$. Also, if $\mu\in \NE(F)$, then $v(\mu,F(\mu)) = 0$ by \autoref{prop: ns_implies_rest_equals_nash}, and therefore $\overline{V}(\mu) = V(\mu) = \storage(\mu,F(\mu)) = 0$ by \eqref{eq: delta-dissipative-2}. Furthermore, if $\mu\in\P(S)\setminus\NE(F)$, then again by \autoref{prop: ns_implies_rest_equals_nash} we have that $v(\mu,F(\mu)) \ne 0$, so $\overline{V}(\mu) = V(\mu) = \storage(\mu,F(\mu)) > 0$ by \eqref{eq: delta-dissipative-2}. \iftoggle{longform}{Therefore, the first two conditions from \autoref{def: lyapunov_function} on $V$ to be a global Lyapunov function for $\NE(F)$ under $\mu \mapsto v(\mu,F(\mu))$ are satisfied.}{}

		Next, since $\overline{\storage}$ and $\overline{F}$ are Fr\'echet differentiable,
		\begin{equation*}
			D\overline{V}(\mu) = \partial_1 \overline{\storage}(\mu,\overline{F}(\mu)) + \partial_2 \overline{\storage}(\mu,\overline{F}(\mu)) \circ D\overline{F}(\mu)
		\end{equation*}
		for all $\mu\in U\cap U'$, and therefore, since $\overline{F}(\mu) = F(\mu)$ for all $\mu\in\P(S)$, it holds for all $\mu\in\P(S)$ that
		\iftoggle{twocol}{
		\begin{equation}
		\begin{aligned}
			D\overline{V}(\mu)v(\mu,F(\mu)) &= \partial_1 \overline{\storage}(\mu,F(\mu)) v(\mu, F(\mu)) \\
			&\qquad + \partial_2 \overline{\storage}(\mu,F(\mu)) D\overline{F}(\mu) v(\mu, F(\mu)) \\
			&\le -\sigma(\mu,F(\mu)) \\
			&\qquad + w(v(\mu,F(\mu)),D\overline{F}(\mu)v(\mu,F(\mu))) \\
			&\le -\sigma(\mu,F(\mu)) \\
			&\le 0,
		\end{aligned}
		\label{eq: decreasing_derivative}
		\end{equation}
		}{
		\begin{equation}
		\begin{aligned}
			D\overline{V}(\mu)v(\mu,F(\mu)) &= \partial_1 \overline{\storage}(\mu,F(\mu)) v(\mu, F(\mu)) + \partial_2 \overline{\storage}(\mu,F(\mu)) D\overline{F}(\mu) v(\mu, F(\mu)) \\
			&\le -\sigma(\mu,F(\mu)) + w(v(\mu,F(\mu)),D\overline{F}(\mu)v(\mu,F(\mu))) \\
			&\le -\sigma(\mu,F(\mu)) \\
			&\le 0,
		\end{aligned}
		\label{eq: decreasing_derivative}
		\end{equation}
		}
		where the first inequality follows from \eqref{eq: delta-dissipative-1}, and the second inequality follows from \eqref{eq: supply_inequality} together with the fact that $v(\mu,F(\mu))\in T\P(S)$. Hence, $V$ is indeed a global Lyapunov function for $\NE(F)$ under $\mu \mapsto v(\mu,F(\mu))$, so $\NE(F)$ is weakly Lyapunov stable under the EDM \eqref{eq: edm}.

        \paragraph{Proof of global weak attraction\iftoggle{longform}{.}{}}
        Now suppose that the $\delta$-dissipativity of $v$ is strict and that the additional hypotheses of \autoref{ass: game_derivative} and \autoref{ass: dynamics_map} both hold.\iftoggle{longform}{ By \autoref{lem: strict_lyapunov}, it suffices to show that $V$ is a strict global Lyapunov function for $\NE(F)$ under $\mu\mapsto v(\mu,F(\mu))$ (according to \autoref{def: lyapunov_function}).}{ It suffices to show that $V$ is a strict global Lyapunov function for $\NE(F)$ under $\mu\mapsto v(\mu,F(\mu))$ (see Appendix B in our online technical report \citet{anderson2025dissipativity-long}).} This amounts to proving that $\mu \mapsto D\overline{V}(\mu)v(\mu,F(\mu))$ is weakly continuous and that $D\overline{V}(\mu)v(\mu,F(\mu)) < 0$ for all $\mu\in\P(S)\setminus\NE(F)$. Indeed, the continuity condition holds by \autoref{lem: weakly_continuous_derivative}, which we state and prove in \autoref{sec: proofs}.

		Next, if $\mu\in\P(S)\setminus\NE(F)$, then \autoref{prop: ns_implies_rest_equals_nash} gives that $v(\mu,F(\mu)) \ne 0$ so $\sigma(\mu,F(\mu)) > 0$ by \eqref{eq: delta-dissipative-3}, implying that $D\overline{V}(\mu)v(\mu,F(\mu)) < 0$ for all such $\mu$ by \eqref{eq: decreasing_derivative}. Hence, $V$ is indeed a strict global Lyapunov function for $\NE(F)$ under $\mu\mapsto v(\mu,F(\mu))$, so $\NE(F)$ is globally weakly attracting under the EDM \eqref{eq: edm}.
	\end{proof}

	\begin{lemmaE}[][all end,text link=]
		\label{lem: weakly_continuous_derivative}
		\iftoggle{longform}{The map }{}$\mu\mapsto D\overline{V}(\mu)v(\mu,F(\mu))$ is weakly continuous.
	\end{lemmaE}

	\begin{proofE}
		Since $S$ is a metric space, the weak topology on $\P(S)$ is metrizable \citep[Theorem~11.3.3]{dudley2002real}. Therefore, the weak topology on $\P(S)$ is first-countable and hence functions with domain $\P(S)$ are weakly continuous if they are weakly sequentially continuous. Thus, to prove the claim, it suffices to show that
		\begin{equation*}
			D\overline{V}(\mu_n)v(\mu_n,F(\mu_n)) \to D\overline{V}(\mu)v(\mu,F(\mu))
		\end{equation*}
		whenever $\mu_n \to \mu$ weakly. To this end, let $\{\mu_n \in \P(S) : n\in\N\}$ be a sequence that converges weakly to $\mu\in\P(S)$. Then we have that
		\iftoggle{twocol}{
		\begin{equation}
		\begin{aligned}
			&|D\overline{V}(\mu_n)v(\mu_n,F(\mu_n)) - D\overline{V}(\mu)v(\mu,F(\mu))| \\
			&\quad = |D\overline{V}(\mu_n)v(\mu_n,\overline{F}(\mu_n)) - D\overline{V}(\mu)v(\mu,\overline{F}(\mu))| \\
												  &\quad \le |D\overline{V}(\mu)(v(\mu_n,\overline{F}(\mu_n)) - v(\mu,\overline{F}(\mu)))| \\
												  &\quad \quad+|(D\overline{V}(\mu_n) - D\overline{V}(\mu))v(\mu,\overline{F}(\mu))| \\
												  &\quad \quad+|(D\overline{V}(\mu_n) - D\overline{V}(\mu))(v(\mu_n,\overline{F}(\mu_n)) - v(\mu,\overline{F}(\mu)))|
		\end{aligned}
		\label{eq: weak_continuity_bound}
		\end{equation}
		}{
		\begin{equation}
		\begin{aligned}
			|D\overline{V}(\mu_n)v(\mu_n,F(\mu_n)) - D\overline{V}(\mu)v(\mu,F(\mu))| &= |D\overline{V}(\mu_n)v(\mu_n,\overline{F}(\mu_n)) - D\overline{V}(\mu)v(\mu,\overline{F}(\mu))| \\
												  &\le |D\overline{V}(\mu)(v(\mu_n,\overline{F}(\mu_n)) - v(\mu,\overline{F}(\mu)))| \\
												  &\qquad+|(D\overline{V}(\mu_n) - D\overline{V}(\mu))v(\mu,\overline{F}(\mu))| \\
												  &\qquad+|(D\overline{V}(\mu_n) - D\overline{V}(\mu))(v(\mu_n,\overline{F}(\mu_n)) \iftoggle{thesis}{ \\ &\qquad\qquad }{}- v(\mu,\overline{F}(\mu)))|.
		\end{aligned}
		\label{eq: weak_continuity_bound}
		\end{equation}
		}
		Assume for the time being that every Fr\'echet derivative $D\overline{V}(\nu)$ is weakly continuous, and that $D\overline{V}$ is weakly continuous on $\P(S)$. Then, under this assumption, it holds that
		\begin{equation*}
			|D\overline{V}(\mu)(v(\mu_n,\overline{F}(\mu_n)) - v(\mu,\overline{F}(\mu)))| \to 0,
		\end{equation*}
		since $v(\mu_n,\overline{F}(\mu_n)) \to v(\mu,\overline{F}(\mu))$ weakly, as $\overline{F}$ is weakly continuous and $v$ is continuous with respect to the weak-$\infty$ topology on its domain and the weak topology on its codomain. Furthermore,
		\iftoggle{twocol}{
		\begin{align*}
			&|(D\overline{V}(\mu_n) - D\overline{V}(\mu))v(\mu,\overline{F}(\mu))| \\
			&\qquad \le \|D\overline{V}(\mu_n) - D\overline{V}(\mu)\|_{\M(S)^*} \|v(\mu,\overline{F}(\mu))\|_\TV \to 0,
		\end{align*}
		}{
		\begin{equation*}
			|(D\overline{V}(\mu_n) - D\overline{V}(\mu))v(\mu,\overline{F}(\mu))| \le \|D\overline{V}(\mu_n) - D\overline{V}(\mu)\|_{\M(S)^*} \|v(\mu,\overline{F}(\mu))\|_\TV \to 0,
		\end{equation*}
		}
		since $D\overline{V}(\mu_n) \to D\overline{V}(\mu)$ in the dual space $\M(S)^*$ with associated operator norm $\|\cdot\|_{\M(S)^*}$ induced by the total variation norm on $\M(S)$, as $D\overline{V} \colon U\cap U' \to \M(S)^*$ is weakly continuous on $\P(S)$ by our above assumption. Finally,
		\iftoggle{twocol}{
		\begin{multline*}
			|(D\overline{V}(\mu_n) - D\overline{V}(\mu))(v(\mu_n,\overline{F}(\mu_n)) - v(\mu,\overline{F}(\mu)))| \\
			\begin{aligned}
			& \le \|D\overline{V}(\mu_n) - D\overline{V}(\mu)\|_{\M(S)^*} \|v(\mu_n,\overline{F}(\mu_n)) \\
			&\qquad - v(\mu,\overline{F}(\mu))\|_\TV \\
			& \le 2\|D\overline{V}(\mu_n) - D\overline{V}(\mu)\|_{\M(S)^*} \sup_{\nu\in\P(S)}\|v(\nu,\overline{F}(\nu))\|_\TV \\
			& \to 0,
			\end{aligned}
		\end{multline*}
		}{
		\iftoggle{thesis}{
		\begin{align*}
			& |(D\overline{V}(\mu_n) - D\overline{V}(\mu))(v(\mu_n,\overline{F}(\mu_n)) - v(\mu,\overline{F}(\mu)))| \\
			&\qquad\le \|D\overline{V}(\mu_n) - D\overline{V}(\mu)\|_{\M(S)^*} \|v(\mu_n,\overline{F}(\mu_n)) - v(\mu,\overline{F}(\mu))\|_\TV \\
			& \qquad \le 2\|D\overline{V}(\mu_n) - D\overline{V}(\mu)\|_{\M(S)^*} \sup_{\nu\in\P(S)}\|v(\nu,\overline{F}(\nu))\|_\TV \\
			& \qquad \to 0,
		\end{align*}
		}{
		\begin{align*}
			& |(D\overline{V}(\mu_n) - D\overline{V}(\mu))(v(\mu_n,\overline{F}(\mu_n)) - v(\mu,\overline{F}(\mu)))| \\
			&\qquad \le \|D\overline{V}(\mu_n) - D\overline{V}(\mu)\|_{\M(S)^*} \|v(\mu_n,\overline{F}(\mu_n)) - v(\mu,\overline{F}(\mu))\|_\TV \\
			&\qquad \le 2\|D\overline{V}(\mu_n) - D\overline{V}(\mu)\|_{\M(S)^*} \sup_{\nu\in\P(S)}\|v(\nu,\overline{F}(\nu))\|_\TV \\
			&\qquad \to 0,
		\end{align*}
		}
		}
		since again $D\overline{V}(\mu_n) \to D\overline{V}(\mu)$ in $\M(S)^*$ by the weak continuity assumption on $D\overline{V}$, and \iftoggle{thesis}{$\sup_{\nu\in\P(S)}\|v(\nu,\overline{F}(\nu))\|_\TV \le \sup_{(\nu,g)\in \P(S) \times \overline{F}(\P(S))} \|v(\nu,g)\|_\TV \le M$}{\iftoggle{twocol}{$\sup_{\nu\in\P(S)}\|v(\nu,\overline{F}(\nu))\|_\TV \le \sup_{(\nu,g)\in \P(S) \times \overline{F}(\P(S))} \|v(\nu,g)\|_\TV \le M$}{$$\sup_{\nu\in\P(S)}\|v(\nu,\overline{F}(\nu))\|_\TV \le \sup_{(\nu,g)\in \P(S) \times \overline{F}(\P(S))} \|v(\nu,g)\|_\TV \le M$$}} for some finite $M\in[0,\infty)$ by the $\|\cdot\|_\TV$-boundedness of $v$ on weak-$\infty$ compact subsets of $\P(S) \times C(S)$. Therefore, under the above assumptions, it must be that
		\begin{equation*}
			D\overline{V}(\mu_n)v(\mu_n,F(\mu_n)) \to D\overline{V}(\mu)v(\mu,F(\mu)),
		\end{equation*}
		which is what was to be proven. Thus, it remains to prove the above assumptions, namely, that every Fr\'echet derivative $D\overline{V}(\nu)$ is weakly continuous, and that $D\overline{V}$ is weakly continuous on $\P(S)$.

		Let us first prove that $D\overline{V}(\mu)\colon \M(S) \to \R$ is weakly continuous for all $\mu\in U\cap U'$. Let $\mu\in U\cap U'$. Since $\partial_2 \overline{\storage}(\mu,\rho) \colon C(S) \to \R$ is continuous (with respect to the topology on $C(S)$ induced by $\|\cdot\|_\infty$) for all $\rho\in C(S)$ by definition of the Fr\'echet derivative, and since $D\overline{F}(\mu)\colon \M(S) \to C(S)$ is weakly continuous under the hypotheses of the theorem, it holds that the composition $\partial_2 \overline{\storage}(\mu,\overline{F}(\mu)) \circ D\overline{F}(\mu) \colon \M(S) \to \R$ is weakly continuous. Since $\partial_1 \overline{\storage}(\mu,\overline{F}(\mu)) \colon \M(S) \to \R$ is also weakly continuous under the hypotheses of the theorem, we conclude that
		\begin{equation*}
			D\overline{V}(\mu) = \partial_1 \overline{\storage}(\mu,\overline{F}(\mu)) + \partial_2 \overline{\storage}(\mu,\overline{F}(\mu)) \circ D\overline{F}(\mu)
		\end{equation*}
		is weakly continuous, proving the first assumption.

		Finally, let us prove the remaining assumption, namely, that $D\overline{V} \colon U\cap U' \to \M(S)^*$ is weakly continuous on $\P(S)$ (that is, continuous with respect to the weak topology on its domain $U\cap U' \subseteq \M(S)$ and the topology on its codomain $\M(S)^*$ induced by the operator norm $\|\cdot\|_{\M(S)^*}$). Once again, since we are considering weak continuity of a function on $\P(S)$, where the weak topology is first-countable, it suffices to prove weak sequential continuity. Let $\{\mu_n \in \P(S) : n\in\N\}$ be a sequence that converges weakly to $\mu\in\P(S)$. Then
		\iftoggle{twocol}{
		\begin{align*}
			&\|D\overline{V}(\mu_n) - D\overline{V}(\mu)\|_{\M(S)^*} \\
			&\qquad \le \left\lVert \partial_1 \overline{\storage}(\mu_n,\overline{F}(\mu_n)) - \partial_1 \overline{\storage}(\mu,\overline{F}(\mu)) \right\rVert_{\M(S)^*} \\
										&\qquad+ \lVert \partial_2 \overline{\storage}(\mu_n,\overline{F}(\mu_n)) \circ D\overline{F}(\mu_n) \\
										&\qquad \qquad - \partial_2 \overline{\storage}(\mu,\overline{F}(\mu)) \circ D\overline{F}(\mu) \rVert_{\M(S)^*}.
		\end{align*}
		}{
		\begin{align*}
			\|D\overline{V}(\mu_n) - D\overline{V}(\mu)\|_{\M(S)^*} &\le \left\lVert \partial_1 \overline{\storage}(\mu_n,\overline{F}(\mu_n)) - \partial_1 \overline{\storage}(\mu,\overline{F}(\mu)) \right\rVert_{\M(S)^*} \\
										&\iftoggle{thesis}{\quad}{\qquad}+ \left\lVert \partial_2 \overline{\storage}(\mu_n,\overline{F}(\mu_n)) \circ D\overline{F}(\mu_n) - \partial_2 \overline{\storage}(\mu,\overline{F}(\mu)) \circ D\overline{F}(\mu) \right\rVert_{\M(S)^*}.
		\end{align*}
		}
		It is clear that the first term in the above upper bound converges to $0$ due to the weak-$\infty$ continuity of $(\nu,\rho)\mapsto \partial_1 \overline{\storage}(\nu,\rho)$ together with the weak continuity of $\overline{F}$. Further upper-bounding the second term in a similar manner to the bound \eqref{eq: weak_continuity_bound} and appealing to the finiteness of $\|\varphi\|_{\TV}$ and $\|\psi\|_{\M(S)^*}$ for $\varphi\in C(S)^* = \M(S)$ and $\psi\in\M(S)^*$ together with the weak continuity of $\overline{F}$ and $D\overline{F}$ as well as the weak-$\infty$ continuity of $(\nu,\rho) \mapsto \partial_2 \overline{\storage}(\nu,\rho)$ yields that the second term converges to $0$ as well. Thus, $D\overline{V}(\mu_n) \to D\overline{V}(\mu)$ in $\M(S)^*$, so $D\overline{V}$ is indeed weakly continuous on $\P(S)$.
	\end{proofE}

	It is easy to see that \autoref{thm: delta-dissipative} generalizes the first main result in \citet{arcak2021dissipativity}, i.e., our \autoref{thm: delta-dissipative} recovers Theorem~1 in \citet{arcak2021dissipativity} when $S$ is finite. We will see in \autoref{sec: monotone_games} that \autoref{thm: delta-dissipative} also recovers other recent stability results for special types of games.


	\subsection{Specialization to Monotone Games}
	\label{sec: monotone_games}

	In this section, we consider the special class of ``monotone games,'' which are sometimes also referred to as ``stable games,'' ``contractive games,'' and ``negative semidefinite games'' in the literature. For a thorough analysis of monotone games over finite strategy sets, see \citet{hofbauer2009stable,fox2013population,park2019population}, and for works considering monotone games with an infinite number of strategies, see \citet{hofbauer2009brown,cheung2014pairwise,lahkar2015logit,lahkar2022generalized}. \hl{The latter works are all restricted to special types of dynamics, e.g., BNN, pairwise comparison, logit, and perturbed best response dynamics. In contrast, our stability result (\autoref{thm: monotone-passive}) for monotone games derived in this section holds more broadly for the class of $\delta$-passive dynamics (see \autoref{def: delta-passive} to come), which constitutes a property that may be verified modularly for various instances of particular dynamics.}

	\begin{definition}
		\label{def: monotone_game}
		A game $F\colon \P(S) \to C(S)$ is \emph{monotone} if
		\begin{equation}
			\left<F(\mu) - F(\nu) , \mu - \nu\right> \le 0
			\label{eq: monotone_game}
		\end{equation}
		for all $\mu,\nu\in\P(S)$. If, additionally, the inequality \eqref{eq: monotone_game} holds strictly for all $\mu,\nu\in\P(S)$ such that $\mu\ne \nu$, then $F$ is \emph{strictly monotone}.
	\end{definition}

	Many games in practice are monotone, e.g., random matching in two-player symmetric zero-sum games \citep[Example~4]{cheung2014pairwise}, contests \citep[Example~5]{hofbauer2009brown}, and the war of attrition \citep[Example~2.4]{hofbauer2009stable}{\color{black}, which we discussed in \autoref{ex: war_of_attrition}}. {\color{black}In general, monotonicity of games facilitates convergence of evolutionary dynamics towards Nash equilibria. Intuitively, this is because monotone games can be thought of as having payoffs that always guide the population in the direction of increasing the average payoff. This intuition is made precise in Kachurovskii's theorem \citep[Chapter~II.7,~Proposition~7.4]{showalter1996monotone}, which establishes a correspondence between (increasing) monotone operators and gradients of convex functions, the latter of which are well-known to direct flows toward global minima. With this in mind, a monotone game can be thought of as a generalization of (the negative of) a Fr\'echet derivative of a convex ``energy'' functional.}
    
    \iftoggle{longform}{\hl{ Characterizations of the equilibria of monotone games are given in \autoref{sec: supplementary}, e.g., the convexity of $\NE(F)$.}}{} \iftoggle{thesis}{The following notion is closely related to that of monotonicity, as we will see in \autoref{lem: monotone_sde}, and will serve as the link between monotonicity and the inequality \eqref{eq: supply_inequality}.

	\begin{definition}
		\label{def: sde}
		A game $F\colon \P(S) \to C(S)$ that extends to a continuously Fr\'echet differentiable map $\overline{F} \colon U' \to C(S)$ defined on a strongly open set $U'\subseteq \M(S)$ containing $\P(S)$ is said to \emph{satisfy self-defeating externalities} if
		\begin{equation*}
			\left<D\overline{F}(\mu)\nu,\nu\right> \le 0 ~ \text{for all $\mu\in\P(S)$ and all $\nu\in T\P(S)$}.
		\end{equation*}
	\end{definition}

	\begin{lemma}[{\citealp[Lemma~3]{cheung2014pairwise}}]
		\label{lem: monotone_sde}
		Consider a game $F\colon \P(S) \to C(S)$ that extends to a continuously Fr\'echet differentiable map $\overline{F} \colon U' \to C(S)$ defined on a strongly open set $U'\subseteq \M(S)$ containing $\P(S)$. It holds that $F$ is monotone if and only if $F$ satisfies self-defeating externalities.
	\end{lemma}}{\hl{The following lemma serves as the link between monotonicity and the inequality \eqref{eq: supply_inequality}, and is key in proving \autoref{thm: monotone-passive} to come.

	\begin{lemma}[{\citealp[Lemma~3]{cheung2014pairwise}}]
		\label{lem: monotone_sde}
		Consider a game $F\colon \P(S) \to C(S)$ that extends to a continuously Fr\'echet differentiable map $\overline{F} \colon U' \to C(S)$ defined on a strongly open set $U'\subseteq \M(S)$ containing $\P(S)$. It holds that $F$ is monotone if and only if
		\begin{equation}
			\left<D\overline{F}(\mu)\nu,\nu\right> \le 0 ~ \text{for all $\mu\in\P(S)$ and all $\nu\in T\P(S)$}.
			\label{eq: sde}
		\end{equation}
	\end{lemma}}}

	We now show that our general dissipativity theory can be applied to monotone games to recover recent stability results in the literature. We start with the following specialization of $\delta$-dissipativity, which generalizes the notion of $\delta$-passivity introduced in \citet{fox2013population} for the finite-strategy setting.

	\begin{definition}
		\label{def: delta-passive}
		A map $v\colon\P(S)\times C(S) \to \M(S)$ is \emph{$\delta$-passive} if it is $\delta$-dissipative with supply rate $w \colon (\mu,\eta) \mapsto \left<\eta,\mu\right> = \int_S \eta d\mu$. The map $v$ is \emph{strictly $\delta$-passive} if it is strictly $\delta$-dissipative with such supply rate $w$.
	\end{definition}

	As is the case with $\delta$-dissipativity, we see that $\delta$-passivity is solely a property of the nonautonomous portion of the evolutionary dynamics defined by the dynamics map $v$. We remark that $\delta$-passivity is common in practice. {\color{black}For instance, our two earlier examples of the BNN and impartial pairwise comparison dynamics are both strictly $\delta$-passive, as the following result shows.}

	\iftoggle{longform}{
	\begin{propositionE}[][end,restate,text link=]
	}{
	\begin{propositionE}[][end,text link=]
	}
		\label{prop: bnn_pairwise_are_delta-passive}
		If $v\colon \P(S) \times C(S) \to \M(S)$ is the dynamics map for either the BNN dynamics of \autoref{ex: bnn} or the impartial pairwise comparison dynamics of \autoref{ex: pcd}, then $v$ is strictly $\delta$-passive.
	\end{propositionE}

	\iftoggle{longform}{
	\begin{proofE}
		We prove the result for the two dynamics separately.

		\paragraph{BNN dynamics\iftoggle{longform}{.}{}}
		Consider the BNN dynamics of \autoref{ex: bnn}. We have that
		\begin{equation*}
			v(\mu,\rho)(B) = \int_B \max\{0,\left<\rho,\delta_s\right> - \left<\rho,\mu\right>\} d\lambda(s) - \mu(B) \int_S \max\{0,\left<\rho,\delta_s\right> - \left<\rho,\mu\right>\} d\lambda(s)
		\end{equation*}
		for all $\mu\in\P(S)$, all $\rho\in C(S)$, and all $B\in\B(S)$. Define $\overline{\storage} \colon \M(S) \times C(S) \to \R$ and $\sigma \colon \P(S) \times C(S) \to \R$ by
		\begin{align*}
			\overline{\storage}(\mu,\rho) &= \frac{1}{2}\int_S \max\{0,\left<\rho,\delta_s\right> - \left<\rho,\mu\right>\}^2 d\lambda(s), \\
			\sigma(\mu,\rho) &= \left<\rho,v(\mu,\rho)\right> \int_S \max\{0,\left<\rho,\delta_s\right> - \left<\rho,\mu\right>\} d\lambda(s).
		\end{align*}
		Notice that $\overline{\storage}(\mu,\rho)$ and $\sigma(\mu,\rho)$ are finite for all $\mu\in\M(S)$ and all $\rho\in C(S)$, since $s\mapsto \max\{0,\left<\rho,\delta_s\right> - \left<\rho,\mu\right>\}$ and $s\mapsto \max\{0,\left<\rho,\delta_s\right> - \left<\rho,\mu\right>\}^2$ are continuous and $S$ is compact. Also notice that $\overline{\storage}(\mu,\rho) \ge 0$ for all $\mu\in\M(S)$ and all $\rho\in C(S)$. Thus, we may define $\storage \colon \P(S) \times C(S) \to \Rpl$ by the restriction of $\overline{\storage}$ to the domain $\P(S)\times C(S) \subseteq \M(S) \times C(S)$. We claim that $\sigma$ and $\storage$ are appropriate maps to prove the strict $\delta$-passivity of $v$.

		To this end, first note that $\M(S)$ is strongly open, $\overline{\storage}$ is weak-$\infty$-continuous, $\overline{\storage}$ is Fr\'echet differentiable, $(\mu,\rho)\mapsto \partial_1 \overline{\storage}(\mu,\rho)$ and $(\mu,\rho)\mapsto \partial_2 \overline{\storage}(\mu,\rho)$ are weak-$\infty$-continuous, and every partial Fr\'echet derivative $\partial_1 \overline{\storage}(\mu,\rho)$ is weakly continuous. All that remains to prove are \eqref{eq: delta-dissipative-1} with $w\colon (\mu,\eta)\mapsto \left<\eta,\mu\right>$, \eqref{eq: delta-dissipative-2}, \eqref{eq: delta-dissipative-3}, and that $\sigma \ge 0$.

		Let $\mu\in\P(S)$ and $\rho\in C(S)$. It holds that $\storage(\mu,\rho) = 0$ if and only if
		\begin{equation}
			\int_S \max\{0,\left<\rho,\delta_s\right> - \left<\rho,\mu\right>\}^2 d\lambda(s) = 0.
			\label{eq: zero_storage_bnn}
		\end{equation}
		Since $s\mapsto \max\{0,\left<\rho,\delta_s\right> - \left<\rho,\mu\right>\}^2$ is a continuous real-valued function on $S$, the preimage $U\coloneqq \{s\in S : \max\{0,\left<\rho,\delta_s\right> - \left<\rho,\mu\right>\}^2 > 0\}$ is open. Therefore, if $U$ is nonempty, it contains some $s'\in S$, and hence since $\lambda$ has full support, $s'$ must be an element of $\supp(\lambda)$, implying that $\lambda(U) > 0$. This in turn would imply that $\int_U \max\{0,\left<\rho,\delta_s\right> - \left<\rho,\mu\right>\}^2 d\lambda(s) > 0$ as the Lebesgue integral of a positive function over a set of positive measure is positive. However, this would contradict \eqref{eq: zero_storage_bnn}. Thus, $\storage(\mu,\rho) = 0$ if and only if
		\begin{equation*}
			\max\{0,\left<\rho,\delta_s\right> - \left<\rho,\mu\right>\}^2 = 0 ~ \text{for all $s\in S$},
		\end{equation*}
		which holds if and only if
		\begin{equation}
			\rho(s) \le \left<\rho,\mu\right> ~ \text{for all $s\in S$}.
			\label{eq: bnn_passivity_nash}
		\end{equation}
		It is clear that, if $\left<\rho,\nu\right> \le \left<\rho,\mu\right>$ for all $\nu\in \P(S)$, then \eqref{eq: bnn_passivity_nash} holds. Conversely, if \eqref{eq: bnn_passivity_nash} holds, then $\left<\rho,\nu\right> = \int_S \rho(s) d\nu(s) \le \int_S \left<\rho,\mu\right>d\nu(s) = \left<\rho,\mu\right>$ for all $\nu\in\P(S)$, and thus by Nash stationarity of $v$ (\autoref{prop: bnn_pcd_satisfy_ns}) we conclude that $\storage(\mu,\rho) = 0$ if and only if
		\begin{equation*}
			v(\mu,\rho) = 0,
		\end{equation*}
		which proves \eqref{eq: delta-dissipative-2}.

		Again let $\mu\in\P(S)$ and $\rho\in C(S)$. If $v(\mu,\rho) = 0$, then certainly $\sigma(\mu,\rho) = 0$ due to linearity of $\left<\rho,\cdot\right>$. Notice that
		\iftoggle{thesis}{
		\begin{align*}
			\left<\rho,v(\mu,\rho)\right> &= \int_S \rho(s') d(v(\mu,\rho))(s') \\
			&= \int_S \rho(s') \max\{0,\left<\rho,\delta_{s'}\right> - \left<\rho,\mu\right>\} d\lambda(s') \\
			&\qquad - \int_S \max\{0,\left<\rho,\delta_s\right> - \left<\rho,\mu\right>\} d\lambda(s) \int_S \rho(s') d\mu(s') \\
			&= \int_S \left(\rho(s') - \int_S \rho(\tilde{s}) d\mu(\tilde{s})\right)\max\{0,\left<\rho,\delta_{s'}\right> - \left<\rho,\mu\right>\} d\lambda(s') \\
			&= \int_S (\left<\rho,\delta_{s'}\right> - \left<\rho,\mu\right>)\max\{0,\left<\rho,\delta_{s'}\right> - \left<\rho,\mu\right>\} d\lambda(s').
		\end{align*}
		}{
		\begin{align*}
			\left<\rho,v(\mu,\rho)\right> &= \int_S \rho(s') d(v(\mu,\rho))(s') \\
			&= \int_S \rho(s') \max\{0,\left<\rho,\delta_{s'}\right> - \left<\rho,\mu\right>\} d\lambda(s') - \int_S \max\{0,\left<\rho,\delta_s\right> - \left<\rho,\mu\right>\} d\lambda(s) \int_S \rho(s') d\mu(s') \\
			&= \int_S \left(\rho(s') - \int_S \rho(\tilde{s}) d\mu(\tilde{s})\right)\max\{0,\left<\rho,\delta_{s'}\right> - \left<\rho,\mu\right>\} d\lambda(s') \\
			&= \int_S (\left<\rho,\delta_{s'}\right> - \left<\rho,\mu\right>)\max\{0,\left<\rho,\delta_{s'}\right> - \left<\rho,\mu\right>\} d\lambda(s').
		\end{align*}
		}
		Notice that $(\left<\rho,\delta_{s'}\right> - \left<\rho,\mu\right>)\max\{0,\left<\rho,\delta_{s'}\right> - \left<\rho,\mu\right>\} \ge 0$ for all $s'\in S$ and hence $\left<\rho,v(\mu,\rho)\right> \ge 0$. Furthermore, by the usual arguments based on continuity and nonnegativity of the integrand together with full support of $\lambda$, we see that
		\begin{equation*}
			\left<\rho,v(\mu,\rho)\right> = 0
		\end{equation*}
		if and only if
		\begin{equation*}
			\rho(s') = \left<\rho,\delta_{s'}\right> \le \left<\rho,\mu\right> ~ \text{for all $s'\in S$},
		\end{equation*}
		which, as shown above, holds true if and only if $v(\mu,\rho) = 0$. Furthermore, notice that by the same arguments,
		\begin{equation*}
			\int_S \max\{0,\left<\rho,\delta_s\right> - \left<\rho,\mu\right>\} d\lambda(s) \ge 0,
		\end{equation*}
		with equality holding if and only if $v(\mu,\rho) = 0$. Thus,
		\begin{equation*}
			\sigma(\mu,\rho) = \left<\rho,v(\mu,\rho)\right>\int_S \max\{0,\left<\rho,\delta_s\right> - \left<\rho,\mu\right>\} d\lambda(s) \ge 0,
		\end{equation*}
		with equality holding if and only if $v(\mu,\rho) = 0$. This proves \eqref{eq: delta-dissipative-3}.

		All that remains to be proven is \eqref{eq: delta-dissipative-1} with $w\colon (\mu,\eta) \mapsto \left<\eta,\mu\right>$. Let $\mu\in\P(S)$, $\rho \in C(S)$, and $\eta\in C(S)$. Define $\tau \colon \R\to \Rpl$ by $\tau(r) = \max\{0,r\}^2$, so that $\tau'(r) = 2\max\{0,r\}$ and $\overline{\storage}(\mu,\rho) = \frac{1}{2}\int_S \tau(\left<\rho,\delta_s\right> - \left<\rho,\mu\right>) d\lambda(s)$. Computing the first partial Fr\'echet derivative of $\overline{\storage}$ using the chain rule yields that
		\begin{align*}
			\partial_1 \overline{\storage}(\mu,\rho)v(\mu,\rho) &= \frac{1}{2}\int_S \tau'(\left<\rho,\delta_s\right> - \left<\rho,\mu\right>)(-\left<\rho,v(\mu,\rho)\right>) d\lambda(s) \\
			&= -\left<\rho,v(\mu,\rho)\right>\int_S \max\{0,\left<\rho,\delta_s\right> - \left<\rho,\mu\right>\}d\lambda(s) \\
			&= -\sigma(\mu,\rho).
		\end{align*}
		Computing the second partial Fr\'echet derivative of $\overline{\storage}$ using the chain rule yields that
		\begin{align*}
			\partial_2 \overline{\storage}(\mu,\rho) \eta &= \frac{1}{2}\int_S \tau'(\left<\rho,\delta_s\right> - \left<\rho,\mu\right>)(\left<\eta,\delta_s\right> - \left<\eta,\mu\right>) d\lambda(s) \\
			&= \int_S (\left<\eta,\delta_s\right> - \left<\eta,\mu\right>)\max\{0,\left<\rho,\delta_s\right> - \left<\rho,\mu\right>\} d\lambda(s) \\
			&= \int_S \left(\eta(s) - \int_S \eta(\tilde{s}) d\mu(\tilde{s})\right)\max\{0,\left<\rho,\delta_s\right> - \left<\rho,\mu\right>\} d\lambda(s) \\
			&= \int_S \eta(s) \max\{0,\left<\rho,\delta_s\right> - \left<\rho,\mu\right>\}d\lambda(s) \iftoggle{thesis}{\\ & \qquad}{} - \int_S \eta(\tilde{s})d\mu(\tilde{s}) \int_S \max\{0,\left<\rho,\delta_s\right> - \left<\rho,\mu\right>\} d\lambda(s) \\
			&= \int_S \eta(s) d(v(\mu,\rho))(s) \\
			&= \left<\eta,v(\mu,\rho)\right> \\
			&= w(v(\mu,\rho),\eta).
		\end{align*}
		Thus, altogether we find that
		\begin{equation*}
			\partial_1 \overline{\storage}(\mu,\rho)v(\mu,\rho) + \partial_2 \overline{\storage}(\mu,\rho) \eta = -\sigma(\mu,\rho) + w(v(\mu,\rho),\eta),
		\end{equation*}
		which shows that \eqref{eq: delta-dissipative-1} holds and hence concludes the proof for the BNN dynamics.

		\paragraph{Impartial pairwise comparison dynamics\iftoggle{longform}{.}{}}
		Consider the impartial pairwise comparison dynamics of \autoref{ex: pcd}. We have that
		\begin{equation*}
			v(\mu,\rho)(B) = \int_B \int_S \gamma(s,s',\rho) d\mu(s) d\lambda(s') - \int_B \int_S \gamma(s',s,\rho) d\lambda(s) d\mu(s')
		\end{equation*}
		for all $\mu\in\P(S)$, all $\rho\in C(S)$, and all $B\in\B(S)$. Since the pairwise comparison dynamics under consideration are impartial, it holds that for all $s'\in S$, there exists some continuous function $\phi_{s'} \colon \R\to\Rpl$ such that
		\begin{equation*}
			\gamma(s,s',\rho) = \phi_{s'}(\rho(s') - \rho(s))
		\end{equation*}
		for all $s\in S$ and all $\rho\in C(S)$. For all $s'\in S$, define $\tau_{s'}\colon \R \to \Rpl$ by
		\begin{equation*}
			\tau_{s'}(r) = \int_{[0,r]} \phi_{s'}(u) du,
		\end{equation*}
		where we see that $\tau_{s'}(r) = 0$ whenever $r < 0$, since we take $[0,r] = \emptyset$ in such cases. Notice that $\tau_{s'}$ is strictly increasing on $[0,\infty)$ since $\phi_{s'}(u) > 0$ for all $u>0$: let $u>0$, let $s,s'\in S$ be such that $s\ne s'$, and let $\rho\in C(S)$ be such that $\rho(s')-\rho(s) = u > 0$ (which exists by Urysohn's lemma and the fact that $S$ is a metric space and hence normal), so that, by sign-preservation, we have that $\sign(\phi_{s'}(u)) = \sign(\phi_{s'}(\rho(s')-\rho(s))) = \sign(\gamma(s,s',\rho)) = \sign(\max\{0,\rho(s')-\rho(s)\}) = 1$. Define $\overline{\storage} \colon \M(S) \times C(S) \to \R$ and $\sigma \colon \P(S) \times C(S) \to \R$ by
		\begin{align*}
			\overline{\storage}(\mu,\rho) &= \int_S \int_S \tau_s(\rho(s) - \rho(s')) d\lambda(s) d\mu(s'), \\
			\sigma(\mu,\rho) &= -\overline{\storage}(v(\mu,\rho),\rho).
		\end{align*}
		Notice that $\overline{\storage}(\mu,\rho)$ is finite for all $\mu\in\M(S)$ and all $\rho\in C(S)$ since $(s,s')\mapsto \tau_s(\rho(s)-\rho(s'))$ is continuous and $S$ is compact. Also notice that $\overline{\storage}(\mu,\rho) \ge 0$ for all $\mu\in\P(S)$ and all $\rho\in C(S)$. Thus, we may define $\storage \colon \P(S) \times C(S) \to \Rpl$ by the restriction of $\overline{\storage}$ to the domain $\P(S) \times C(S) \subseteq \M(S) \times C(S)$. We claim that $\sigma$ and $\storage$ are appropriate maps to prove the strict $\delta$-passivity of $v$.

		To this end, first note that $\M(S)$ is strongly open, $\overline{\storage}$ is weak-$\infty$-continuous, $\overline{\storage}$ is Fr\'echet differentiable, $(\mu,\rho)\mapsto \partial_1 \overline{\storage}(\mu,\rho)$ and $(\mu,\rho)\mapsto \partial_2 \overline{\storage}(\mu,\rho)$ are weak-$\infty$-continuous, and every partial Fr\'echet derivative $\partial_1 \overline{\storage}(\mu,\rho)$ is weakly continuous. All that remains to prove are \eqref{eq: delta-dissipative-1} with $w\colon (\mu,\eta)\mapsto \left<\eta,\mu\right>$, \eqref{eq: delta-dissipative-2}, \eqref{eq: delta-dissipative-3}, and that $\sigma \ge 0$.

		Let $\mu\in\P(S)$ and $\rho\in C(S)$. It holds that $\storage(\mu,\rho) = 0$ if and only if
		\begin{equation*}
			\int_S \int_S \tau_s(\rho(s) - \rho(s')) d\lambda(s) d\mu(s') = 0,
		\end{equation*}
		which holds if and only if
		\begin{equation*}
			\tau_s(\rho(s) - \rho(s')) = 0 ~ \text{for all $s\in S$ and all $s'\in\supp(\mu)$},
		\end{equation*}
		since $\lambda$ has full support, $s\mapsto \tau_s(\rho(s) - \rho(s'))$ is nonnegative and continuous for all $s' \in S$, and $s' \mapsto \int_S \tau_s(\rho(s) - \rho(s')) d\lambda(s)$ is nonnegative and continuous (which follows from compactness of $S$ together with the dominated convergence theorem). Since, for all $s\in S$, it holds that $\tau_s$ is strictly increasing on $[0,\infty)$ and $\tau_s(0) = 0$, it must be that $\storage(\mu,\rho) = 0$ if and only if
		\begin{equation*}
			\rho(s) \le \rho(s') ~ \text{for all $s\in S$ and all $s'\in \supp(\mu)$}.
		\end{equation*}
		Therefore, by \autoref{lem: equivalent_nash-general}, it holds that $\left<\rho,\nu\right> \le \left<\rho,\mu\right>$ for all $\nu\in\P(S)$. Hence, by Nash stationarity of $v$ (\autoref{prop: bnn_pcd_satisfy_ns}), it holds that $\storage(\mu,\rho) = 0$ if and only if
		\begin{equation*}
			v(\mu,\rho) = 0,
		\end{equation*}
		which proves \eqref{eq: delta-dissipative-2}.

		Again let $\mu\in\P(S)$ and $\rho\in C(S)$. If $v(\mu,\rho) = 0$, then\iftoggle{thesis}{}{ certainly} $\sigma(\mu,\rho) = -\overline{\storage}(v(\mu,\rho),\rho) = 0$ due to linearity of $\overline{\storage}(\cdot,\rho)$. Writing out $\sigma(\mu,\rho)$, we find that
		\begin{align*}
			\sigma(\mu,\rho) &= -\int_S\int_S \tau_s(\rho(s) - \rho(s')) d\lambda(s) d(v(\mu,\rho))(s') \\
					 &= -\int_S \left(\int_S \tau_s(\rho(s) - \rho(s')) d\lambda(s)\right)\left(\int_S \gamma(s,s',\rho) d\mu(s)\right) d\lambda(s') \\
					 &\qquad + \int_S \left(\int_S \tau_s(\rho(s) - \rho(s')) d\lambda(s)\right)\left(\int_S \gamma(s',s,\rho) d\lambda(s)\right) d\mu(s') \\
					 &= \int_S\int_S \gamma(s',s,\rho) \int_S \left(\tau_{\tilde{s}}(\rho(\tilde{s}) - \rho(s')) - \tau_{\tilde{s}}(\rho(\tilde{s}) - \rho(s))\right) d\lambda(\tilde{s}) d\lambda(s) d\mu(s').
		\end{align*}
		For all $s,s'\in S$ such that $\rho(s) \le \rho(s')$, it holds by sign-preservation that $\sign(\gamma(s',s,\rho)) = \sign(\max\{0,\rho(s)-\rho(s')\}) = 0$, and therefore $\gamma(s',s,\rho) = 0$ for all such $s,s'$. On the other hand, if $s,s'\in S$ are such that $\rho(s) > \rho(s')$, then $\sign(\gamma(s',s,\rho)) = \sign(\max\{0,\rho(s)-\rho(s')\}) = 1$, implying that $\gamma(s',s,\rho) > 0$. Furthermore, in this case with $\rho(s) > \rho(s')$, we see that $\rho(\tilde{s}) - \rho(s) < \rho(\tilde{s}) - \rho(s')$ for all $\tilde{s}\in S$, and therefore $\tau_{\tilde{s}}(\rho(\tilde{s}) - \rho(s')) \ge \tau_{\tilde{s}}(\rho(\tilde{s}) - \rho(s))$ for all $\tilde{s} \in S$ by the fact that every $\tau_{\tilde{s}}$ is nondecreasing. Thus, we immediately see that
		\begin{equation*}
			\sigma(\mu,\rho) \ge 0.
		\end{equation*}
		We furthermore see that if $\sigma(\mu,\rho) = 0$, then
		\begin{equation*}
			\gamma(s',s,\rho)\int_S \left(\tau_{\tilde{s}}(\rho(\tilde{s}) - \rho(s')) - \tau_{\tilde{s}}(\rho(\tilde{s}) - \rho(s))\right) d\lambda(\tilde{s}) = 0 ~ \text{for all $s\in S$ and all $s'\in\supp(\mu)$}
		\end{equation*}
		by the usual arguments based on continuity and nonnegativity of the integrand together with full support of $\lambda$. Thus, let $s\in S$ and $s'\in\supp(\mu)$. Either $\gamma(s',s,\rho) = 0$, or $\int_S \left(\tau_{\tilde{s}}(\rho(\tilde{s}) - \rho(s')) - \tau_{\tilde{s}}(\rho(\tilde{s}) - \rho(s))\right) d\lambda(\tilde{s}) = 0$. In the former case, it must be that $\rho(s) \le \rho(s')$, for otherwise $\phi_s(\rho(s) - \rho(s')) > 0$, which would contradict the fact that $\phi_s(\rho(s) - \rho(s')) = \gamma(s',s,\rho) = 0$. Suppose that the latter case holds. Then either $\rho(s) \le \rho(s')$ or $\rho(s) > \rho(s')$. If $\rho(s) > \rho(s')$, then, as argued above, we find that $\tau_{\tilde{s}}(\rho(\tilde{s}) - \rho(s')) - \tau_{\tilde{s}}(\rho(\tilde{s}) - \rho(s)) \ge 0$ for all $\tilde{s}\in S$, and hence by the usual arguments based on continuity and nonnegativity of the integrand together with the full support of $\lambda$, we conclude that $\tau_{\tilde{s}}(\rho(\tilde{s}) - \rho(s')) = \tau_{\tilde{s}}(\rho(\tilde{s}) - \rho(s))$ for all $\tilde{s}\in S$. In this case, by the fact that every $\tau_{\tilde{s}}$ is strictly increasing on $[0,\infty)$ and $\rho(s') \ne \rho(s)$, it must be the case that, for all $\tilde{s}\in S$, we have that $\rho(\tilde{s}) - \rho(s') \le 0$ and $\rho(\tilde{s}) - \rho(s) \le 0$. But these two inequalities cannot hold simultaneously, as they would imply that $\rho(s) \le \rho(s')$ and $\rho(s') \le \rho(s)$, which contradicts the fact that $\rho(s) > \rho(s')$ in the case under consideration. Hence, we conclude that, when $\sigma(\mu,\rho) = 0$, it must hold that
		\begin{equation*}
			\rho(s) \le \rho(s') ~ \text{for all $s\in S$ and all $s'\in\supp(\mu)$}.
		\end{equation*}
		Thus, by \autoref{lem: equivalent_nash-general}, we find that $\left<\rho,\nu\right> \le \left<\rho,\mu\right>$ for all $\nu\in\P(S)$, and therefore by Nash stationarity of $v$ (\autoref{prop: bnn_pcd_satisfy_ns}), it holds that $v(\mu,\rho) = 0$ whenever $\sigma(\mu,\rho) = 0$. This proves \eqref{eq: delta-dissipative-3}.

		All that remains to be proven is \eqref{eq: delta-dissipative-1} with $w\colon (\mu,\eta) \mapsto \left<\eta,\mu\right>$. Let $\mu\in\P(S)$, $\rho\in C(S)$, and $\eta\in C(S)$. Since $\overline{\storage}(\cdot,\rho)$ is linear, it is immediate that $D(\overline{\storage}(\cdot,\rho))(\mu) = \overline{\storage}(\cdot,\rho)$, which implies that
		\begin{equation*}
			\partial_1 \overline{\storage}(\mu,\rho) v(\mu,\rho) = \overline{\storage}(v(\mu,\rho),\rho) = -\sigma(\mu,\rho).
		\end{equation*}
		Furthermore, computing the second partial Fr\'echet derivative of $\overline{\storage}$ using the chain rule yields that
		\begin{equation*}
			\partial_2 \overline{\storage}(\mu,\rho) \eta = \int_S \int_S \tau_s'(\rho(s)-\rho(s'))(\eta(s)-\eta(s')) d\lambda(s) d\mu(s'),
		\end{equation*}
		where the derivatives of the functions $\tau_s \colon \R\to\R_+$ are computed via the fundamental theorem of calculus:
		\begin{equation*}
			\tau_s'(r) = \frac{d}{dr} \int_{[0,r]} \phi_s(u) du = \phi_s(r).
		\end{equation*}
		By impartiality of the pairwise comparison dynamics under consideration, we find that
		\begin{align*}
			\partial_2 \overline{\storage}(\mu,\rho) \eta &= \int_S \int_S \gamma(s',s,\rho) (\eta(s) - \eta(s')) d\lambda(s) d\mu(s') \\
			&= \int_S \eta(s) \int_S \gamma(s',s,\rho) d\mu(s') d\lambda(s) - \int_S \eta(s') \int_S \gamma(s',s,\rho) d\lambda(s) d\mu(s') \\
			&= \int_S \eta(s') \int_S \gamma(s,s',\rho) d\mu(s) d\lambda(s') - \int_S \eta(s') \int_S \gamma(s',s,\rho) d\lambda(s) d\mu(s') \\
			&= \int_S \eta(s') d(v(\mu,\rho))(s') \\
			&= \left<\eta,v(\mu,\rho)\right> \\
			&= w(v(\mu,\rho),\eta).
		\end{align*}
		Thus, altogether we find that
		\begin{equation*}
			\partial_1 \overline{\storage}(\mu,\rho)v(\mu,\rho) + \partial_2 \overline{\storage}(\mu,\rho) \eta = -\sigma(\mu,\rho) + w(v(\mu,\rho),\eta),
		\end{equation*}
		which shows that \eqref{eq: delta-dissipative-1} holds and hence concludes the proof.
	\end{proofE}
	}{}

	The proof of \autoref{prop: bnn_pairwise_are_delta-passive} relies on generalizing and combining the proof techniques of \citet[Theorem~4.5]{fox2013population}, \citet[Theorem~3]{hofbauer2009brown}, and \citet[Theorem~4]{cheung2014pairwise}. We write the proof in full detail in \iftoggle{longform}{\autoref{sec: proofs}}{our online technical report \citet{anderson2025dissipativity-long}}.

	Finally, we give our reduction of \autoref{thm: delta-dissipative} to the case of monotone games.

	\begin{theorem}
		\label{thm: monotone-passive}
		Consider a game $F \colon \P(S) \to C(S)$\iftoggle{thesis}{,}{ and} let $v\colon\P(S) \times C(S) \to T\P(S)$\iftoggle{thesis}{, and assume that \autoref{ass: exist_edm_solution} holds}{}. Furthermore, assume that \autoref{ass: game} holds and that the extension $\overline{F}$ is continuously Fr\'echet differentiable. If $v$ is Nash stationary, $v$ is $\delta$-passive, and $F$ is monotone, then $\NE(F)$ is weakly Lyapunov stable under the EDM \eqref{eq: edm}. If, additionally, \autoref{ass: game_derivative} and \autoref{ass: dynamics_map} both hold and $v$ is strictly $\delta$-passive, then $\NE(F)$ is globally weakly attracting under the EDM \eqref{eq: edm}.
	\end{theorem}

	\hl{
	\begin{proof}
		Suppose that $v$ is Nash stationary, $v$ is $\delta$-passive, and $F$ is monotone. Let $w\colon \M(S) \times C(S) \to \R$ be defined by $w(\mu,\eta) = \left<\eta,\mu\right>$. Then it holds that $v$ is $\delta$-dissipative with supply rate $w$. Furthermore, by \autoref{lem: monotone_sde}, $F$ satisfies \iftoggle{thesis}{self-defeating externalities}{\eqref{eq: sde}}, and therefore
		\iftoggle{twocol}{
		\begin{equation*}
			w(\nu,D\overline{F}(\mu)\nu) = \left<D\overline{F}(\mu)\nu,\nu\right> \le 0
		\end{equation*}
		for all $\mu\in\P(S)$ and all $\nu\in T\P(S)$.}{
		\begin{equation*}
			w(\nu,D\overline{F}(\mu)\nu) = \left<D\overline{F}(\mu)\nu,\nu\right> \le 0 ~ \text{for all $\mu\in\P(S)$ and all $\nu\in T\P(S)$.}
		\end{equation*}
		} Hence, by \autoref{thm: delta-dissipative}, it holds that $\NE(F)$ is weakly Lyapunov stable under the EDM \eqref{eq: edm}. The fact that $\NE(F)$ is globally weakly attracting under the EDM \eqref{eq: edm} given the additional hypotheses of \autoref{ass: game_derivative} and \autoref{ass: dynamics_map} is immediate from \autoref{thm: delta-dissipative}.
	\end{proof}
	}

	Notice the modularity of \autoref{thm: monotone-passive}: to prove stability of the interconnected EDM \eqref{eq: edm}, we may analyze the Nash stationarity and $\delta$-passivity of the nonautonomous portion of the dynamics defined by the dynamics map $v$ independently from the monotonicity of the system's feedback defined by the game $F$. This allows for the direct proof of stability for the entire class of monotone games $F$ given some dynamics map $v$ that is known to be Nash stationary and $\delta$-passive. For example, \autoref{thm: monotone-passive} together with \autoref{prop: bnn_pairwise_are_delta-passive} recovers the key stability results for BNN dynamics and impartial pairwise comparison dynamics over infinite strategy sets, proven in \citet[Theorem~3]{hofbauer2009brown} and \citet[Theorem~4]{cheung2014pairwise}, respectively. This recovery is formally stated below.

	\iftoggle{longform}{
	\begin{corollaryE}[][end,restate,text link=]
	}{
	\begin{corollaryE}[][end,text link=]
	}
		\label{cor: recovering_monotone}
		Consider a game $F\colon\P(S)\to C(S)$\iftoggle{thesis}{,}{ and} let $v\colon\P(S)\times C(S) \to T\P(S)$\iftoggle{thesis}{, and assume that \autoref{ass: exist_edm_solution} holds}{}. Furthermore, assume that \autoref{ass: game} holds and that the extension $\overline{F}$ is continuously Fr\'echet differentiable. If $F$ is monotone and $v$ is the dynamics map for either the BNN dynamics of \autoref{ex: bnn} or the impartial pairwise comparison dynamics of \autoref{ex: pcd}, then $\NE(F)$ is weakly Lyapunov stable under the EDM \eqref{eq: edm}. If, additionally, \autoref{ass: game_derivative} holds, then $\NE(F)$ is globally weakly attracting under the EDM \eqref{eq: edm}.
	\end{corollaryE}

	\iftoggle{longform}{
	\begin{proofE}
		Notice that weak Lyapunov stability of $\NE(F)$ follows immediately from \autoref{thm: monotone-passive} together with \autoref{prop: bnn_pcd_satisfy_ns} and \autoref{prop: bnn_pairwise_are_delta-passive}. Furthermore, global weak attraction of $\NE(F)$ under \autoref{ass: game_derivative} follows by additionally noting that, for both the BNN dynamics and the impartial pairwise comparison dynamics, $v$ satisfies the appropriate continuity conditions of \autoref{ass: dynamics_map} and $v$ is $\|\cdot\|_\TV$-bounded on weak-$\infty$ compact subsets of $\P(S)\times C(S)$ (the latter condition of which follows from the fact that $v(\mu,\rho)(B) \le 4\|\rho\|_\infty$ for all $\mu\in\P(S)$ and all $\rho\in C(S)$ for the BNN dynamics and that the conditional switch rate $\gamma$ is assumed bounded for the pairwise comparison dynamics).
	\end{proofE}
	}{}

	\iftoggle{longform}{
	\subsection{Extension to Dynamic Payoff Models}
	\label{sec: dynamic_payoff}

	In this section, we consider the case where, instead of static payoff feedback given by $\rho(t) = F(\mu(t))$, as in the EDM \eqref{eq: edm}, the payoff itself has dynamics. \iftoggle{longform}{In doing so, we will consider the derivatives $\dot{\rho}(t)$ of a payoff $\rho \colon [0,\infty) \to C(S)$ (see \autoref{sec: differentiability}). Since $C(S)$ is a Banach space, it holds that $\dot{\rho}(t) \in C(S)$ whenever it exists.}{{\color{black}For the purposes of conciseness and enhanced readability, we give informal versions of the key definitions and results throughout this section, and omit proofs. We refer the interested reader to our online technical report \citet{anderson2025dissipativity-long} for the formal definitions, results, and proofs of this section.}}

    \iftoggle{longform}{
	\begin{definition}
		\label{def: dpedm}
		Let $\mu_0\in\P(S)$, $\rho_0\in C(S)$, $v \colon \P(S) \times C(S) \to \M(S)$, and $u \colon \P(S) \times C(S) \to C(S)$. The differential equation
		\begin{equation}
			\begin{aligned}
				\dot{\mu}(t) &= v(\mu(t),\rho(t)), \\
				\dot{\rho}(t) &= u(\mu(t),\rho(t)), \\
				\mu(0) &= \mu_0, \\
				\rho(0) &= \rho_0,
			\end{aligned}
			\label{eq: dpedm}
		\end{equation}
		is called a \emph{dynamic payoff evolutionary dynamics model (DPEDM)}. The measure $\mu_0$ is called the \emph{initial state}, the function $\rho_0$ is called the \emph{initial payoff}, the mapping $v$ is called the \emph{dynamics map}, and the mapping $u$ is called the \emph{payoff map}. A pair $(\mu,\rho)$ with strongly differentiable $\mu\colon [0,\infty) \to \P(S)$ and differentiable $\rho \colon [0,\infty) \to C(S)$ satisfying \eqref{eq: dpedm} is called a \emph{solution to the DPEDM}.
	\end{definition}
    }{
    \begin{definition}[Informal]
		\label{def: dpedm}
		Let $\mu_0\in\P(S)$, $\rho_0\in C(S)$, $v \colon \P(S) \times C(S) \to \M(S)$, and $u \colon \P(S) \times C(S) \to C(S)$. The differential equation
		\begin{equation}
			\begin{aligned}
				\dot{\mu}(t) &= v(\mu(t),\rho(t)), \\
				\dot{\rho}(t) &= u(\mu(t),\rho(t)), \\
				\mu(0) &= \mu_0, \\
				\rho(0) &= \rho_0,
			\end{aligned}
			\label{eq: dpedm}
		\end{equation}
		is called a \emph{dynamic payoff evolutionary dynamics model (DPEDM)}. The mapping $v$ is called the \emph{dynamics map}, and the mapping $u$ is called the \emph{payoff map}.
	\end{definition}
    }

    {\color{black}\autoref{def: dpedm} may be interpreted as replacing the static feedback controller in the interconnection \autoref{fig: diagram} with a more general controller that has its own underlying dynamics. An example of such a dynamic payoff mechanism occurs when player's exhibit time delays in their strategy revisions (c.f., \autoref{sec: smoothing_dynamics}).} Similar to the case for general EDMs, we will \iftoggle{thesis}{}{always }assume that unique solutions to the DPEDM \eqref{eq: dpedm} exist.

	\iftoggle{thesis}{
	\begin{assumption}
		\label{ass: dpedm_existence}
		Consider a dynamics map $v\colon \P(S) \times C(S) \to \M(S)$ and a payoff map $u \colon \P(S) \times C(S) \to C(S)$. For every initial state $\mu_0 \in \P(S)$ and initial payoff $\rho_0 \in C(S)$, there exists a unique solution $(\mu,\rho)$ to the DPEDM \eqref{eq: dpedm}.
	\end{assumption}
	}{}

        \iftoggle{longform}{
	For games $F\colon \P(S) \to C(S)$ that extend to a Fr\'echet differentiable map $\overline{F} \colon U' \to C(S)$ defined on a strongly open set $U' \subseteq \M(S)$ containing $\P(S)$, the EDM \eqref{eq: edm} is a special case of the DPEDM \eqref{eq: dpedm} with $u \colon (\mu,\rho)\mapsto D\overline{F}(\mu) v(\mu,\rho)$. Since the payoff map in a DPEDM is no longer defined by a static game, the inequality \eqref{eq: supply_inequality} and notions of monotonicity are no longer applicable when characterizing the ``energy supplied'' to the population by the payoff. Instead, we turn to notions of ``antidissipativity.'' The following definition extends such notions from those introduced in \citet{fox2013population} for finite strategy sets to the setting of infinite $S$.
    }{
    Since the payoff map in a DPEDM is no longer defined by a static game, the inequality \eqref{eq: supply_inequality} and notions of monotonicity are no longer applicable when characterizing the ``energy supplied'' to the population by the payoff. {\color{black}Instead, we turn to notions of ``antidissipativity,'' extending the notions introduced in \citet{fox2013population} for finite strategy sets to the setting of infinite $S$. For the formal definitions of antidissipativity in the case of infinite $S$, we refer to our online technical report \citet{anderson2025dissipativity-long}. Intuitively, a payoff map $u$ is called \emph{(strictly) $\delta$-antidissipative with supply rate $\tilde{w}$} if it is (strictly) $\delta$-dissipative with supply rate $-\tilde{w}$ (albeit the notions are technically defined for maps with different codomains). Similarly, $u$ is called \emph{$\delta$-antipassive} if it is $\delta$-antidissipative with supply rate $\tilde{w} \colon (\mu,\nu) \mapsto \int_S \eta d\mu$.}
    }

    \iftoggle{longform}{
	\begin{definition}
		\label{def: delta-antidissipative}
		A map $u \colon \P(S) \times C(S) \to C(S)$ is \emph{$\delta$-antidissipative with supply rate\iftoggle{thesis}{\linebreak}{} $\tilde{w} \colon \M(S) \times C(S) \to \R$} if there exist $\gamma \colon \P(S) \times C(S) \to \Rpl$ and $\Gamma \colon \P(S) \times C(S) \to \Rpl$ that extends to a map $\overline{\Gamma} \colon \tilde{U}\times C(S) \to \R$ with strongly open $\tilde{U}\subseteq \M(S)$ containing $\P(S)$, such that the following conditions hold:
		\begin{enumerate}
			\item $\overline{\Gamma}$ is weak-$\infty$-continuous.
			\item $\overline{\Gamma}$ is Fr\'echet differentiable.
			\item For all strongly differentiable $\mu \colon [0,\infty) \to \P(S)$, all $\rho_0 \in C(S)$, all solutions\iftoggle{thesis}{\linebreak}{} $\rho\colon[0,\infty) \to C(S)$ to the differential equation $\dot{\rho}(t) = u(\mu(t),\rho(t))$ with $\rho(0) = \rho_0$, and all $t\in[0,\infty)$, it holds that
				\iftoggle{twocol}{
				\begin{equation}
					\begin{aligned}
					& \partial_1 \overline{\Gamma}(\mu(t),\rho(t)) \dot{\mu}(t) + \partial_2 \overline{\Gamma}(\mu(t),\rho(t)) u(\mu(t),\rho(t)) \\
					&\qquad \le -\gamma(\mu(t),\rho(t)) - \tilde{w}(\dot{\mu}(t),u(\mu(t),\rho(t))).
					\end{aligned}
					\label{eq: delta-antidissipative-1}
				\end{equation}
				}{
				\begin{equation}
					\partial_1 \overline{\Gamma}(\mu(t),\rho(t)) \dot{\mu}(t) + \partial_2 \overline{\Gamma}(\mu(t),\rho(t)) u(\mu(t),\rho(t)) \le -\gamma(\mu(t),\rho(t)) - \tilde{w}(\dot{\mu}(t),u(\mu(t),\rho(t))).
					\label{eq: delta-antidissipative-1}
				\end{equation}
				}
			\item For all $\mu\in\P(S)$ and all $\rho \in C(S)$, it holds that
				\begin{equation}
					\Gamma(\mu,\rho) = 0 ~ \text{if and only if} ~ u(\mu,\rho) = 0.
					\label{eq: delta-antidissipative-2}
				\end{equation}
		\end{enumerate}
		If, additionally, $(\mu,\rho) \mapsto \partial_1 \overline{\Gamma}(\mu,\rho)$ and $(\mu,\rho) \mapsto \partial_2 \overline{\Gamma}(\mu,\rho)$ are weak-$\infty$-continuous, every partial Fr\'echet derivative $\partial_1 \overline{\Gamma}(\mu,\rho)$ is weakly continuous, and
		\begin{equation}
			\gamma(\mu,\rho) = 0 ~ \text{if and only if} ~ u(\mu,\rho) = 0
			\label{eq: delta-antidissipative-3}
		\end{equation}
		for all $\mu\in\P(S)$ and all $\rho \in C(S)$, then $u$ is \emph{strictly $\delta$-antidissipative with supply rate $\tilde{w}$}.
	\end{definition}

	Notice that $\delta$-antidissipativity is a property solely of a payoff map $u$, and not of any particular dynamics map $v$. One may intuitively think of $\delta$-antidissipativity with supply rate $\tilde{w}$ as $\delta$-dissipativity with supply rate $-\tilde{w}$, albeit the notions are defined for maps with different codomains. We may also define a similar notion that is analogous to $\delta$-passivity.

	\begin{definition}
		\label{def: delta-antipassive}
		A map $u\colon\P(S) \times C(S) \to C(S)$ is \emph{$\delta$-antipassive} if it is $\delta$-antidissipative with supply rate $\tilde{w} \colon (\mu,\eta) \mapsto \left<\eta,\mu\right> = \int_S \eta d\mu$. The map $u$ is \emph{strictly $\delta$-antipassive} if it is strictly $\delta$-antidissipative with such supply rate $\tilde{w}$.
	\end{definition}
        }{}

	\citet[Theorem~4.3]{fox2013population} show that every monotone game over a finite strategy set induces $\delta$-antipassive payoff dynamics,\footnote{Technically, they show $\delta$-antipassivity in the sense of input-output mappings, which slightly differs from the notion of $\delta$-antipassivity of payoff maps used in \iftoggle{thesis}{this thesis}{our paper}.} so $\delta$-antipassivity may be viewed as a generalization of monotonicity to the dynamic payoff setting. \iftoggle{longform}{Before moving on to our generalization of \autoref{thm: delta-dissipative} to the setting of DPEDMs, we remark that \autoref{def: dpedm} does not immediately come equipped with any notion of a game, and hence has no inherent game-theoretic notion of equilibria. The following definition serves to link dynamic payoffs to games, namely, by ensuring that payoffs represent a valid static game at steady state.}{{\color{black}We now report our generalization of the stability result in \autoref{thm: delta-dissipative} to the setting of DPEDMS.}}

    \iftoggle{longform}{
	\begin{definition}
		\label{def: payoff_stationary}
		Consider a game $F\colon \P(S) \to C(S)$. A map $u\colon \P(S) \times C(S) \to C(S)$ is \emph{$F$-payoff stationary} if, for all $\mu\in\P(S)$ and all $\rho\in C(S)$,
		\begin{equation*}
			u(\mu,\rho) = 0
		\end{equation*}
		implies that
		\begin{equation*}
			\rho = F(\mu).
		\end{equation*}
	\end{definition}

	As was the case in the static payoff setting, we need some technical regularity conditions in order to apply Lyapunov theory. \hl{The appropriate conditions are distinct from those for static payoffs (i.e., we no longer need to assume \autoref{ass: game}, \autoref{ass: game_derivative}, or \autoref{ass: dynamics_map} in what follows).} We list the new conditions below, and then state our \hl{stability theorem} for DPEDMs.

	\begin{assumption}
		\label{ass: dpedm_compactness}
		Consider a dynamics map $v\colon \P(S) \times C(S) \to \M(S)$ and a payoff map $u\colon \P(S) \times C(S) \to C(S)$. There exists a compact set $K \subseteq C(S)$ such that, for all initial states $\mu_0\in\P(S)$ and all initial payoffs $\rho_0\in K$, the solution $(\mu,\rho)$ to the DPEDM \eqref{eq: dpedm} satisfies $\rho(t) \in K$ for all $t\in [0,\infty)$.
	\end{assumption}

	\autoref{ass: dpedm_compactness} can be viewed as a ``positive invariance'' condition on the payoff dynamics. Such an assumption on the bounded evolutions of the payoffs is standard in related works (cf., \citealt{kara2023pairwise}) and is necessary to employ Lyapunov theory.

	\begin{assumption}
		\label{ass: weak-infty_continuity}
		The dynamics map $v\colon \P(S) \times C(S) \to \M(S)$ is continuous with respect to the weak-$\infty$ topology on its domain and the weak topology on its codomain. Furthermore, the payoff map $u\colon \P(S) \times C(S) \to C(S)$ is weak-$\infty$-continuous.
	\end{assumption}
        {\color{black}\autoref{ass: weak-infty_continuity} can be thought of as a generalization of \autoref{ass: game} and \autoref{ass: dynamics_map} to the dynamic payoff setting. 
        }
        }{}

        \iftoggle{longform}{
	\begin{theorem}
		\label{thm: dpedm}
		Consider a weakly continuous game $F\colon \P(S) \to C(S)$, let $v \colon \P(S) \times C(S) \to T\P(S)$, and let $u \colon \P(S) \times C(S) \to C(S)$.\iftoggle{thesis}{ Assume that \autoref{ass: dpedm_existence} holds.}{}\iftoggle{thesis}{ Furthermore, assume}{ Assume} that \autoref{ass: dpedm_compactness} holds with some compact $K\subseteq C(S)$ containing $F(\NE(F))$, and that \autoref{ass: weak-infty_continuity} holds. If $v$ is Nash stationary and $\delta$-dissipative with supply rate $w\colon \M(S) \times C(S) \to \R$ and $u$ is $F$-payoff stationary and $\delta$-antidissipative with supply rate $\tilde{w} \ge w$, then
		\begin{equation*}
		P\coloneqq \{(\mu,\rho) \in \P(S) \times C(S) : v(\mu,\rho) = 0, ~ u(\mu,\rho) = 0\}
		\end{equation*}
		is a subset of $\NE(F) \times F(\NE(F))$ and is weak-$\infty$-Lyapunov stable under the DPEDM \eqref{eq: dpedm}. If, additionally, the $\delta$-dissipativity of $v$ and the $\delta$-antidissipativity of $u$ are both strict and $v$ is $\|\cdot\|_\TV$-bounded on $\P(S) \times K$, then $P$ is weak-$\infty$-attracting under the DPEDM \eqref{eq: dpedm} from every $(\mu_0,\rho_0) \in \P(S) \times K$.
	\end{theorem}
        }{
        {\color{black}
        \begin{theorem}[Informal]
		\label{thm: dpedm}
		Consider a weakly continuous game $F\colon \P(S) \to C(S)$, and consider the DPEDM \eqref{eq: dpedm} with some adequately continuous maps $v \colon \P(S) \times C(S) \to T\P(S)$ and $u \colon \P(S) \times C(S) \to C(S)$. If the following three conditions hold, then every stationary population state of the DPEDM is a Nash equilibrium of $F$, and the set of such stationary population states is Lyapunov stable:
        \begin{enumerate}
            \item The payoffs $\rho(t)$ evolve in a compact positively invariant set containing $F(\NE(F))$.
            \item The dynamics map $v$ is Nash stationary and $\delta$-dissipative with supply rate $w$.
            \item The rest points of the payoff map $u$ occur at payoffs $\rho=F(\mu)$ of the game $F$, and $u$ is $\delta$-antidissipative with supply rate $\tilde{w} \ge w$.
        \end{enumerate}
        If, in addition to the above three conditions, the dynamics map $v$ is bounded and the dissipativity properties of $v$ and $u$ are both strict, then the stationary population states (which are Nash equilibria of $F$) are attracting from every initial condition of the DPEDM.
	\end{theorem}
        }
        }

	\iftoggle{longform}{
	\begin{proof}
        We break the proof up into two parts: first, the proof of weak-$\infty$ Lyapunov stability, and second, the proof of global weak-$\infty$ attraction.
        
        \paragraph{Proof of weak-$\infty$ Lyapunov stability\iftoggle{longform}{.}{}}
        Since $v$ is $\delta$-dissipative with supply rate $w\colon \M(S) \times C(S) \to \R$, there exist $\sigma \colon \P(S) \times C(S) \to \Rpl$ and $\storage\colon\P(S) \times C(S) \to \Rpl$ with $\storage$ having an appropriate extension $\overline{\storage} \colon U\times C(S) \to \R$ as in \autoref{def: delta-dissipative}. Furthermore, since $u$ is $\delta$-antidissipative with supply rate $\tilde{w}\colon \M(S) \times C(S) \to \R$, there exist $\gamma \colon \P(S) \times C(S) \to \Rpl$ and $\Gamma\colon\P(S) \times C(S) \to \Rpl$ with $\Gamma$ having an appropriate extension $\overline{\Gamma} \colon \tilde{U}\times C(S) \to \R$ as in \autoref{def: delta-antidissipative}. Define $V \colon \P(S)\times K \to \Rpl$ by $V(\mu,\rho) = \storage(\mu,\rho) + \Gamma(\mu,\rho)$. Consider $P = \{(\mu,\rho) \in \P(S) \times C(S) : v(\mu,\rho) = 0, ~ u(\mu,\rho) = 0\}$. If $(\mu,\rho) \in P$, then $u(\mu,\rho) = 0$, implying that $\rho = F(\mu)$ by $F$-payoff stationarity, and hence $v(\mu,F(\mu)) = 0$, so $\mu \in \NE(F)$ by Nash stationarity. Thus, $P \subseteq \NE(F) \times F(\NE(F)) \subseteq \P(S) \times K$. By \autoref{prop: ne_closed} and \autoref{prop: weakly_continuous_payoff}, $\NE(F)$ is weakly compact, and hence $F(\NE(F))$ is compact as $F$ is weakly continuous. Since $v$ is continuous with respect to the weak-$\infty$ topology on its domain and the weak topology on its codomain, and $u$ is weak-$\infty$-continuous, it holds that $P = v^{-1}(\{0\})\cap u^{-1}(\{0\})$ is weak-$\infty$-closed, and hence must be weak-$\infty$-compact as well as $\NE(F) \times F(\NE(F))$ is.\iftoggle{longform}{ Thus, by \autoref{lem: lyapunov}, it suffices to show that $V$ is a global Lyapunov function for $P$ under $(\mu,\rho) \mapsto (v(\mu,\rho),u(\mu,\rho))$ (according to \autoref{def: lyapunov_function}).}{ Thus, it suffices to show that $V$ is a global Lyapunov function for $P$ under $(\mu,\rho) \mapsto (v(\mu,\rho),u(\mu,\rho))$ (see Appendix B in our online technical report \citet{anderson2025dissipativity-long}).} Let $\overline{V}\colon U\cap \tilde{U} \times C(S) \to \R$ be defined by $\overline{V}(\mu,\rho) = \overline{\storage}(\mu,\rho) + \overline{\Gamma}(\mu,\rho)$. Note that $U\cap \tilde{U}$ is strongly open and contains $\P(S)$, and that $\overline{V}$ is weak-$\infty$-continuous and Fr\'echet differentiable since $\overline{\storage}$ and $\overline{\Gamma}$ are. Also note that $\overline{V}(\mu,\rho) = V(\mu,\rho)$ for all $(\mu,\rho) \in \P(S) \times K$. Also, if $(\mu,\rho) \in P$, then $v(\mu,\rho) = 0$ and $u(\mu,\rho) = 0$, so $\overline{V}(\mu,\rho) = 0$ by \eqref{eq: delta-dissipative-2} and \eqref{eq: delta-antidissipative-2}. Furthermore, if $(\mu,\rho) \in (\P(S) \times K) \setminus P$, then again by \eqref{eq: delta-dissipative-2} and \eqref{eq: delta-antidissipative-2} we have that $\overline{V}(\mu,\rho) > 0$. \iftoggle{longform}{Therefore, the first two conditions from \autoref{def: lyapunov_function} on $V$ to be a global Lyapunov function for $P$ under $(v,u)$ are satisfied.}{}

		Next, it holds for all $(\mu,\rho) \in \P(S) \times K$ that
		\iftoggle{twocol}{
		\begin{equation}
		\begin{aligned}
			D\overline{V}(\mu,\rho)(v(\mu,\rho),u(\mu,\rho)) &= \partial_1 \overline{\storage}(\mu,\rho) v(\mu, \rho) \\
			&\qquad + \partial_2 \overline{\storage}(\mu,\rho) u(\mu, \rho) \\
			&\qquad + \partial_1 \overline{\Gamma}(\mu,\rho) v(\mu, \rho) \\
			&\qquad + \partial_2 \overline{\Gamma}(\mu,\rho) u(\mu, \rho) \\
			&\le -\sigma(\mu,\rho) \\
			&\qquad + w(v(\mu,\rho),u(\mu,\rho)) \\
			&\qquad - \gamma(\mu,\rho) \\
			&\qquad - \tilde{w}(v(\mu,\rho),u(\mu,\rho)) \\
			&\le -\sigma(\mu,F(\mu)) - \gamma(\mu,\rho) \\
			&\le 0,
		\end{aligned}
		\label{eq: decreasing_derivative-dpedm}
		\end{equation}
		}{
		\begin{equation}
		\begin{aligned}
			D\overline{V}(\mu,\rho)(v(\mu,\rho),u(\mu,\rho)) &= \partial_1 \overline{\storage}(\mu,\rho) v(\mu, \rho) + \partial_2 \overline{\storage}(\mu,\rho) u(\mu, \rho) \\
			&\qquad + \partial_1 \overline{\Gamma}(\mu,\rho) v(\mu, \rho) + \partial_2 \overline{\Gamma}(\mu,\rho) u(\mu, \rho) \\
			&\le -\sigma(\mu,\rho) + w(v(\mu,\rho),u(\mu,\rho)) - \gamma(\mu,\rho) - \tilde{w}(v(\mu,\rho),u(\mu,\rho)) \\
			&\le -\sigma(\mu,F(\mu)) - \gamma(\mu,\rho) \\
			&\le 0,
		\end{aligned}
		\label{eq: decreasing_derivative-dpedm}
		\end{equation}
		}
		where the first inequality follows from \eqref{eq: delta-dissipative-1} and \eqref{eq: delta-antidissipative-1}, and the second inequality follows from the fact that $\tilde{w} \ge w$. Hence, $V$ is indeed a global Lyapunov function for $P$ under $(v,u)$, so $P$ is weak-$\infty$-Lyapunov stable under the DPEDM \eqref{eq: dpedm}.

        \paragraph{Proof of global weak-$\infty$ attraction\iftoggle{longform}{.}{}}
        Now suppose that the $\delta$-dissipativity of $v$ and the $\delta$-antidissipativity of $u$ are both strict, and that $v$ is $\|\cdot\|_\TV$-bounded on $\P(S)\times K$.\iftoggle{longform}{ By \autoref{lem: strict_lyapunov}, it suffices to show that $V$ is a strict global Lyapunov function for $P$ under $(\mu,\rho) \mapsto (v(\mu,\rho),u(\mu,\rho))$ (according to \autoref{def: lyapunov_function}).}{ It suffices to show that $V$ is a strict global Lyapunov function for $P$ under $(\mu,\rho) \mapsto (v(\mu,\rho),u(\mu,\rho))$ (see Appendix B in our online technical report \citet{anderson2025dissipativity-long}).} This amounts to proving that $(\mu,\rho) \mapsto D\overline{V}(\mu,\rho)(v(\mu,\rho),u(\mu,\rho))$ is weak-$\infty$-continuous and that $D\overline{V}(\mu,\rho)(v(\mu,\rho),u(\mu,\rho)) < 0$ for all $(\mu,\rho)\in(\P(S) \times K)\setminus P$. Indeed, the continuity condition holds by \autoref{lem: weakly_continuous_derivative-dpedm}, which we state and prove in \autoref{sec: proofs}.

		Next, if $(\mu,\rho)\in(\P(S)\times K)\setminus P$, then $v(\mu,\rho) \ne 0$ or $u(\mu,\rho) \ne 0$, so $\sigma(\mu,\rho) > 0$ or $\gamma(\mu,\rho) > 0$ by \eqref{eq: delta-dissipative-3} and \eqref{eq: delta-antidissipative-3}, implying that $D\overline{V}(\mu,\rho)(v(\mu,\rho),u(\mu,\rho)) < 0$ for all such $(\mu,\rho)$ by \eqref{eq: decreasing_derivative-dpedm}. Hence, $V$ is indeed a strict global Lyapunov function for $P$ under $(v,u)$, so $P$ is globally weak-$\infty$-attracting under the DPEDM \eqref{eq: dpedm} from $K$.
	\end{proof}

	\begin{lemmaE}[][all end,text link=]
		\label{lem: weakly_continuous_derivative-dpedm}
		\iftoggle{longform}{The map }{}$(\mu,\rho)\mapsto D\overline{V}(\mu,\rho)(v(\mu,\rho),u(\mu,\rho))$ is weakly continuous.
	\end{lemmaE}

	\begin{proofE}
		The result follows from a nearly identical analysis as in the proof of \autoref{lem: weakly_continuous_derivative} with minor changes. In particular, it follows from the $\|\cdot\|_\TV$-boundedness of $v$ on $\P(S)\times K$, the weak-$\infty$-to-weak continuity of $v$, the weak-$\infty$ continuity of $u$, the weak continuity of every $\partial_1 \overline{\storage}(\mu,\rho)$ and every $\partial_1 \overline{\Gamma}(\mu,\rho)$, and the weak-$\infty$ continuity of the maps $(\mu,\rho)\mapsto \partial_1 \overline{\storage}(\mu,\rho)$, $(\mu,\rho)\mapsto \partial_2 \overline{\storage}(\mu,\rho)$, $(\mu,\rho)\mapsto \partial_1 \overline{\Gamma}(\mu,\rho)$, and $(\mu,\rho)\mapsto \partial_2 \overline{\Gamma}(\mu,\rho)$.
	\end{proofE}
        }{}

    \iftoggle{longform}{
	The set $P$ in \autoref{thm: dpedm} corresponds to the set of rest points of the DPEDM \eqref{eq: dpedm}. The result shows that, under the appropriate regularity conditions, the DPEDM has\iftoggle{thesis}{ dynamically}{} stable rest points whenever the dynamics map is $\delta$-dissipative and the payoff map is $\delta$-antidissipative, and the incoming energy supply rate to the dynamics is less than that of the payoffs. Since, under the hypotheses of the theorem, it holds that $P \subseteq \{(\mu,\rho) \in \P(S) \times C(S) : \rho = F(\mu), ~ \mu\in\NE(F)\} \subseteq \NE(F) \times F(\NE(F))$, the result shows convergence of the $\mu$-component of the trajectory $(\mu,\rho)$ to $\NE(F)$, and convergence of the $\rho$-component to the corresponding static payoff given by the game $F$.

    As was the case in the static payoff setting, the appropriate continuity properties, as well as the new positive invariance assumption on the payoff evolution, are technical conditions needed in \autoref{thm: dpedm} in order to rigorously apply Lyapunov theory. Similar to the intuition of the compactness condition on $\NE(F)$ for stability of EDMs, the compact positive invariance condition of \autoref{ass: dpedm_compactness} is key in ensuring that the payoffs in a DPEDM do not drift towards some equilibrium ``out at infinity.''
    }{
    {\color{black}Intuitively, \autoref{thm: dpedm} shows that, under the appropriate regularity conditions, the DPEDM has stable Nash equilibria whenever the dynamics and payoffs are dissipative, and the incoming energy supply rate to the dynamics is less than that of the payoffs. As was the case in the static payoff setting, the appropriate continuity properties, as well as the new positive invariance assumption on the payoff evolution, are technical conditions needed in order to rigorously apply Lyapunov theory. Similar to the intuition of the compactness condition on $\NE(F)$ for stability of EDMs, the above compact positive invariance condition is key in ensuring that the payoffs in a DPEDM do not drift towards some equilibrium ``out at infinity.''}
    }

    {\color{black}    
    Similar to the static payoff setting, it is easy to see that our dissipativity-based result \autoref{thm: dpedm} may be specialized to the case \iftoggle{longform}{of $\delta$-passive dynamics maps coupled with $\delta$-antipassive payoff maps}{in which the dynamics and payoffs are passive}, resulting in analogues to \autoref{thm: monotone-passive} and \autoref{cor: recovering_monotone} for the dynamic payoff setting. In particular, the latter specialization yields the following result, which is stronger than \citet[Theorem~3]{hofbauer2009brown} and \citet[Theorem~4]{cheung2014pairwise}, as it allows for $\delta$-antipassive dynamic payoffs, which are more general than the static payoffs $\rho(t) = F(\mu(t))$ with monotone $F$ considered in these prior works.
    }

    \iftoggle{longform}{
	\iftoggle{longform}{
	\begin{corollaryE}[][end,restate,text link=]
	}{
	\begin{corollaryE}[][end,text link=]
	}
		\label{cor: dpedm-bnn_and_pcd}
		Consider a weakly continuous game $F\colon \P(S) \to C(S)$, let $v\colon \P(S) \times C(S) \to T\P(S)$, and let $u\colon \P(S) \times C(S) \to C(S)$.\iftoggle{thesis}{ Assume that \autoref{ass: dpedm_existence} holds.}{}\iftoggle{thesis}{ Furthermore, assume}{ Assume} that \autoref{ass: dpedm_compactness} holds with some compact $K\subseteq C(S)$ containing $F(\NE(F))$, and that \autoref{ass: weak-infty_continuity} holds. If $v$ is the dynamics map for either the BNN dynamics of \autoref{ex: bnn} or the impartial pairwise comparison dynamics of \autoref{ex: pcd} and $u$ is $F$-payoff stationary and strictly $\delta$-antipassive, then
		\begin{equation*}
		P \coloneqq \{(\mu,\rho) \in \P(S) \times C(S) : v(\mu,\rho) = 0, ~ u(\mu,\rho) = 0\}
		\end{equation*}
		is a subset of $\NE(F) \times F(\NE(F))$ and is weak-$\infty$-Lyapunov stable under the DPEDM \eqref{eq: dpedm} and weak-$\infty$-attracting under the DPEDM \eqref{eq: dpedm} from every $(\mu_0,\rho_0)\in \P(S) \times K$.
	\end{corollaryE}

	\iftoggle{longform}{
	\begin{proofE}
		The proof follows analogously to that of \autoref{cor: recovering_monotone}.
	\end{proofE}
	}{}
        }
        {
        {\color{black}
        \begin{corollary}[Informal]
            \label{cor: dpedm-bnn_and_pcd}
            Consider a weakly continuous game $F\colon \P(S) \to C(S)$, and consider the DPEDM \eqref{eq: dpedm} with some adequately continuous payoff map $u \colon \P(S) \times C(S) \to C(S)$ and with $v \colon \P(S) \times C(S) \to T\P(S)$ being the dynamics map for either the BNN dynamics of \autoref{ex: bnn} or the impartial pairwise comparison dynamics of \autoref{ex: pcd}. If the following two conditions hold, then every stationary population state of the DPEDM is a Nash equilibrium of $F$, and the set of such stationary population states is both Lyapunov stable and attracting:
            \begin{enumerate}
                \item The payoffs $\rho(t)$ evolve in a compact positively invariant set containing $F(\NE(F))$.
                \item The rest points of the payoff map $u$ occur at payoffs $\rho=F(\mu)$ of the game $F$, and $u$ is strictly $\delta$-antipassive.
            \end{enumerate}
        \end{corollary}
        }
        }

        {\color{black} \autoref{cor: dpedm-bnn_and_pcd} allows one to easily assess the stability of the popular BNN and impartial pairwise comparison dynamics, in the case where the payoffs take a general dynamic form. We study a practical example of such dynamic payoffs coupled with the BNN dynamics in \autoref{sec: smoothing_dynamics}, where smoothing effects, such as time delays in player strategy revisions, result in suppressing short-term variations in the payoff evolution.}
	}{} 

	\section{Case Studies}
\label{sec: examples}

\subsection{War of Attrition---Failure of Finite Approximations}
\label{sec: failure_finite}

\hl{In this section, we tie together various results from the literature to generate an example in which approximations of an infinite-dimensional evolutionary game are proven to be\iftoggle{thesis}{ dynamically}{} stable via finite-dimensional dissipativity analysis, yet the true underlying dynamics do not weakly converge to the set of Nash equilibria. This finding shows that one cannot in general use finite-dimensional dissipativity theory to assess the stability of evolutionary games over infinite strategy sets, further motivating our studies directly concerned with the infinite-dimensional regime.}

Consider endowing the {\color{black}war of attrition game $F$ of \autoref{ex: war_of_attrition}} with the BNN dynamics of \autoref{ex: bnn}, where the reference measure $\lambda$ is Lebesgue. We will now show that the prior dissipativity results of \citet{arcak2021dissipativity} guarantee that finite-strategy approximations of this evolutionary game asymptotically converge to their unique Nash equilibrium. Despite this, we will find that the infinite-dimensional dynamics do not \hl{weakly} converge to the unique Nash equilibrium $\mu^\star$, \hl{justifying} the need for direct consideration of dissipativity theory over infinite strategy sets, as we have done in this \iftoggle{thesis}{thesis}{paper}.

Let $n\in\N$ and consider a finite approximation of the strategy set given by $S_n = \{s_1,\dots,s_n\} \subseteq S$, with $s_1< s_2 <\cdots < s_n$. Restricting the game $F$ to the set of measures
\begin{equation*}
	\D(S_n) \coloneqq \left\{\sum_{i=1}^n x_i \delta_{s_i} \in \P(S) : x\in\Delta^{n-1}\right\}
\end{equation*}
with $\Delta^{n-1} \coloneqq \{x\in\Rn : \text{$x_i \ge 0$ for all $i$}, ~ \sum_{i=1}^n x_i = 1\}$ yields the finite-dimensional approximation $\hat{F}_n \colon \Delta^{n-1} \to \Rn$ given by
\begin{equation*}
	(\hat{F}_n(x))_i \coloneqq F\left(\sum_{j=1}^n x_j \delta_{s_j}\right)(s_i) = \sum_{j=1}^n x_j f(s_i,s_j).
\end{equation*}
Thus, the finite-dimensional game may be written as
\begin{equation*}
	\hat{F}_n(x) = A_n x,
\end{equation*}
where
\iftoggle{twocol}{
\begin{align*}
	A_n \coloneqq& \begin{bmatrix} f(s_1,s_1) & f(s_1,s_2) & \cdots & f(s_1,s_n) \\
	f(s_2,s_1) & f(s_2,s_2) & \cdots & f(s_2,s_n) \\
	\vdots & \vdots & \ddots & \vdots \\
	f(s_n,s_1) & f(s_n,s_2) & \cdots & f(s_n,s_n)
\end{bmatrix} \\
	=&
\begin{bmatrix}
	\frac{V}{2} - s_1 & -s_1 & \cdots & -s_1 \\
	V - s_1 & \frac{V}{2} - s_2 & \cdots & -s_2 \\
	\vdots & \vdots & \ddots & \vdots \\
	V - s_1 & V - s_2 & \cdots & \frac{V}{2} - s_n
\end{bmatrix}
\in \R^{n\times n}.
\end{align*}
}{
\begin{equation*}
	A_n \coloneqq \begin{bmatrix} f(s_1,s_1) & f(s_1,s_2) & \cdots & f(s_1,s_n) \\
	f(s_2,s_1) & f(s_2,s_2) & \cdots & f(s_2,s_n) \\
	\vdots & \vdots & \ddots & \vdots \\
	f(s_n,s_1) & f(s_n,s_2) & \cdots & f(s_n,s_n)
\end{bmatrix} =
\begin{bmatrix}
	\frac{V}{2} - s_1 & -s_1 & \cdots & -s_1 \\
	V - s_1 & \frac{V}{2} - s_2 & \cdots & -s_2 \\
	\vdots & \vdots & \ddots & \vdots \\
	V - s_1 & V - s_2 & \cdots & \frac{V}{2} - s_n
\end{bmatrix}
\in \R^{n\times n}.
\end{equation*}
}

This finite-dimensional game $\hat{F}_n$ is monotone \citep{hofbauer2009stable}. The corresponding finite-dimensional BNN dynamics are given by
\iftoggle{twocol}{
\begin{align*}
	\dot{x}_i(t) &= \max\{0, (\hat{F}_n(x(t)))_i - x(t)^\top \hat{F}_n(x(t))\} \\
	&\qquad - x_i(t) \sum_{i=1}^n \max\{0, (\hat{F}_n(x(t)))_i - x(t)^\top \hat{F}_n(x(t))\}
\end{align*}
}{
\begin{equation*}
	\dot{x}_i(t) = \max\{0, (\hat{F}_n(x(t)))_i - x(t)^\top \hat{F}_n(x(t))\} - x_i(t) \sum_{i=1}^n \max\{0, (\hat{F}_n(x(t)))_i - x(t)^\top \hat{F}_n(x(t))\}
\end{equation*}
}
for all $i\in\{1,\dots,n\}$ (cf., \citealt[Example~4.3.4]{sandholm2010population}). These finite-dimensional BNN dynamics are Nash stationary and $\delta$-passive \citep{arcak2021dissipativity}, and therefore since $\hat{F}_n$ is monotone and admits a continuously differentiable extension defined on $\Rn$ (given by the linear map defined by $A_n$), \citet[Theorem~1]{arcak2021dissipativity} asserts that $\NE(\hat{F}_n)$ is globally asymptotically stable under these dynamics.

Based on the above analysis, one may hope that $\NE(F) = \{\mu^\star\}$ is also globally weakly attracting under the infinite-dimensional EDM \eqref{eq: edm}. However, despite $F$ being monotone and $v$ being Nash stationary and strictly $\delta$-passive, this is not the case, as \citet[Example~6]{hofbauer2009brown} {\color{black}argues} that this infinite-dimensional dynamic does not weakly converge to $\mu^\star$. The intuition for this lack of convergence given by \citet{hofbauer2009brown} is that the BNN dynamics under Lebesgue measure cannot generate mass at $s = T$, but that a point mass at $s=T$ is present in the equilibrium distribution $\mu^\star$.

{\color{black}We now verify this nonconvergence by simulating the BNN dynamics with $T=2$ and $V=1$ for a variety of initial population states $\mu_0$, and plotting the integral values $\int_S g d\mu(t)$ for a bounded continuous function $g$ against time. The function $g$ is specified to be given by
\begin{equation*}
    g(s) = \max\left\{ 0 , \frac{1}{T - s^\star}(s - s^\star) \right\},
\end{equation*}
where $s^\star = T - V/2 = 3/2$. Notice that $g(s) \in [0,1]$ for all $s\in S = [0,T]$. In this case, the integral values reduce to
\begin{equation*}
    \int_S g d\mu(t) = \int_{[0,T)} g d\mu(t),
\end{equation*}
whenever $\mu_0(\{T\}) = 0$ (since then, the BNN dynamics give that $\mu(t)(\{T\}) = 0$ for all time $t$). Notice that, using the above finite approximation, this gives that
\iftoggle{twocol}{
\begin{equation*}
    \int_S g d\mu(t) \approx \int_{[0,T)} g d \left( \sum_{i=1}^n x_i(t) \delta_{s_i} \right) = \hspace*{-0.5em} \sum_{\substack{i=1 \\ s^\star \le s_i < T}}^n \hspace*{-0.5em} x_i(t) \frac{s_i - s^\star}{T - s^\star}.
\end{equation*}
}{
\begin{equation*}
    \int_S g d\mu(t) \approx \int_{[0,T)} g d \left( \sum_{i=1}^n x_i(t) \delta_{s_i} \right) = \sum_{\substack{i=1 \\ s^\star \le s_i < T}}^n x_i(t) \frac{s_i - s^\star}{T - s^\star}.
\end{equation*}
}
These integral values are plotted in \autoref{fig: compare_int} for simulated BNN dynamics starting from each of the following initial states (which all satisfy $\mu_0(\{T\})=0$):
\begin{enumerate}
    \item Gaussian $\mu_0$ with mean $1$ and variance $0.1$.
    \item Half-uniform $\mu_0$ with distribution function
    \begin{equation*}
        \mu_0([0,s]) = \begin{cases}
            0 & \text{if $s\in[0,T/2)$}, \\
            \frac{1}{T/2}(s-T/2) & \text{if $s\in[T/2,T]$}.
        \end{cases}
    \end{equation*}
    \item Dirac $\mu_0 = \delta_{T/2}$.
    \item Multiple Diracs $\mu_0 = \frac{1}{3}(\delta_{T/4} + \delta_{T/2} + \delta_{3T/4})$.
\end{enumerate}
It is observed that the integral values all converge to zero, yet the Nash equilibrium $\mu^\star$ gives
\begin{equation*}
    \int_S g d\mu^\star = e^{-s^\star / V} = e^{-3/2} \approx 0.2231.
\end{equation*}
Thus, the infinite-dimensional BNN dynamics do not converge weakly to the Nash equilibrium $\mu^\star$, despite the corresponding finite-strategy game asymptotically converging to its Nash equilibrium.
}

\begin{figure}[ht]
    \centering
%
\definecolor{mycolor1}{rgb}{0.00000,0.44700,0.74100}%
\definecolor{mycolor2}{rgb}{0.85000,0.32500,0.09800}%
\definecolor{mycolor3}{rgb}{0.92900,0.69400,0.12500}%
\definecolor{mycolor4}{rgb}{0.49400,0.18400,0.55600}%
\definecolor{mycolor5}{rgb}{0.46600,0.67400,0.18800}%
\begin{tikzpicture}

\begin{axis}[%
width=0.475\linewidth,
height=0.375\linewidth,
at={(0\linewidth,0\linewidth)},
scale only axis,
xmin=0,
xmax=10,
xtick={0,5,10},
xlabel style={font=\color{white!15!black}},
xlabel={Time},
ymin=0,
ymax=0.25,
ytick={0,0.1,0.2,0.3},
yticklabels={{0},{0.1},{0.2},{}},
ylabel style={font=\color{white!15!black}},
ylabel={Integral value},
axis background/.style={fill=white},
title style={font=\normalsize},xlabel style={font=\normalsize},ylabel style={font=\normalsize},legend style={font=\small},ticklabel style={font=\normalsize}
]
\addplot [color=mycolor1, dotted, line width=2.5pt, forget plot]
  table[row sep=crcr]{%
0	0.223130160148429\\
10	0.223130160148429\\
};
\addplot [color=mycolor2, line width=2.5pt, forget plot]
  table[row sep=crcr]{%
0	0.0144360731384197\\
0.0500000000000007	0.141781481239368\\
0.0999999999999996	0.150526478911358\\
0.15	0.155830543630699\\
0.199999999999999	0.160641465659683\\
0.25	0.164323275618296\\
0.300000000000001	0.166804778047824\\
0.35	0.168034512818942\\
0.4	0.168097199677872\\
0.449999999999999	0.167257773523255\\
0.5	0.165793229751751\\
0.550000000000001	0.164024364510448\\
0.65	0.160196520953793\\
0.699999999999999	0.158156706595923\\
0.800000000000001	0.153654244849349\\
0.85	0.151794824584952\\
0.9	0.150339849854758\\
0.949999999999999	0.149030685413553\\
1	0.147482740502968\\
1.05	0.145706133794716\\
1.15	0.141848567776938\\
1.3	0.135996577899762\\
1.45	0.130311515249073\\
1.6	0.124759283770153\\
1.7	0.121288195910063\\
1.8	0.118045793196174\\
1.9	0.115004781350757\\
2	0.112124598646792\\
2.1	0.10937933258618\\
2.2	0.106760451543087\\
2.3	0.104248271317882\\
2.4	0.101853535895035\\
2.5	0.0995755046456033\\
2.6	0.0974164936467705\\
2.7	0.0953537194164067\\
2.85	0.0924270253779049\\
3	0.089665064025457\\
3.15	0.0870507597837555\\
3.3	0.0845847773371737\\
3.45	0.0822580578870546\\
3.6	0.0800698829176429\\
3.75	0.0780005210813481\\
3.95	0.0754008635760854\\
4.15	0.0729607975535966\\
4.35	0.0706620546091088\\
4.55	0.0685102270447722\\
4.75	0.0665052860466719\\
5.05	0.0637139718812811\\
5.3	0.0615489177691178\\
5.55	0.0595154006818976\\
5.8	0.0576185902827362\\
6.1	0.0555178009376309\\
6.4	0.0535747390029435\\
6.65	0.0520334711457355\\
7	0.0500259031977048\\
7.3	0.0484231174714971\\
7.6	0.0469343489800575\\
8.05	0.0448636280837302\\
8.45	0.0431593379732309\\
8.85	0.0415742600121156\\
9.5	0.0392772621343962\\
10	0.0376530275183615\\
};
\addplot [color=mycolor3, line width=2.5pt, forget plot]
  table[row sep=crcr]{%
0	0.24306750812775\\
0.0500000000000007	0.158419481077855\\
0.0999999999999996	0.15466583559102\\
0.15	0.160083707471969\\
0.199999999999999	0.164070553775124\\
0.25	0.165322894497351\\
0.300000000000001	0.16485774217132\\
0.35	0.163618327246901\\
0.5	0.159320042648943\\
0.550000000000001	0.15789214701236\\
0.6	0.156363733784231\\
0.65	0.154663703118137\\
0.75	0.150979353410627\\
0.85	0.147104921163274\\
1	0.1410254884883\\
1.1	0.13723367812867\\
1.2	0.133644745581616\\
1.35	0.128507692001238\\
1.6	0.120054742750938\\
1.7	0.116674031295277\\
1.8	0.113558559737136\\
1.9	0.110659171471104\\
2.05	0.106542989117971\\
2.2	0.102648018387949\\
2.3	0.100357251009589\\
2.4	0.0982663479539134\\
2.5	0.0961985854357046\\
2.65	0.0929156735587142\\
2.75	0.0909076029917486\\
2.85	0.0891131864219012\\
2.95	0.0875143366262048\\
3.2	0.0838398260519053\\
3.5	0.0795220981951914\\
3.7	0.0767817579796244\\
3.9	0.0742306653714575\\
4.1	0.071848759563883\\
4.3	0.0696009688340471\\
4.5	0.0674926223963777\\
4.7	0.0655078373265852\\
4.9	0.0636499407757203\\
5.15	0.0614996301063311\\
5.45	0.0590987567381465\\
5.75	0.0568453545653185\\
6.05	0.0547668261961682\\
6.3	0.0531666822333641\\
6.65	0.0510862115009392\\
7.25	0.0478304222648287\\
7.6	0.0461385813014612\\
8.05	0.0441266631755273\\
8.45	0.0424733591988584\\
8.85	0.0409412378432989\\
9.35	0.0391923658860165\\
9.9	0.0374156832911243\\
10	0.037111911832568\\
};
\addplot [color=mycolor4, line width=2.5pt, forget plot]
  table[row sep=crcr]{%
0	0\\
0.0500000000000007	0.142674405418576\\
0.0999999999999996	0.151610682441804\\
0.15	0.157123317568038\\
0.199999999999999	0.161360334571311\\
0.25	0.163362698889154\\
0.300000000000001	0.163530394804326\\
0.35	0.162571376038215\\
0.4	0.161004289043216\\
0.800000000000001	0.146405311000061\\
0.85	0.14444294428394\\
0.9	0.142402058523624\\
1.1	0.133911334892204\\
1.2	0.129934901767754\\
1.3	0.126187629313451\\
1.4	0.122626035176777\\
1.5	0.119205879837354\\
1.7	0.112528690002025\\
1.8	0.109317531209117\\
1.9	0.106304348068985\\
2	0.103459542857406\\
2.1	0.100734701640171\\
2.2	0.0981267603264779\\
2.3	0.0956369286313308\\
2.4	0.0932740770444198\\
2.5	0.0910389982575879\\
2.6	0.0889319041378336\\
2.7	0.0869475889315883\\
2.8	0.0850769635400042\\
2.95	0.0824545130882672\\
3.1	0.0800116373279831\\
3.25	0.0777268876386792\\
3.4	0.0755717587585991\\
3.55	0.0735412067212078\\
3.75	0.0710047509133922\\
3.95	0.0686183307997794\\
4.2	0.0657756330586814\\
4.4	0.0636185404138487\\
4.6	0.0616277446907691\\
4.8	0.0598266585471201\\
5.05	0.0577639180961835\\
5.35	0.055447196176635\\
5.65	0.0532859060480995\\
5.9	0.0516179664006682\\
6.2	0.0497693802493\\
6.55	0.0477731150945999\\
6.95	0.0456387426470464\\
7.25	0.0441607964291553\\
7.6	0.0425998982008267\\
8.15	0.0403336946564181\\
8.5	0.0390061363199639\\
8.9	0.0376539222102661\\
9.7	0.0351205691920971\\
10	0.034253141847552\\
};
\addplot [color=mycolor5, line width=2.5pt, forget plot]
  table[row sep=crcr]{%
0	0.00133868808567605\\
0.0500000000000007	0.133482283244408\\
0.0999999999999996	0.144303539997882\\
0.15	0.148763701195403\\
0.199999999999999	0.152289387245041\\
0.25	0.154729968994058\\
0.300000000000001	0.155680653930851\\
0.35	0.155511268483677\\
0.4	0.154784853373243\\
0.5	0.152937262663198\\
0.550000000000001	0.151762217410523\\
0.6	0.150288752644604\\
0.65	0.14850314149975\\
0.75	0.144387174202336\\
0.949999999999999	0.136129991860955\\
1	0.134081248978275\\
1.05	0.132115211419148\\
1.1	0.130260929165159\\
1.2	0.126796790980727\\
1.35	0.121862700941641\\
1.8	0.107408989785249\\
1.9	0.104383028921946\\
2	0.10149950930208\\
2.1	0.098741998252196\\
2.2	0.0961176214950346\\
2.3	0.0936275393516937\\
2.4	0.091272118405838\\
2.5	0.089053006913991\\
2.6	0.0869578760676646\\
2.7	0.0849894819210544\\
2.8	0.0831377581786477\\
2.9	0.0813941918736631\\
3.05	0.0789631113434375\\
3.2	0.0767216552124914\\
3.35	0.0746366215673451\\
3.55	0.0720254769596895\\
3.75	0.0695483109057076\\
4	0.0666157685341791\\
4.25	0.0638669487603742\\
4.45	0.0618202098093743\\
4.65	0.0599280398122612\\
4.85	0.0581802536707325\\
5.05	0.0565573711004319\\
5.3	0.0546864518396752\\
5.6	0.0526094845610992\\
5.9	0.0506690492962605\\
6.25	0.0485643526764061\\
6.75	0.0457400966990242\\
7.05	0.0441974603872417\\
7.4	0.0425753378236386\\
7.7	0.0413067265816824\\
8.15	0.0395841579591281\\
8.65	0.037771280342902\\
9.1	0.0362783620288738\\
9.6	0.0347617756116811\\
10	0.0336387515512904\\
};
\end{axis}

\begin{axis}[%
width=0.613\linewidth,
height=0.46\linewidth,
at={(-0.08\linewidth,-0.051\linewidth)},
scale only axis,
xmin=0,
xmax=1,
ymin=0,
ymax=1,
axis line style={draw=none},
ticks=none,
axis x line*=bottom,
axis y line*=left,
title style={font=\normalsize},xlabel style={font=\normalsize},ylabel style={font=\normalsize},legend style={font=\small},ticklabel style={font=\normalsize}
]
\end{axis}
\end{tikzpicture}%
 \\
    {
	\vspace*{0.5\baselineskip}
	$\begin{array}{l}
	 \begin{tikzpicture}\draw[color=mycolor1, line width=2.5pt, dotted] (0,0) -- (0.5,0); \end{tikzpicture} ~\text{\small $\int_S g d\mu^\star$ for Nash equilibrium $\mu^\star$} \\
     \begin{tikzpicture}\draw[color=mycolor2, line width=2.5pt] (0,0) -- (0.5,0); \end{tikzpicture} ~\text{\small $\int_S g d\mu(t)$ with $\mu_0$ Guassian} \\
     \begin{tikzpicture}\draw[color=mycolor3, line width=2.5pt] (0,0) -- (0.5,0); \end{tikzpicture} ~\text{\small  $\int_S g d\mu(t)$ with $\mu_0$ half-uniform} \\
     \begin{tikzpicture}\draw[color=mycolor4, line width=2.5pt] (0,0) -- (0.5,0); \end{tikzpicture} ~\text{\small  $\int_S g d\mu(t)$ with $\mu_0$ Dirac} \\
     \begin{tikzpicture}\draw[color=mycolor5, line width=2.5pt] (0,0) -- (0.5,0); \end{tikzpicture} ~\text{\small  $\int_S g d\mu(t)$ with $\mu_0$ multiple Diracs}
	\end{array}$
	}
    \caption{{\color{black}The infinite-strategy BNN dynamics do not converge to the Nash equilibrium since $\int_S g d\mu(t) \not\to \int_S g d\mu^\star$ as $t\to\infty$ for the bounded continuous function $g(s) = \max\left\{0, \frac{1}{T-s^\star}(s-s^\star)\right\}$.}}
    \label{fig: compare_int}
\end{figure}

Our theoretical results pinpoint two key underlying technical conditions being violated in this example. In particular, $F$ is not weakly continuous since $f$ is not continuous, and furthermore, there exists $\mu\in\P(S)$ such that $F(\mu) \notin C(S)$, implying that $F$ does not even have codomain $C(S)$. Such continuity conditions are key assumptions in our stability results. This breakdown of dissipativity-based stability guarantees when moving ``from finite to infinite'' demonstrates the importance in carefully identifying the technical conditions under which infinite-dimensional stability may be guaranteed, as we have done in our main results of \autoref{sec: dissipativity_theory}.

\subsection{Continuous War of Attrition}
\label{sec: continuous_war_of_attrition}

The function $f$ defining the war of attrition game in \autoref{sec: failure_finite} can be equivalently written as
\begin{equation*}
	f(s,s') = V \Theta(s-s') - s\Theta(s'-s) - s'\Theta(s-s'),
\end{equation*}
where $\Theta \colon \R \to \R$ is the step function given by
\begin{equation*}
	\Theta(x) = \begin{aligned}
		\begin{cases}
			0 & \text{if $x < 0$}, \\
			\frac{1}{2} & \text{if $x=0$}, \\
			1 & \text{if $x>0$}.
		\end{cases}
	\end{aligned}
\end{equation*}
\citet{iyer2016evolutionary} propose a smoothed variant of the war of attrition by replacing the discontinuous step function $\Theta$ by the logistic function $\Theta_\alpha \colon \R \to \R$ given by
\begin{equation*}
	\Theta_\alpha(x) = \frac{1}{1+e^{-\alpha x}},
\end{equation*}
where $\alpha > 0$ is the smoothing parameter. However, in doing so, it is unclear whether the resulting game is monotone, where the difficulty arises when analyzing the values of $\int_S \int_S (s\Theta_\alpha(s'-s)+s'\Theta_\alpha(s-s')) d\mu(s') d\nu(s)$ for various $\mu,\nu\in\P(S)$.

In this example, we propose a relaxed variant of the game in \citet{iyer2016evolutionary} in which we only modify the war of attrition to be continuous, rather than smooth. This is accomplished by noting that
\begin{equation*}
	s\Theta(s'-s)+s'\Theta(s-s') = \min\{s,s'\}
\end{equation*}
is already a continuous function of $(s,s') \in S\times S$, and therefore the only term that should be replaced in $f(s,s')$ is $V\Theta(s-s')$, as it is where the discontinuity appears. To do so, let $\tilde{\Theta} \colon \R \to \R$ be a Lipschitz continuous function such that $0 \le \tilde{\Theta}(X) \le 1$ and $\tilde{\Theta}(x) + \tilde{\Theta}(-x) = 1$ for all $x\in \R$. For example, one may use the logistic function $\tilde{\Theta} = \Theta_\alpha$, or even a piecewise linear approximation of the step function $\Theta$ given by
\begin{equation*}
	\tilde{\Theta}(x) = \begin{aligned}
		\begin{cases}
			0 & \text{if $x < x_0$}, \\
			\frac{x}{2x_0} + \frac{1}{2} & \text{if $x\in[-x_0,x_0]$}, \\
			1 & \text{if $x>x_0$}.
		\end{cases}
		\end{aligned}
\end{equation*}
Then, we consider the game given by
\begin{equation}
	\begin{aligned}
		F_\mu(s) &\coloneqq \int_S \tilde{f}(s,s') d\mu(s'), \\
		\tilde{f}(s,s') &\coloneqq V \tilde{\Theta}(s-s') - \min\{s,s'\}.
	\end{aligned}
	\label{eq: continuous_war}
\end{equation}

We refer to our variant $F$ as the ``continuous war of attrition.'' The closer $\tilde{\Theta}$ approximates the step function $\Theta$, the closer the continuous war of attrition approximates the original form of the war of attrition. We give two new results: 1) the continuous war of attrition is a monotone game, and 2) the continuous war of attrition is weakly Lyapunov stable and globally weakly attracting under the BNN and impartial pairwise comparison dynamics. It is easily verified that indeed $F(\mu) \in C(S)$ for all $\mu\in\P(S)$, that $F$ satisfies both \autoref{ass: game} and \autoref{ass: game_derivative} with the extension $\overline{F}$ being defined by $\tilde{f}$ as well, and that $\overline{F}$ is continuously Fr\'echet differentiable. We now present our results.

\begin{theorem}
	\label{thm: continuous_war-monotone}
	The continuous war of attrition game $F \colon \P(S) \to C(S)$ defined by \eqref{eq: continuous_war} is monotone.
\end{theorem}

\hl{
\begin{proof}
	In this proof, we denote the indicator function on a set $A \subseteq \R$ by $\chi_A \colon \R \to \R$, where
	\begin{equation*}
		\chi_A(t) = \begin{aligned}
		\begin{cases}
			1 & \text{if $t \in A$}, \\
			0 & \text{if $t \notin A$}.
		\end{cases}
		\end{aligned}
	\end{equation*}

	Let $\mu,\nu \in\P(S)$. It holds that
	\iftoggle{twocol}{
	\begin{align*}
		& 2 \int_S \int_S \tilde{\Theta}(s-s') d(\mu-\nu)(s') d(\mu-\nu)(s) \\
		&\qquad = \int_S \int_S \tilde{\Theta}(s-s') d(\mu-\nu)(s') d(\mu-\nu)(s) \\
		&\qquad\qquad + \int_S \int_S \tilde{\Theta}(s'-s) d(\mu-\nu)(s') d(\mu-\nu)(s) \\
		&\qquad = \int_S \int_S \Big(\tilde{\Theta}(s-s') \\
		&\qquad\qquad + \tilde{\Theta}(s'-s)\Big) d(\mu-\nu)(s') d(\mu-\nu)(s) \\
		&\qquad = \int_S \int_S d(\mu-\nu)(s') d(\mu-\nu)(s) \\
		& \qquad = ((\mu-\nu)(S))^2 \\
		& \qquad =0,
	\end{align*}
	}{
	\begin{align*}
		& 2 \int_S \int_S \tilde{\Theta}(s-s') d(\mu-\nu)(s') d(\mu-\nu)(s) \\
		&\qquad = \int_S \int_S \tilde{\Theta}(s-s') d(\mu-\nu)(s') d(\mu-\nu)(s) + \int_S \int_S \tilde{\Theta}(s'-s) d(\mu-\nu)(s') d(\mu-\nu)(s) \\
		&\qquad = \int_S \int_S \left(\tilde{\Theta}(s-s') + \tilde{\Theta}(s'-s)\right) d(\mu-\nu)(s') d(\mu-\nu)(s) \\
		&\qquad = \int_S \int_S d(\mu-\nu)(s') d(\mu-\nu)(s) \\
		& \qquad = ((\mu-\nu)(S))^2 \\
		& \qquad =0,
	\end{align*}
	}
	since $(\mu-\nu)(S) = \mu(S) - \nu(S) = 0$. Therefore,
	\begin{equation*}
		\int_S \int_S \tilde{\Theta}(s-s')d(\mu-\nu)(s')d(\mu-\nu)(s) = 0.
	\end{equation*}
	Next, we note that
	\iftoggle{twocol}{
	\begin{align*}
		& \int_S \int_S \min\{s,s'\} d\mu(s') d\nu(s) \\
		&\quad = \int_S \int_S \int_{[0,\min\{s,s'\}]} dt d\mu(s') d\nu(s) \\
		&\quad = \int_S \int_S \int_{[0,\infty)} \chi_{\{t' \in \R : t' \le \min\{s,s'\}\}}(t) dt d\mu(s') d\nu(s) \\
		&\quad = \int_S \int_S \int_{[0,\infty)} \chi_{\{t' \in \R : t'\le s\}}(t) \chi_{\{t'\in\R : t\le s'\}}(t) dt d\mu(s') d\nu(s) \\
		&\quad = \int_S \int_S \int_{[0,\infty)} \chi_{\{\tilde{s} \in S : \tilde{s} \ge t\}}(s) \chi_{\{\tilde{s} \in S : \tilde{s} \ge t\}}(s') dt d\mu(s') d\nu(s) \\
		&\quad = \int_{[0,\infty)} \int_S \chi_{\{\tilde{s} \in S : \tilde{s} \ge t\}}(s') d\mu(s')\int_S \chi_{\{\tilde{s} \in S : \tilde{s} \ge t\}}(s) d\nu(s) dt \\
		&\quad = \int_{[0,\infty)} \mu(S \cap [t,\infty)) \nu(S\cap [t,\infty)) dt.
	\end{align*}
	}{
	\begin{align*}
		\int_S \int_S \min\{s,s'\} d\mu(s') d\nu(s) &= \int_S \int_S \int_{[0,\min\{s,s'\}]} dt d\mu(s') d\nu(s) \\
		&= \int_S \int_S \int_{[0,\infty)} \chi_{\{t' \in \R : t' \le \min\{s,s'\}\}}(t) dt d\mu(s') d\nu(s) \\
		&= \int_S \int_S \int_{[0,\infty)} \chi_{\{t' \in \R : t'\le s\}}(t) \chi_{\{t'\in\R : t\le s'\}}(t) dt d\mu(s') d\nu(s) \\
		&= \int_S \int_S \int_{[0,\infty)} \chi_{\{\tilde{s} \in S : \tilde{s} \ge t\}}(s) \chi_{\{\tilde{s} \in S : \tilde{s} \ge t\}}(s') dt d\mu(s') d\nu(s) \\
		&= \int_{[0,\infty)} \int_S \chi_{\{\tilde{s} \in S : \tilde{s} \ge t\}}(s') d\mu(s')\int_S \chi_{\{\tilde{s} \in S : \tilde{s} \ge t\}}(s) d\nu(s) dt \\
		&= \int_{[0,\infty)} \mu(S \cap [t,\infty)) \nu(S\cap [t,\infty)) dt.
	\end{align*}
	}
	Therefore, we find that
	\iftoggle{twocol}{
	\begin{align*}
		& \int_S \int_S \min\{s,s'\} d(\mu-\nu)(s') d(\mu-\nu)(s) \\
		&\quad = \int_{[0,\infty)} \Big(\mu(S\cap[t,\infty))^2 \\
		&\quad\quad - 2\mu(S\cap[t,\infty))\nu(S\cap[t,\infty)) + \nu(S\cap[t,\infty))^2\Big) dt \\
		&\quad = \int_{[0,\infty)} \left(\mu(S\cap[t,\infty)) - \nu(S\cap[t,\infty))\right)^2 dt.
	\end{align*}
	}{
	\begin{align*}
		& \int_S \int_S \min\{s,s'\} d(\mu-\nu)(s') d(\mu-\nu)(s) \\
		&\qquad = \int_{[0,\infty)} \left(\mu(S\cap[t,\infty))^2 - 2\mu(S\cap[t,\infty))\nu(S\cap[t,\infty)) + \nu(S\cap[t,\infty))^2\right) dt \\
		&\qquad = \int_{[0,\infty)} \left(\mu(S\cap[t,\infty)) - \nu(S\cap[t,\infty))\right)^2 dt.
	\end{align*}
	}
	Thus, overall, it holds that
	\iftoggle{twocol}{
	\begin{align*}
		& \left<F(\mu) - F(\nu), \mu-\nu\right> \\
		&\qquad = \int_S(F_\mu(s) - F_\nu(s)) d(\mu-\nu)(s) \\
		&\qquad = \int_S \int_S \tilde{f}(s,s') d(\mu-\nu)(s') d(\mu-\nu)(s) \\
		&\qquad = V\int_S \int_S \tilde{\Theta}(s-s')d(\mu-\nu)(s')d(\mu-\nu)(s) \\
		&\qquad\qquad - \int_S \int_S \min\{s,s'\} d(\mu-\nu)(s') d(\mu-\nu)(s) \\
		&\qquad = -\int_{[0,\infty)} \left(\mu(S\cap[t,\infty)) - \nu(S\cap[t,\infty))\right)^2 dt \\
		&\qquad \le 0.
	\end{align*}
	}{
	\begin{align*}
		\left<F(\mu) - F(\nu), \mu-\nu\right> &= \int_S(F_\mu(s) - F_\nu(s)) d(\mu-\nu)(s) \\
		&= \int_S \int_S \tilde{f}(s,s') d(\mu-\nu)(s') d(\mu-\nu)(s) \\
		&= V\int_S \int_S \tilde{\Theta}(s-s')d(\mu-\nu)(s')d(\mu-\nu)(s) \\
		&\qquad - \int_S \int_S \min\{s,s'\} d(\mu-\nu)(s') d(\mu-\nu)(s) \\
		&= -\int_{[0,\infty)} \left(\mu(S\cap[t,\infty)) - \nu(S\cap[t,\infty))\right)^2 dt \\
		&\le 0.
	\end{align*}
	}
	Hence, $F$ is monotone.
\end{proof}
}

\autoref{thm: continuous_war-monotone} allows us to immediately apply our dissipativity theory to conclude that indeed the continuous war of attrition exhibits global stability on the infinite strategy set $S$, unlike the original version of the game:

\iftoggle{longform}{
\begin{corollaryE}[][end,restate,text link=]
}{
\begin{corollaryE}[][end,text link=]
}
	\label{cor: cts_war_stable}
	Consider the continuous war of attrition game $F\colon\P(S) \to C(S)$ defined by \eqref{eq: continuous_war}. If $v \colon \P(S) \times C(S) \to T\P(S)$ is the dynamics map for either the BNN dynamics of \autoref{ex: bnn} or the impartial pairwise comparison dynamics of \autoref{ex: pcd}\iftoggle{thesis}{ and if \autoref{ass: exist_edm_solution} holds}{}, then $\NE(F)$ is weakly Lyapunov stable and globally weakly attracting under the EDM \eqref{eq: edm}.
\end{corollaryE}

\iftoggle{longform}{
\begin{proofE}
	This is immediate from \autoref{cor: recovering_monotone} together with \autoref{thm: continuous_war-monotone} and the fact that $F$ satisfies all of the appropriate regularity conditions.
\end{proofE}
}{}

In \iftoggle{longform}{\autoref{fig: continuous_war-bnn-unif}}{\autoref{fig: continuous_war-bnn-gauss}}, we display a computer simulation illustrating the stability of the continuous war of attrition game \eqref{eq: continuous_war} with $T = 2$, $V = 1$, $\tilde{\Theta} = \Theta_\alpha$, and $\alpha=100$, under the BNN dynamics. \hl{The simulation is carried out using the discretization technique described in \autoref{sec: failure_finite}, which we know respects the true stability of the infinite-dimensional dynamics due to \autoref{cor: cts_war_stable}.} \iftoggle{longform}{The initial population state in \autoref{fig: continuous_war-bnn-unif} is the uniform distribution on $S=[0,2]$. We see that the distribution function values $\mu(t)([0,s])$ converge in time towards those of a distribution closely resembling $\mu^\star$, the unique Nash equilibrium of the (discontinuous) war of attrition given in \eqref{eq: war_nash}. Upon increasing $\alpha$, this limiting distribution function even more closely approximates that of $\mu^\star$. The simulation is repeated in \autoref{fig: continuous_war-bnn-gauss} using a Gaussian initial population state with mean $1$ and variance $0.1$. The same convergent behavior is observed.}{The initial population state in \autoref{fig: continuous_war-bnn-gauss} is a Gaussian distribution on $S=[0,2]$ with mean $1$ and variance $0.1$. We see that the distribution function values $\mu(t)([0,s])$ converge in time towards those of a distribution closely resembling $\mu^\star$, the unique Nash equilibrium of the (discontinuous) war of attrition given in \eqref{eq: war_nash}. Upon increasing $\alpha$, this limiting distribution function even more closely approximates that of $\mu^\star$. Repeating the simulation for initial population states given by the uniform distribution and Dirac distributions yields the same convergent behavior.}

\iftoggle{longform}{
\begin{figure}[ht]
    \centering
    \input{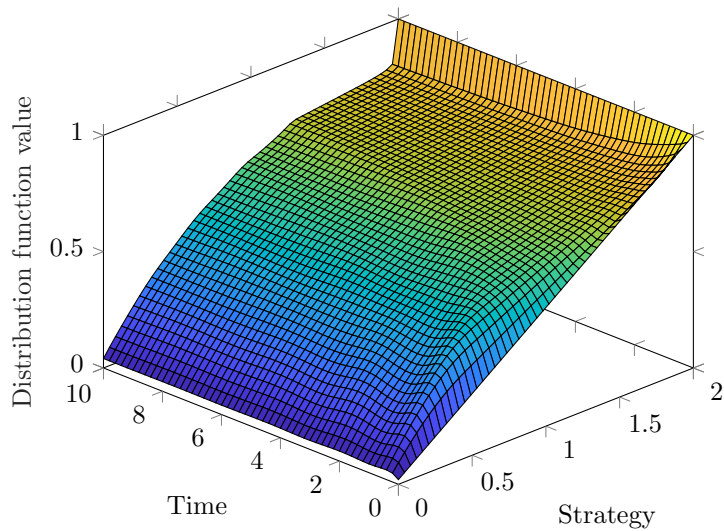}
    \caption{Evolution of the distribution function $s \mapsto \mu(t)([0,s])$ for continuous war of attrition on $S = [0,2]$ under BNN dynamics with uniform initial distribution $\mu_0$.}
    \label{fig: continuous_war-bnn-unif}
\end{figure}
}{}

\begin{figure}[ht]
    \centering
    \input{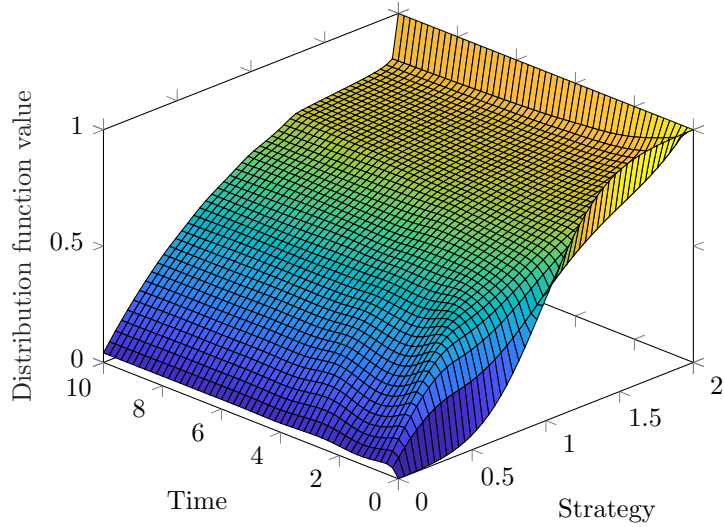}
    \caption{Evolution of the distribution function $s \mapsto \mu(t)([0,s])$ for continuous war of attrition on $S = [0,2]$ under BNN dynamics with Gaussian initial distribution $\mu_0$ (mean $1$, variance $0.1$).}
    \label{fig: continuous_war-bnn-gauss}
\end{figure}

\subsection{Smoothing Dynamics}
\label{sec: smoothing_dynamics}

\iftoggle{longform}{In this section, we consider the DPEDM \eqref{eq: dpedm} with dynamic payoffs.}{{\color{black}In this section, we consider an even more general form of evolutionary dynamics than \eqref{eq: edm}, in which the payoffs no longer take the static form $\rho(t) = F(\mu(t))$, but rather they have their own inherent dynamics. In our online technical report \citet{anderson2025dissipativity-long}, we extend our dissipativity theory to such settings.}} \iftoggle{longform}{Specifically, we}{Here, we specifically} consider smoothing dynamics, which occur when short-term variations in an evolutionary game's payoffs are suppressed, e.g., by the time delay between when a player receives payoff information and when they revise their strategy \citep{fox2013population,arcak2021dissipativity}. Formally, the smoothing dynamics\iftoggle{longform}{ DPEDM}{} corresponding to a game $F \colon \P(S) \to C(S)$ \iftoggle{longform}{is}{are} given by
\begin{equation*}
	\begin{aligned}
		\dot{\mu}(t) &= v(\mu(t),\rho(t)), \\
		\dot{\rho}(t) &= \lambda\left(F(\mu(t)) - \rho(t)\right), \\
		\mu(0) &= \mu_0, \\
		\rho(0) &= \rho_0,
	\end{aligned}
\end{equation*}
where $\lambda > 0$ is the smoothing parameter.\iftoggle{longform}{ Notice that $u(\mu,\rho) = \lambda\left(F(\mu) - \rho\right) = 0$ if and only if $\rho = F(\mu)$, so $u$ \iftoggle{longform}{is $F$-payoff stationary}{{\color{black}satisfies the rest point assumption in \autoref{thm: dpedm}}}.}{}

Even in the case of finite strategy sets, the incorporation of smoothing dynamics may turn a\iftoggle{thesis}{ dynamically}{} stable evolutionary process into an unstable one \citep{park2019population}; smoothing the payoff dynamics does not necessarily help with closed-loop stability. This may also be the case in our setting of infinite strategy sets. Indeed, for the continuous war of attrition game of \autoref{sec: continuous_war_of_attrition} with $T=2$, $V=1$, $\tilde{\Theta} = \Theta_\alpha$, and $\alpha = 100$, together with the BNN dynamics map $v$ and $\lambda = 1$, we see in \autoref{fig: cts_war-smoothing-good} that the smoothing has caused the evolution {\color{black}to lose its asymptotic stability} (\hl{the persistent oscillations are verified numerically at times $t>10^4$}).

\begin{figure}[ht]
    \centering
    \input{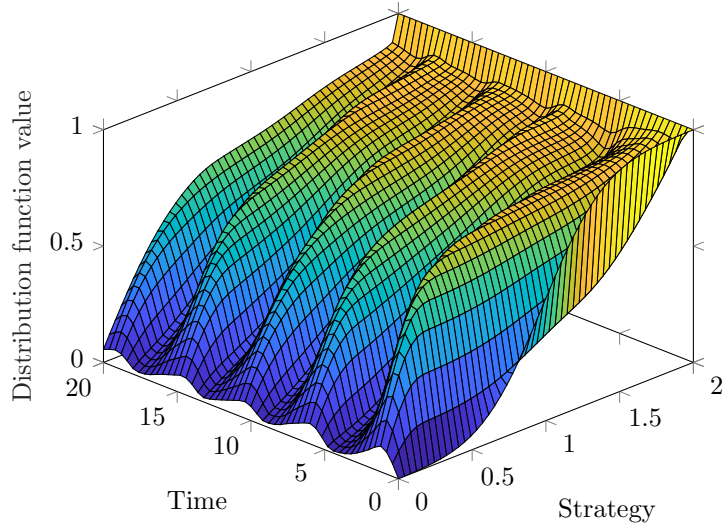}
    \caption{Evolution of the distribution function $s \mapsto \mu(t)([0,s])$ for continuous war of attrition game under BNN dynamics with smoothing, \hl{together with Gaussian initial distribution $\mu_0$ (mean $1$, variance $0.1$) and initial payoff $\rho_0 = F(\mu_0)$.}}
    \label{fig: cts_war-smoothing-good}
\end{figure}

Next, we consider the smoothing dynamics corresponding to a different game, namely, that given by
\begin{equation*}
	\begin{aligned}
		F_\mu(s) &\coloneqq \int_S f(s,s') d\mu(s'), \\
		f(s,s') &= \cos(2\pi s) - \cos(2\pi s').
	\end{aligned}
\end{equation*}
We will refer to this as the ``cosine game.'' It is easy to see that $\left<F_\nu,\nu\right> = 0$ for all $\nu\in\M(S)$, and in particular this shows that $F$ is monotone. For finite $S$, \citet{fox2013population} show that the smoothing dynamics corresponding to games satisfying $\left<F_\nu,\nu\right> \le 0$ for all $\nu\in\M(S)$ \iftoggle{longform}{are $\delta$-antipassive}{{\color{black}are what is called ``$\delta$-antipassive,''}} under an invertibility condition\iftoggle{longform}{.}{, and that this leads to closed-loop stability when coupled with $\delta$-passive dynamics maps (such as those of the BNN and impartial pairwise comparison dynamics).} \iftoggle{longform}{Therefore, one may suspect based on our \autoref{thm: dpedm} and \autoref{cor: dpedm-bnn_and_pcd} that this DPEDM with a $\delta$-passive dynamics map (such as that of BNN or impartial pairwise comparison) results in closed-loop\iftoggle{thesis}{ dynamic}{} stability.}{{\color{black}In our online technical report \citet{anderson2025dissipativity-long}, we use our generalized dissipativity theory to show that such a combination of $\delta$-passive dynamics with $\delta$-antipassive payoffs also gives rise to closed-loop stability in the infinite-strategy setting. This suggests that the cosine game smoothing dynamics under consideration, together with BNN or impartial pairwise comparison dynamics, will result in closed-loop stability.}} We numerically find that this is indeed the case for simulated dynamics with smoothing parameter $\lambda = 0.5$ and Gaussian initial population state $\mu_0$ with mean $1$ and variance $0.1$. \autoref{fig: cos-no_smoothing} shows the evolution of the population state without smoothing (i.e., for the EDM \eqref{eq: edm} with static feedback\iftoggle{longform}{}{ $\rho(t) = F(\mu(t))$}), \autoref{fig: cos-smooth-good} shows the evolution for smoothing with initial payoff $\rho_0 = F(\mu_0)$, and \autoref{fig: cos-smooth-bad} shows the evolution for smoothing with initial payoff given by $\rho_0(s) = -s^2$.

\begin{figure}[ht]
    \centering
    \input{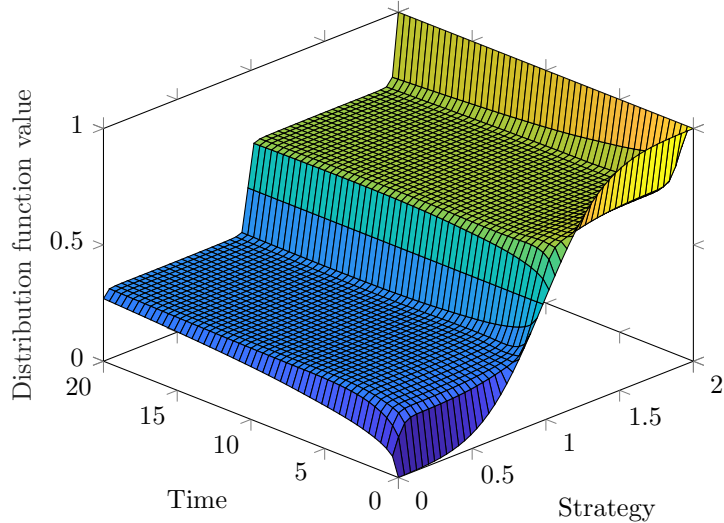}
    \caption{Evolution of the distribution function $s \mapsto \mu(t)([0,s])$ for the cosine game under BNN dynamics with static feedback, together with Gaussian initial distribution $\mu_0$ (mean $1$, variance $0.1$).}
    \label{fig: cos-no_smoothing}
\end{figure}

\begin{figure}[ht]
    \centering
    \input{figures/cos-smooth-good.tikz}
    \caption{Evolution of the distribution function $s \mapsto \mu(t)([0,s])$ for the cosine game under BNN dynamics with smoothing, together with Gaussian initial distribution $\mu_0$ (mean $1$, variance $0.1$) and initial payoff $\rho_0 = F(\mu_0)$.}
    \label{fig: cos-smooth-good}
\end{figure}

\begin{figure}[ht]
    \centering
    \input{figures/cos-smooth-bad.tikz}
    \caption{Evolution of the distribution function $s \mapsto \mu(t)([0,s])$ for the cosine game under BNN dynamics with smoothing, together with Gaussian initial distribution $\mu_0$ (mean $1$, variance $0.1$) and initial payoff $\rho_0(s) = -s^2$.}
    \label{fig: cos-smooth-bad}
\end{figure}

All evolutions appear to exhibit asymptotic stability towards a Nash equilibrium; it is easy to verify that $\delta_0$, $\delta_1$, and $\delta_2$ are all Nash equilibria of $F$, and hence the convex combination $\frac{1}{3}\delta_0 + \frac{1}{3}\delta_1 + \frac{1}{3}\delta_2$ is as well. Interestingly, in the case of \autoref{fig: cos-smooth-bad} where the initial payoff is uninformative of the game's structure, the population state initially approaches a different Nash equilibrium, namely $\delta_0$, before the system overcomes the time delay of smoothing and begins approaching $\frac{1}{3}\delta_0+\frac{1}{3}\delta_1+\frac{1}{3}\delta_2$. However, running the simulations for a longer time horizon shows that all of these evolutions actually end up adjusting their mass distributions to coincide with an even different Nash equilibrium, that being $\frac{1}{2}\delta_0 + \frac{1}{2}\delta_2$ (cf., \autoref{fig: cos-no_smoothing-long} for the static feedback case).

\begin{figure}[ht]
    \centering
    \input{figures/cos-no_smoothing-long.tikz}
    \caption{Long-time evolution of the distribution function $s \mapsto \mu(t)([0,s])$ for the cosine game under BNN dynamics with static feedback, together with Gaussian initial distribution $\mu_0$ (mean $1$, variance $0.1$).}
    \label{fig: cos-no_smoothing-long}
\end{figure}

	\section{Conclusions}
\label{sec: conclusions}

In this \iftoggle{thesis}{thesis}{paper}, we extend notions from dissipativity theory to evolutionary games with an infinite number of strategies. Our novel\iftoggle{thesis}{ dynamic}{} stability results for games evolving under $\delta$-dissipative evolutionary dynamics provide a complete characterization of the technical conditions under which such stability is guaranteed. \iftoggle{longform}{We both specialize our theory to monotone games, and extend our theory to $\delta$-dissipative evolutionary dynamics coupled with $\delta$-antidissipative dynamic feedback payoffs.}{} \hl{Our new framework and results are applicable to much broader classes of games and dynamics than past works, recovering a handful of prior stability guarantees as special cases. This breadth is illustrated through case studies including a newly proposed variant of the classical war of attrition game.} Interesting directions for future research include the development of sufficient conditions for $\delta$-dissipativity\iftoggle{longform}{ and $\delta$-antidissipativity}{} from properties of a system's finite-strategy approximations, the identification and analysis of game-theoretic models and applications falling within the scope of our framework, and extensions to games with multiple populations. \hl{Another open problem of interest is the generalization of the invertibility requirement used in \citet[Theorem~4.6]{fox2013population} to our \iftoggle{longform}{setting of maps between Banach spaces in order to prove $\delta$-antipassivity of payoff maps generated by smoothing of monotone games.}{{\color{black}Banach space setting in order to prove stability for the smoothing dynamics of monotone infinite-strategy games.}}} {\color{black}Finally, rigorously relating {\color{black}our dissipativity theory to the stability of discrete-time continuous-strategy evolutionary games poses an interesting future research direction.}}

	\appendix

	\iftoggle{proofs}{
	
	\section{Proofs}
	\label{sec: proofs}
	\printProofs

		\section{Supplementary Definitions and Results}
	\label{sec: supplementary}

	\subsection{Differentiation in Banach Spaces}

	Here, we review formal definitions for the notions of differentiability used throughout this \iftoggle{thesis}{thesis}{paper}.

	\begin{definition}
		\label{def: derivative}
		Let $(X,\|\cdot\|)$ be a Banach space. A mapping $x \colon [0,\infty) \to X$ is \emph{differentiable at $t=0$} if there exists $\dot{x}(0)\in X$ such that
		\begin{equation*}
			\lim_{\epsilon \downarrow 0}\left\lVert \frac{x(\epsilon) - x(0)}{\epsilon} - \dot{x}(0) \right\rVert = 0,
		\end{equation*}
		and is \emph{differentiable at $t\in(0,\infty)$} if there exists $\dot{x}(t) \in X$ such that
		\begin{equation*}
			\lim_{\epsilon \to 0}\left\lVert \frac{x(t+\epsilon) - x(t)}{\epsilon} - \dot{x}(t) \right\rVert = 0,
		\end{equation*}
		and in either of these cases, $\dot{x}(t)$ is called the \emph{derivative of $x$ at $t$}. A mapping $x \colon [0,\infty) \to X$ that is differentiable at $t=0$ and at every $t\in(0,\infty)$ is called \emph{differentiable}.
	\end{definition}

	\begin{definition}
		\label{def: strong_derivative}
		A mapping $\mu \colon [0,\infty) \to \M(S)$ is \emph{strongly differentiable at $t\in[0,\infty)$} if $\mu$ is differentiable at $t$ with respect to the norm $\|\cdot\|_\TV$ on the Banach space $\M(S)$.
	\end{definition}

	A strong derivative $\dot{\mu}(t)$ of $\mu$ at $t$, if it exists, is necessarily unique. The qualifier ``strong'' is used to emphasize that $\dot{\mu}(t)$ is defined in terms of convergence with respect to the strong topology. Note that if $\mu$ is strongly differentiable, then it is continuous with respect to the strong topology. In this case, since every weakly open set is strongly open, it must also be that $\mu$ is weakly continuous.

	\begin{definition}
		\label{def: frechet}
		Let $(X,\|\cdot\|_X)$ and $(Y,\|\cdot\|_Y)$ be Banach spaces and let $U \subseteq X$ be open. A mapping $f\colon U\to Y$ is called \emph{Fr\'echet differentiable at $x\in U$} if there exists a bounded linear operator $Df(x) \colon X \to Y$ such that
		\begin{equation*}
			\lim_{\epsilon \to 0} \frac{\|f(x+\epsilon) - f(x) - Df(x)\epsilon\|_Y}{\|\epsilon\|_X} = 0,
		\end{equation*}
		and in this case $Df(x)$ is called the \emph{Fr\'echet derivative of $f$ at $x$}. A mapping $f\colon U\to Y$ that is Fr\'echet differentiable at every $x\in U$ is called \emph{Fr\'echet differentiable}.
	\end{definition}

	Throughout this work, we consider maps $f\colon U\to Y$ with $U\subseteq X$ where $(X,\|\cdot\|_X)$ and $(Y,\|\cdot\|_Y)$ may be $(\R,|\cdot|)$, $(\M(S),\|\cdot\|_\TV)$, or $(C(S),\|\cdot\|_\infty)$. Fr\'echet differentiability is always with respect to one of the norms $|\cdot|$, $\|\cdot\|_\TV$, or $\|\cdot\|_\infty$ in this work. The particular norm is clear from context. We remark that $\mu\colon [0,\infty) \to \M(S)$ is strongly differentiable on $(0,\infty)$ if and only if it is Fr\'echet differentiable on $(0,\infty)$. In this case, the strong derivative coincides with the Fr\'echet derivative \hl{under the identification of} $\dot{\mu}(t)$ with the linear map $D\mu(t) \colon \R \to \M(S)$ defined by usual multiplication; $D\mu(t) \colon \epsilon \mapsto \dot{\mu}(t) \epsilon$. We similarly identify $D\mu(0)$ with $\dot{\mu}(0)$ when $\mu$ is strongly differentiable at $0$. As is the case with strong derivatives, Fr\'echet derivatives are unique when they exist.

	Partial Fr\'echet differentiation is defined as follows:

	\begin{definition}
		\label{def: partial_frechet}
		Let $(X,\|\cdot\|_X)$, $(Y,\|\cdot\|_Y)$, and $(Z,\|\cdot\|_Z)$ be Banach spaces and let $U\subseteq X$ and $V\subseteq Y$ be open. Let $(x,y)\in U\times Y$ and assume that $f(\cdot,y)\colon U\to Z$ and $f(x,\cdot)\colon V\to Z$ are Fr\'echet differentiable. The \emph{first partial Fr\'echet derivative of $f$ at $(x,y)$} is the bounded linear operator $\partial_1 f(x,y) \colon X\to Z$ defined by
		\begin{equation*}
			\partial_1 f(x,y) = D(f(\cdot,y))(x).
		\end{equation*}
		Similarly, the \emph{second partial Fr\'echet derivative of $f$ at $(x,y)$} is the bounded linear operator $\partial_2 f(x,y) \colon Y\to Z$ defined by
		\begin{equation*}
			\partial_2 f(x,y) = D(f(x,\cdot))(y).
		\end{equation*}
	\end{definition}

	\subsection{Alternative Notions of Equilibrium in Population Games}

	Aside from the notion of a Nash equilibrium, another commonly used notion of static stability within evolutionary game theory is the following, due to \citet{smith1974theory}.

	\begin{definition}
		\label{def: ess}
		A population state $\mu\in\P(S)$ is an \emph{evolutionarily stable state (ESS) of the game $F\colon \P(S) \to C(S)$} if, for all $\nu\in\P(S)\setminus\{\mu\}$, there exists $\epsilon(\nu) \in (0,1]$ such that for all $\eta\in (0,\epsilon(\nu)]$ it holds that
		\begin{equation}
			h^F_{\nu:\mu}(\eta) \coloneqq E_F(\nu,(1-\eta)\mu+\eta \nu) - E_F(\mu,(1-\eta)\mu+\eta \nu) < 0.
			\label{eq: ess}
		\end{equation}
		The function $h^F_{\nu:\mu}$ is called the \emph{score function of $\nu$ against $\mu$}, and the value $\epsilon(\nu)$ is called an \emph{invasion barrier for $\mu$ against $\nu$}.
	\end{definition}

	Intuitively, a population state $\mu\in\P(S)$ is evolutionarily stable whenever the average mean payoff to a mutated population $\nu$ is lower given payoffs defined by a small mutation $(1-\eta)\mu+\eta \nu$ towards it, i.e., the population is not incentivized to continue evolving towards any mutant population given a small fluctuation towards it. Perhaps less commonly used is the following relaxation of evolutionary stability---yet, it becomes important in the study of monotone games to be defined later.

	\begin{definition}
		\label{def: nss}
		A population state $\mu\in\P(S)$ is a \emph{neutrally stable state (NSS) of the game $F\colon\P(S) \to C(S)$} if, for all $\nu\in\P(S)$, there exists $\epsilon(\nu)\in(0,1]$ such that for all $\eta\in(0,\epsilon(\nu)]$ it holds that
		\begin{equation*}
			h^F_{\nu:\mu}(\eta) \le 0.
		\end{equation*}
		Such a value $\epsilon(\nu)$ is called a \emph{neutrality barrier for $\mu$ against $\nu$}.
	\end{definition}

	The following proposition shows that, under a mild condition, neutral stability (and hence evolutionary stability) is stronger than stability in the sense of Nash.

	\begin{proposition}[]
		\label{prop: nss_implies_nash}
		Let $\mu\in\P(S)$ be a NSS of the game $F\colon\P(S)\to C(S)$. If $h^F_{\nu:\mu}$ is right-continuous at $0$ for all $\nu\in\P(S)$, then $\mu$ is a Nash equilibrium of the game $F$.
	\end{proposition}

	\begin{proof}
		Let $\mu\in\P(S)$ be a NSS of the game $F\colon\P(S)\to C(S)$. Suppose that $h^F_{\nu:\mu}$ is right-continuous at $0$ for all $\nu\in\P(S)$. Let $\nu\in\P(S)$. Then, there exists $\epsilon(\nu)\in(0,1]$ such that
		\begin{equation*}
			h^F_{\nu:\mu}(\eta) = E_F(\nu,(1-\eta)\mu+\eta \nu) - E_F(\mu,(1-\eta)\mu+\eta \nu) \le 0
		\end{equation*}
		for all $\eta\in(0,\epsilon(\nu)]$. Thus, by the right-continuity of $h^F_{\nu:\mu}$, it holds that
		\begin{equation*}
			E_F(\nu,\mu)-E_F(\mu,\mu) = h^F_{\nu:\mu}(0) = \lim_{\eta \downarrow 0} h^F_{\nu:\mu}(\eta) \le 0.
		\end{equation*}
		Since $\nu$ is arbitrary, this proves the claim.
	\end{proof}

	Notice that $h^F_{\nu:\mu}$ is right-continuous at $0$ for all $\mu,\nu\in\P(S)$ whenever $F$ is weakly continuous. The converse of \autoref{prop: nss_implies_nash} is not true in general. However, it can be shown that a Nash equilibrium is an ESS (and hence a NSS) \hl{under additional conditions; see, e.g., \autoref{prop: monotone_nash} and \citet[Theorem~21]{bomze1989game}.}

	Notice that the notions of ESS and NSS are local ones. They can be extended into global notions as follows.

	\begin{definition}
		\label{def: globally_nss_ess}
		A population state $\mu\in\P(S)$ is a \emph{globally neutrally stable state (GNSS) of the game $F\colon\P(S)\to C(S)$} if
		\begin{equation}
			E_F(\nu,\nu) \le E_F(\mu,\nu)
			\label{eq: gnss}
		\end{equation}
		for all $\nu\in\P(S)$. If, additionally, the inequality \eqref{eq: gnss} holds strictly for all $\nu\in\P(S)\setminus\{\mu\}$, then $\mu$ is a \emph{globally evolutionarily stable state (GESS) of the game $F$}.
	\end{definition}

	As one should expect, every GNSS is a NSS, and every GESS is an ESS, as the following result shows.

	\begin{proposition}
		\label{prop: gnss_implies_nss}
		Let $\mu\in\P(S)$. If $\mu$ is a GNSS of the game $F\colon\P(S)\to C(S)$, then it is a NSS of the game $F$. If $\mu$ is a GESS of the game $F$, then it is an ESS of the game $F$.
	\end{proposition}

	\begin{proof}
		Suppose that $\mu\in\P(S)$ is a GNSS of the game $F\colon\P(S)\to C(S)$. Let $\nu\in\P(S)$. Then, since $\mu$ is a GNSS of the game $F$, it holds that
		\begin{equation*}
			E_F((1-\eta)\mu+\eta \nu,(1-\eta)\mu + \eta \nu) - E_F(\mu,(1-\eta)\mu+\eta \nu) \le 0
		\end{equation*}
		for all $\eta\in (0,1]$. By linearity of $E_F$ in its first argument, we find that
		\begin{equation*}
			\eta E_F(\nu,(1-\eta) \mu + \eta \nu) - \eta E_F(\mu,(1-\eta) \mu + \eta \nu) \le 0
		\end{equation*}
		for all $\eta\in(0,1]$. Dividing by $\eta$ proves that $\mu$ is a NSS of the game $F$. The proof that $\mu$ being a GESS implies that $\mu$ is an ESS is identical as above with strict inqualities when considering $\nu\in\P(S)\setminus\{\mu\}$.
	\end{proof}

	If a GESS exists, it must necessarily be the unique Nash equilibrium under a mild regularity condition, as the following proposition shows. Hence, globally evolutionarily stable states are stable in a very strong sense.

	\begin{proposition}
		\label{prop: unique_gess}
		Let $\mu\in\P(S)$ be a GESS of the game $F\colon\P(S)\to C(S)$, and suppose that $h^F_{\nu:\mu}$ is right-continuous at $0$ for all $\nu\in\P(S)$. Then, it holds that $\NE(F) = \{\mu\}$.
	\end{proposition}

	\begin{proof}
		Since $\mu$ is a GESS of the game $F$, it holds that $\mu$ is a NSS of the game $F$, and therefore $\mu\in\NE(F)$ by \autoref{prop: nss_implies_nash}, as $h^F_{\nu:\mu}$ is right-continuous at $0$. For all $\nu\in\P(S)\setminus\{\mu\}$, it holds that $E_F(\nu,\nu) < E(\mu,\nu)$ since $\mu$ is a GESS of the game $F$, and therefore such $\nu$ are not Nash equilibria of the game $F$. This proves that indeed $\NE(F) = \{\mu\}$.
	\end{proof}

	\subsubsection{Equilibria of Monotone Games}

	The following results show that the added structure of monotone games yields more information about the game's equilibria.

	\begin{proposition}
		\label{prop: monotone_nash}
		Suppose that the game $F\colon \P(S)\to C(S)$ is monotone. Then the following all hold:
		\begin{enumerate}
			\item Every Nash equilibrium of the game $F$ is a GNSS of the game $F$.
			\item Every strict Nash equilibrium of the game $F$ is a GESS of the game $F$.
			\item If $F$ is strictly monotone, then every Nash equilibrium of the game $F$ is a GESS of the game $F$.
		\end{enumerate}
	\end{proposition}

	\begin{proof}
	Let $\mu\in\P(S)$ be a Nash equilibrium of the game $F$. Then $\int_S F_\mu d\nu \le \int_S F_\mu d\mu$ for all $\nu\in\P(S)$, so by monotonicity it holds that
		\begin{align*}
			E_F(\nu,\nu) - E_F(\mu,\nu) &= \int_S F_\nu d\nu - \int_S F_\nu d\mu \\
					&= \int_S F_\mu(\nu-\mu) + \int_S(F_\nu-F_\mu)d(\nu-\mu) \\
					&\le 0
		\end{align*}
		for all $\nu\in\P(S)$. Hence, $\mu$ is a GNSS of the game $F$. It is clear that if $\mu$ is a strict Nash equilibrium or if $F$ is strictly monotone, then the above inequality becomes strict for $\nu\in\P(S)\setminus\{\mu\}$ and hence $\mu$ is a GESS of the game $F$ in these cases.
	\end{proof}

	\autoref{prop: monotone_nash} shows that we can ensure a sort of ``global evolutionary stability'' for Nash equilibria in the case of monotone games, whereas in more general games Nash equilibria may only be ``locally'' neutrally or evolutionarily stable, or they may not be neutrally or evolutionarily stable at all.
	
	\begin{corollary}
		\label{cor: unique_nash}
		Suppose that the game $F \colon \P(S) \to C(S)$ is monotone, let $\mu \in \NE(F)$, and assume that $h^F_{\nu : \mu}$ is right-continuous at $0$ for all $\nu\in\P(S)$. If either $\mu$ is a strict Nash equilibrium of $F$ or $F$ is strictly monotone, then $\mu$ is the unique Nash equilibrium of the game $F$.
	\end{corollary}

	\begin{proof}
		This follows directly from \autoref{prop: monotone_nash} together with \autoref{prop: unique_gess}.
	\end{proof}

	\begin{lemma}
		\label{lem: convex_stable_set}
		Consider a game $F\colon\P(S)\to C(S)$ and let $N\subseteq\P(S)$ be an arbitrary set of population states. Let $\stab_N^F \subseteq \P(S)$ denote the set of all population states $\mu\in\P(S)$ such that, for all $\nu\in N$, it holds that
		\begin{equation*}
			E_F(\nu,\nu)\le E_F(\mu,\nu).
		\end{equation*}
		Then, it holds that $\stab_N^F$ is a convex set.
	\end{lemma}

	\begin{proof}
		It holds that
		\begin{align*}
			\stab_N^F &= \{\mu\in\P(S) : \text{$E_F(\nu,\nu) \le E_F(\mu,\nu)$ for all $\nu\in N$}\} \\
			&= \bigcap_{\nu\in N} \{\mu\in\P(S) : E_F(\nu,\nu) \le E_F(\mu,\nu)\}.
		\end{align*}
		Since $E_F$ is linear in its first argument, the set $\{\mu\in\P(S) : E_F(\nu,\nu) \le E_F(\mu,\nu)\}$ is convex for all $\nu\in N$, and therefore the set $\stab_N^F$, being the intersection of convex sets, is also a convex set.
	\end{proof}

	We now show in \autoref{prop: monotone_convex_ne} that the set of Nash equilibria of a monotone game is a convex set under a mild regularity condition. The convexity of $\NE(F)$ rules out the case of isolated Nash equilibria. This result is similar to \citet[Lemma~2]{hofbauer2009brown}, but allows for general nonlinear maps $F$ (whereas their result is derived in the special case that $F(\mu)(s) = \int_S f(s,s') d\mu(s')$ for some function $f\colon S\times S \to \R$).

	\begin{proposition}
		\label{prop: monotone_convex_ne}
		Suppose that the game $F\colon \P(S)\to C(S)$ is monotone. If $h^F_{\nu:\mu}$ is right-continuous at $0$ for every GNSS $\mu\in\P(S)$ of the game $F$ and for all $\nu\in\P(S)$, then $\NE(F)$ is a convex set.
	\end{proposition}

	\begin{proof}
		By \autoref{prop: monotone_nash}, every Nash equilibrium of the game $F$ is a GNSS of the game $F$, and by \autoref{prop: nss_implies_nash} every GNSS of the game $F$ is a Nash equilibrium of the game $F$. Hence, the set of Nash equilibria of the game $F$ equals the set of globally neutrally stable states of the game $F$, so $\NE(F) = \{\mu\in\P(S) : \text{$E_F(\nu,\nu) \le E_F(\mu,\nu)$ for all $\nu\in\P(S)$}\}$. Applying \autoref{lem: convex_stable_set} with $N = \P(S)$ proves the claim.
	\end{proof}

	\subsection{Characteristics and Existence of Solutions to Evolutionary Dynamics}

	Since the population states of our evolutionary game are probability measures, we are primarily concerned with the case where the image of the mapping $\mu \colon [0,\infty) \to \M(S)$ is a subset of $\P(S)$ (so that the curve $t\mapsto \mu(t)$ evolves on the manifold of probability measures). In fact, for such maps, we can characterize their strong derivatives using the tangent space $T\P(S)$.

%

	\begin{proposition}
		\label{prop: tangent_space}
		Let $\mu\colon [0,\infty) \to \M(S)$ be strongly differentiable. If $\mu([0,\infty)) \subseteq \P(S)$, then $\dot{\mu}(t) \in T\P(S)$ for all $t\in[0,\infty)$.
	\end{proposition}

	\begin{proof}
		Suppose that $\mu([0,\infty))\subseteq\P(S)$. Let $t\in(0,\infty)$. Since $\frac{\mu(t+\epsilon) - \mu(t)}{\epsilon}$ converges strongly to $\dot{\mu}(t)$ as $\epsilon \to 0$, it also converges weakly to $\dot{\mu}(t)$ as $\epsilon\to 0$, so
		\begin{equation*}
			\lim_{\epsilon \to 0}\int_S f d \left(\frac{\mu(t+\epsilon) - \mu(t)}{\epsilon}\right) = \int_S f d\dot{\mu}(t)
		\end{equation*}
		for all $f\in C(S)$. In particular, taking $f$ to be the function that is identically $1$ on $S$ yields that
		\begin{equation*}
			\lim_{\epsilon \to 0} \frac{1}{\epsilon}\left(\mu(t+\epsilon)(S) - \mu(t)(S)\right) = \dot{\mu}(t)(S).
		\end{equation*}
		Since $\mu(t)$ and $\mu(t+\epsilon)$ are probability measures for all $\epsilon\in[-t,\infty)$, it holds that $\mu(t+\epsilon)(S) = \mu(t)(S) = 1$ for all such $\epsilon$, and hence $\frac{1}{\epsilon}\left(\mu(t+\epsilon)(S) - \mu(t)(S)\right) = 0$ for all $\epsilon\in[-t,\infty)\setminus\{0\}$. Therefore, it must be that
		\begin{equation*}
			\dot{\mu}(t)(S) = 0,
		\end{equation*}
		so indeed $\dot{\mu}(t) \in T\P(S)$. The case for $t=0$ follows similarly.
	\end{proof}

	\begin{remark}
	\label{rem: tangent_space}
	The proof of \autoref{prop: tangent_space} shows that, upon fixing an arbitrary time $t$, the condition $\dot{\mu}(t)\in T\P(S)$ still holds under a weaker hypothesis. In particular, if $\mu$ is strongly differentiable and $t\in(0,\infty)$ is such that there exists $\epsilon\in(0,t)$ such that $\mu((t-\epsilon,t+\epsilon)) \subseteq \P(S)$, then $\dot{\mu}(t)\in T\P(S)$.
	\end{remark}

	We now briefly discuss characteristics and existence of solutions to the EDM \eqref{eq: edm}. \autoref{prop: tangent_space} shows that if a solution $\mu\colon [0,\infty)\to \P(S)$ to the EDM \eqref{eq: edm} exists, then its strong derivative must satisfy $\dot{\mu}(t) = v(\mu(t),F(\mu(t))) \in T\P(S)$ for all $t\in[0,\infty)$, since the mapping's image satisfies $\mu([0,\infty))\subseteq \P(S)$. The intuition in this case is that the curve $\mu$, which remains in $\P(S)$ for all time, must necessarily have instantaneous velocity vectors that are ``tangent'' to $\P(S)$. This is analogous to the case where $S = \{1,2\}$ so that $\P(S)$ corresponds to the probability simplex $\{\mu \in\Rpl^2 : \mu_1+\mu_2 = 1\}$ in $\R^2$---in this setting it is geometrically obvious that a curve $\mu \colon [0,\infty) \to \P(S)$ must always have velocity vectors in $\{\nu\in\R^2 : \nu_1 + \nu_2 = 0\}$ that keep $\mu(t)$ on the probability simplex.

	Natural questions to ask are when a solution to the EDM \eqref{eq: edm} exists, and when such a solution is unique. These questions have simple answers in the case that the EDM is defined on the entire Banach space $\M(S)$ \citep[Corollary~3.9]{zeidler2013nonlinear}, but our restriction of solutions to the subset $\P(S)$ makes things more difficult. In the case that $S$ is finite, \citet[Theorem~4.4.1]{sandholm2010population} shows that a unique solution exists when $V_F \colon \P(S) \to \M(S)$ defined by $V_F(\mu) = v(\mu,F(\mu))$ is Lipschitz continuous and satisfies that $V_F(\mu)$ is in the tangent cone of $\P(S)$ at $\mu$ for all $\mu\in\P(S)$. However, this proof cannot be directly generalized to the case where $S$ is infinite, as it relies on the existence and uniqueness of closest point projections onto $\P(S)$ (which fails to hold due to non-uniqueness of solutions to $\inf_{\mu\in\P(S)}\|\mu - \nu\|_\TV$ for general $\nu\in\M(S)$, e.g., $\argmin_{\mu\in\P(S)}\|\mu\|_\TV = \P(S)$). Despite these difficulties, some related existence and uniqueness conditions have been proven for differential equations defined on closed subsets of Banach spaces, albeit, they are reliant on technically cumbersome conditions \citep{martin1973differential}. Since our work is focused on the development of\iftoggle{thesis}{ dynamic}{} stability conditions for general EDMs that do in fact possess solutions, we \iftoggle{thesis}{make use of the existence and uniqueness conditions granted by \autoref{ass: exist_edm_solution}}{assume existence and uniqueness of solutions} throughout this \iftoggle{thesis}{thesis}{paper}.

	\subsection{Dynamical Systems in Banach Spaces}

	\begin{definition}
		\label{def: lyapunov_function}
		Consider a Banach space $X$ and a topology $\tau$ on $X$. Let $Y\subseteq X$, let $v \colon Y \to X$, and let $P \subseteq Y$ be $\tau$-compact. A map $V\colon Y \to \Rpl$ is a \emph{global Lyapunov function for $P$ under $v$} if it extends to a $\tau$-continuous Fr\'echet differentiable map $\overline{V} \colon U \to \R$ defined on a norm-open set $U\subseteq X$ containing $Y$ that satisfies the following conditions:
		\begin{enumerate}
			\item $\overline{V}(x) = 0$ for all $x \in P$.
			\item $\overline{V}(x) > 0$ for all $x \in Y \setminus P$.
			\item $D\overline{V}(x) v(x) \le 0$ for all $x \in Y$.
		\end{enumerate}
		If, additionally, the map $x \mapsto D\overline{V}(x) v(x)$ is $\tau$-continuous and $D\overline{V}(x) v(x) < 0$ for all $x \in Y \setminus P$, then $V$ is a \emph{strict global Lyapunov function for $P$ under $v$}.
	\end{definition}

	Notice that the topology $\tau$ in \autoref{def: lyapunov_function} need not coincide with the topology induced by the norm on $X$. Indeed, our dissipativity results for static feedback $\rho(t) = F(\mu(t))$ rely on taking $X = \M(S)$ with $\tau$ being the weak topology and $Y = \P(S)$.

	\begin{lemma}
		\label{lem: lyapunov}
		Consider a Banach space $X$ and a topology $\tau$ on $X$. Let $Y\subseteq X$, let $v \colon Y \to X$, and let $P\subseteq Y$ be $\tau$-compact. If $\tau$ is weaker than the norm topology, $Y$ is $\tau$-compact, and there exists a global Lyapunov function for $P$ under $v$, then $P$ is $\tau$-Lyapunov stable under $v$.
	\end{lemma}

	\begin{proof}
		Suppose that there exists a global Lyapunov function $V\colon Y \to \Rpl$ for $P$ under $v$, and let $\overline{V} \colon U \to \R$ be an appropriate extension as in \autoref{def: lyapunov_function}. Let $Q\subseteq Y$ be relatively $\tau$-open and contain $P$. Then $Q = Y\cap O$ for some $\tau$-open set $O\subseteq X$. Define $\partial_{Y} Q \coloneqq Y \cap \partial O$, where $\partial O$ is the boundary of $O$ in $X$ with respect to $\tau$. It holds that $\partial_Y Q$ is $\tau$-compact since $Y$ is $\tau$-compact and $\partial O$ is $\tau$-closed. Therefore,
		\begin{equation*}
			m \coloneqq \min_{x \in \partial_{Y} Q} \overline{V}(x)
		\end{equation*}
		exists, since $\overline{V}$ is $\tau$-continuous. Notice that, since $\partial O \cap O = \emptyset$, it must be that $\partial_{Y} Q \cap Q = \emptyset$, and therefore $\partial_{Y} Q \cap P = \emptyset$. Hence, since $\overline{V}(x) > 0$ for all $x \in Y \setminus P$, it must be that $\overline{V}(x) > 0$ for all $x \in\partial_{Y} Q$ and thus $m > 0$.

		Now, let
		\begin{equation*}
			R = \{x \in Q : \overline{V}(x) \in (-\infty,m)\}.
		\end{equation*}
		Since $\overline{V}$ is $\tau$-continuous and $(-\infty,m)$ is open, the preimage $\overline{V}^{-1}((-\infty,m))$ is $\tau$-open, and hence $R = Q\cap \overline{V}^{-1}((-\infty,m)) = Y \cap O \cap \overline{V}^{-1}((0,m)) \subseteq Y$ is relatively $\tau$-open. Furthermore, since $P\subseteq Q$ and $P\subseteq V^{-1}((-\infty,m))$ as $\overline{V}(x) = 0$ for all $x\in P$, it holds that $P\subseteq R \subseteq Y$. Let $x \colon [0,\infty) \to Y$ be a solution to the differential equation $\dot{x}(t) = v(x(t))$ with $x(0) = x_0 \in Y$. Suppose that $x_0\in R$. Then, since the Fr\'echet derivative of real-valued functions on $\R$ recovers the usual derivative, we have that
		\begin{equation*}
			\frac{d \overline{V}\circ x}{dt}(t) \epsilon = D(\overline{V}\circ x)(t)\epsilon = (D\overline{V}(x(t)) \circ Dx(t)) \epsilon = D\overline{V}(x(t)) (\epsilon \dot{x}(t)) = \epsilon D\overline{V}(x(t)) v(x(t))
		\end{equation*}
		for all $t \in [0,\infty)$ and all $\epsilon\in\R$, where we have used the chain rule for Fr\'echet differentiation, linearity of Fr\'echet derivatives. Hence,
		\begin{equation*}
			\frac{d \overline{V}\circ x}{dt}(t) = D\overline{V}(x(t)) v(x(t)) \le 0
		\end{equation*}
		for all $t\in[0,\infty)$. Since $\overline{V}$ is $\tau$-continuous and $x$ is $\tau$-continuous since it is necessarily norm-continuous and $\tau$ is weaker than the norm topology, we may apply the mean value theorem to find that $\overline{V}(x(t)) \le \overline{V}(x(0)) < m$ for all $t\in[0,\infty)$. Since $R\subseteq Q$, we conclude that $x(t) \in Q$ for all $t\in [0,\infty)$, so indeed $P$ is $\tau$-Lyapunov stable under $v$.
	\end{proof}

	\begin{lemma}
		\label{lem: strict_lyapunov}
		Consider a Banach space $X$ and a topology $\tau$ on $X$. Let $Y\subseteq X$, let $v \colon Y \to X$, and let $P\subseteq Y$ be $\tau$-compact. Suppose that, for every $x_0\in Y$, there exists a unique solution $x\colon[0,\infty) \to Y$ to the differential equation $\dot{x}(t) = v(x(t))$ with $x(0) = x_0$. If $\tau$ is weaker than the norm topology, $Y$ is $\tau$-compact, and there exists a strict global Lyapunov function for $P$ under $v$, then $P$ is globally $\tau$-attracting under $v$.
	\end{lemma}

	\begin{proof}
		In this proof, we denote the complement of a subset $M \subseteq X$ by $M^c$.

		Suppose that there exists a strict global Lyapunov function $V\colon Y \to \Rpl$ for $P$ under $v$, and let $\overline{V} \colon U\to \R$ be an appropriate extension as in \autoref{def: lyapunov_function}. Let $x_0\in Y$ be arbitrary. Let $Q\subseteq Y$ be relatively $\tau$-open and contain $P$, and let $x\colon [0,\infty) \to Y$ be the unique solution to the differential equation $\dot{x}(t) = v(x(t))$ with $x(0) = x_0$. It suffices to show that there exists $T\in[0,\infty)$ such that
		\begin{equation}
			x(t)\in Q ~ \text{for all $t\in[T,\infty)$}.
			\label{eq: attracting_condition}
		\end{equation}
		Since $V$ is a global Lyapunov function, \autoref{lem: lyapunov} gives that $P$ is $\tau$-Lyapunov stable under $v$, which implies that there exists a relatively $\tau$-open set $R\subseteq Y$ containing $P$ such that $x(t) \in Q$ for all $t\in [0,\infty)$ whenever $x(0)\in R$. By time-invariance of the ordinary differential equation $\dot{x}(t) = v(x(t))$ with $x(0) = x_0$ and uniqueness of its solutions, if there exists $T\in[0,\infty)$ such that $x(T)\in R$, this implies that $x(t) \in Q$ for all $t\in [T,\infty)$. Thus, to prove \eqref{eq: attracting_condition}, it suffices to prove that there exists $T\in[0,\infty)$ such that $x(T)\in R$.

		For the sake of contradiction, suppose that $x(t) \notin R$ for all $t \in [0,\infty)$. Since $R$ is relatively $\tau$-open, $R = Y \cap O$ for some $\tau$-open set $O\subseteq X$, and therefore $Y \setminus R = Y \cap (Y \cap O)^c = Y \cap (Y^c \cup O^c) = Y \cap O^c$ is $\tau$-compact since $O^c$ is $\tau$-closed and $Y$ is $\tau$-compact. Hence,
		\begin{equation*}
			m \coloneqq \max_{y \in Y \setminus R} D\overline{V}(y) v(y)
		\end{equation*}
		exists, since $y \mapsto D\overline{V}(y) v(y)$ is $\tau$-continuous. Since $Y \setminus R \subseteq Y \setminus P$, it must hold that $m<0$ as $V$ is a strict global Lyapunov function. Furthermore, since $x(t) \in Y \setminus R$ for all $t\in[0,\infty)$, it holds that
		\begin{equation*}
			\frac{d\overline{V}\circ x}{dt}(t) = D\overline{V}(x(t))v(x(t)) \le m
		\end{equation*}
		for all $t\in[0,\infty)$. Since $\overline{V}$ is $\tau$-continuous and $x$ is $\tau$-continuous since it is necessarily norm-continuous and $\tau$ is weaker than the norm topology, we may apply the mean value theorem to conclude that, for all $\tau\in (0,\infty)$, there exists $t\in (0,\tau)$ such that
		\begin{equation*}
			\frac{\overline{V}(x(\tau)) - \overline{V}(x(0))}{\tau} = \frac{d \overline{V}\circ x}{dt}(t) \le m,
		\end{equation*}
		and hence
		\begin{equation*}
			\overline{V}(x(\tau)) \le m \tau + \overline{V}(x(0))
		\end{equation*}
		for all $\tau\in (0,\infty)$. Since $m<0$, $m\tau + \overline{V}(x(0)) \to -\infty$ as $\tau \to \infty$, which implies that there exists $\tau\in (0,\infty)$ such that $\overline{V}(x(\tau)) < 0$. Since, for such $\tau$, it holds that $x(\tau)\in Y \setminus R \subseteq Y \setminus P$, this contradicts the property of the global Lyapunov function $V$ that $\overline{V}(y) > 0$ for all $y \in Y \setminus P$. Therefore, the supposition that $x(t) \notin R$ for all $t\in[0,\infty)$ is false, and we conclude that indeed there exists $T\in [0,\infty)$ such that $x(T)\in R$, which completes the proof.
	\end{proof}

	}{}

	\cleardoublepage
	\phantomsection
	\addcontentsline{toc}{section}{References}
	\bibliographystyle{unsrtnat}
	\bibliography{tex/references.bib}
	
\end{document}